\documentclass{amsart}
\usepackage{a4wide}

\usepackage{amsfonts,amssymb,amsmath,mathtools,mathrsfs,bm,stmaryrd,amsthm,dsfont, enumerate,graphicx}
\usepackage[utf8]{inputenc}
\usepackage[hidelinks]{hyperref}
\usepackage[T1]{fontenc}
\usepackage[english]{babel}
\usepackage[capitalise]{cleveref}
\usepackage[normalem]{ulem}

\usepackage[parfill]{parskip}
\usepackage{wrapfig,tikz, tikz-cd}
\usetikzlibrary{arrows}

\usepackage{comment}
\usepackage{hyperref}

\newtheorem{Thm}{Theorem}[section]
\newtheorem{theorema}{Theorem}

\newtheorem{Lem}[Thm]{Lemma}
\newtheorem{Prop}[Thm]{Proposition}
\newtheorem{Cor}[Thm]{Corollary}
\theoremstyle{definition}
\newtheorem{Def}[Thm]{Definition}
\newtheorem{Notation}{Notation}
\newtheorem{Rem}[Thm]{Remark}
\newtheorem{Ex1}[Thm]{Example}

\newtheorem{Assumption}{Assumption}
\newtheorem{Ec}{Exceptional case}

\AtBeginEnvironment{Assumption}{\phantomsection}
\AtBeginEnvironment{Ec}{\phantomsection}

\newcommand{\asref}[1]{\hyperref[#1]{Assumption}}
\newcommand{\ecref}[1]{\hyperref[#1]{Exceptional case}}
\DeclareMathOperator{\HH}{HH}

\DeclareMathOperator{\Der}{Der}

\renewcommand{\Im}{\mathrm{Im}}
\newcommand{\Ker}{\mathrm{Ker}}
\newcommand{\chr}{\mathrm{char}}

\newcommand{\Tip}{\mathrm{Tip}}
\newcommand{\NTip}{\mathrm{NonTip}}
\newcommand{\Supp}{\mathrm{Supp}}

\newcommand{\B}{\mathcal{B}}

\renewcommand{\sl}{\mathfrak{sl}}

\makeatletter
\newcommand{\pll}{\mathord{%
		\mathchoice
		{\pll@@{0.18em}{-0.7em}}   % display
		{\pll@@{0.18em}{-0.7em}}   % text
		{\pll@@{0.12em}{-0.5em}}   % script
		{\pll@@{0.1em}{-0.4em}}    % scriptscript
}}

\newcommand{\pll@@}[2]{%
	\kern#1/\kern#2/\kern#1
}
\makeatother

\renewcommand{\parallel}{\mathrel{/\mkern-5mu/}}
\makeatletter
\newcommand{\notparallel}{%
	\mathrel{\mathpalette\not@parallel\relax}%
}
\newcommand{\not@parallel}[2]{%
	\ooalign{\reflectbox{$\m@th#1\smallsetminus$}\cr\hfil$\m@th#1\parallel$\cr}%
}

\newcommand{\gl}{\mathfrak{gl}}
\newcommand{\pgl}{\mathfrak{pgl}}
\def\e{\mathfrak{e}_1^n}

\def\ap{\mathfrak{ap}_1^n}

\newcommand{\E}{\mathcal{E}_1^n}

\newcommand{\AP}{\mathcal{AP}_1^n}
\newcommand{\M}{\mathcal{M}_1^n}

\begin{document}

\title[On $\HH^1$ of quiver algebras under gluing idempotents]{On the first Hochschild cohomology of finite dimensional quiver algebras under gluing idempotents}
\author[Y. ~Liu]{Yuming Liu}
\address{School of Mathematical Sciences, Laboratory of Mathematics and Complex Systems,
Beijing Normal University, Beijing 100875, P. R. China.}
\email{ymliu@bnu.edu.cn}

\author[L. Rubio y Degrassi]{Lleonard Rubio  y Degrassi}
\address{
Dipartimento di Matematica, Università degli Studi di Padova, Via Trieste 63, 35121 Padova PD, Italy. The corresponding author.}
\email{lleonard@math.unipd.it}

\author[C. ~Wen]{Can Wen}
\address{School of Mathematical Sciences, Laboratory of Mathematics and Complex Systems,
Beijing Normal University, Beijing 100875, P. R. China.}
\email{cwen@mail.bnu.edu.cn}

\begin{abstract}
We compare the Lie algebra structures of the first Hochschild cohomology groups of a  quiver algebra $A$ and a radical embedding $B$ obtained by gluing two idempotents of $A$. Under a mild assumption, we show that the first Hochschild cohomology groups of $A$ and $B$  are either isomorphic as Lie algebras or they differ by a one-dimensional Lie ideal. In particular, in the case of stable equivalences obtained by gluing a source and a sink vertex, we prove that either the first Hochschild cohomology groups of $A$ and $B$ are isomorphic or $\HH^1(B)$ is a central extension of $\HH^1(A)$ by a one-dimensional ideal. As a consequence, we obtain a new invariant under stable equivalences induced by gluing a source and a sink. We also compare the dimensions of $\HH^1(A)$ and $\HH^1(B)$, as well as the centers of $A$ and $B$, when gluing two arbitrary idempotents.
\end {abstract}
\maketitle

\renewcommand{\thefootnote}{\alph{footnote}}
\setcounter{footnote}{-1} \footnote{\it{Mathematics Subject
Classification(2020)}: 16E40, 16G10, 18G65.}
\renewcommand{\thefootnote}{\alph{footnote}}
\setcounter{footnote}{-1} \footnote{ \it{Keywords}:  Hochschild cohomology; Gr{\"{o}}bner basis; Idempotent gluing; Quiver algebra; Stable equivalence.}

\section{Introduction}

Let $k$ be a field. Let $A, B$ be two finite dimensional $k$-algebras and let $\mathrm{rad}(A)$, $\mathrm{rad}(B)$ be the Jacobson radicals of $A$ and $B$, respectively. Let   $\phi: B \to A$ be a radical embedding, that is, an algebra monomorphism such that $\phi(\mathrm{rad}(B))=\mathrm{rad}(A)$. Radical embeddings frequently arise in the study of finite dimensional algebras and their representation theory, for example, in determining the finiteness of the finitistic dimension of algebras \cite{EHIS,X}. If $A$ is basic and $k$ is algebraically closed, then by Xi's observation in \cite[\S 3]{X} we can assume that $B$ is a subalgebra of $A$ obtained by repeatedly gluing two idempotents of $A$. Therefore, the gluing of idempotents plays a pivotal role in the study of radical embeddings.

The gluing of idempotents is also essential in the study of stable equivalences. More precisely, Martinez-Villa proves in \cite{Martinez} that the gluing of a source and a sink induces an equivalence $\underline{\mathrm{mod}}A\stackrel{\sim}{\longrightarrow} \underline{\mathrm{mod}}B$ between the stable module categories modulo projective modules, see also \cite{KL}. Conversely, let $\phi:B \to A$ be a radical embedding obtained by gluing two primitive
idempotents. If $A$ and $B$ are stably equivalent and if the Auslander–Reiten conjecture holds, then $B$ is obtained from $A$ by gluing a source and a sink \cite[Proposition 4.11]{KL}. For this reason, we are particularly interested in this type of gluings.

It is well known that Hochschild cohomology is not functorial, that is, an algebra homomorphism
$\phi: B \to A$, does not give rise to a map from
$\HH^*(A)$ to
$\HH^*(B)$ or from
$\HH^*(B)$ to
$\HH^*(A)$. This makes Hochschild cohomology difficult to compute since it is not possible in general to reduce the study of Hochschild cohomology to   smaller, and potentially easier, algebras. However, there are specific cases for which the functorial properties of Hochschild cohomology have been shown. For example, in the context of fully faithful embeddings of differential graded categories \cite{Keller}. These arise, for example, for derived equivalences \cite{Keller} or stable equivalences of Morita type  \cite{KLZ} \cite{BR3}. In particular, these results imply that the (restricted) Lie algebra structure of the first Hochschild cohomology $\HH^1(A)$ of  an algebra $A$ is an invariant under derived equivalences, and for self-injective algebras, under stable equivalences of Morita type. 

In contrast to the situation for  stable equivalences of Morita type, stable equivalences obtained by gluing idempotents are induced by bimodules that are only projective on one side \cite{KL}. Therefore, no specific invariants are known for these types of stable equivalences beyond those established for general stable equivalences, such as representation dimension \cite{Guo} and representation type \cite{Krause}. For this reason,  a natural question to ask is if $\HH^1(A)$ is an \emph{invariant under stable equivalences induced by gluing a source and a sink}.  If this is not the case, then one could  ask if $\HH^1(A)$ \emph{still has some functoriality properties}, that is, if it is possible to define a (restricted) Lie algebra homomorphism between $\HH^1(A)$ and $\HH^1(B)$. More generally, similar questions could be asked in the case of gluing of two arbitrary idempotents.

The main aim of this work is to address these questions. Let $A$ and $B$ be two finite dimensional quiver algebras  such that $B$ is obtained from $A$ by gluing two arbitrary idempotents. By \cite{RSS,LX}, we can compute  $\HH^{1}(A)$ as the quotient $\Ker(\delta^1_{A})/\Im(\delta^0_{A})$, where $\delta_{A}$ is the differential of a cochain complex $\mathcal{C}_{para}$ which can be described by the  generalized parallel paths method. A similar computation applies to $B$, therefore we can use the complex $\mathcal{C}_{para}$ to compare the Lie algebra structures of $\HH^1(A)$ with $\HH^1(B)$. To make this comparison, we also define, for a fixed gluing of two idempotents, two subspaces $V_{\E}\subseteq \Im(\delta^0_B)$ and $V_{\AP}\subseteq \Ker(\delta^1_B)$, see Definition \ref{special-path} and Definition \ref{elements-in-Z-spp} for further details. When gluing a source and a sink, or equivalently, in the case of a stable equivalence, we have that $V_{\AP}=V_{\E}$. This condition plays a pivotal role in our main theorems.

\begin{theorema}[Theorem \ref{Lie-strucutre2-of-hh1}]
	\label{thm1}
	Let $A$ be a quiver algebra and let $B$ be a radical embedding obtained by gluing two idempotents of $A$. Let $\mathrm{char}(k)$ be zero or big enough and assume $V_{\AP}=V_{\E}$. 
		\begin{enumerate}
		\item If we glue from two different blocks of $A$, then  $\HH^{1}(A)\simeq \HH^{1}(B)$ as  (restricted) Lie algebras.
		\item  If we glue from the same block of $A$, then $\HH^{1}(A)\simeq \HH^{1}(B)/\mathcal{I}$ as  (restricted) Lie algebras, where $\mathcal{I}$ is a one-dimensional (restricted) Lie ideal of $\HH^{1}(B)$.
	\end{enumerate} 
\end{theorema}

As a consequence we obtain: 

\begin{theorema}[Corollary \ref{cor:iso}]
		\label{thm1.2}
	Let $A$ be a quiver algebra and let $B$ be a radical embedding obtained by gluing two idempotents of $A$. Let $\mathrm{char}(k)$ be zero or big enough and assume $V_{\AP}=V_{\E}$. Then
	\[
	\HH^1(A)/\mathrm{rad}(\HH^1(A))\simeq \HH^1(B)/\mathrm{rad}(\HH^1(B)).
	\]
\end{theorema}

In particular, for quiver algebras, we obtain a \emph{new invariant} under stable equivalences induced by gluing a source and a sink.
Theorem \ref{Lie-strucutre-of-hh1} addresses also the case $V_{\AP}\neq V_{\E}$.  In this setting,  we show that $\mathcal{I}$ is not a Lie ideal and we give an exact commutative diagram which relates $\HH^1(A)$ and $\HH^1(B)$. More general conditions for the validity of the above theorem can be found in Assumption \ref{assump}.  In the particular case of stable equivalences induced by idempotent gluing, we obtain the following result:

\begin{theorema}[Theorem \ref{thm:centralext}, Corollary \ref{product-decomp-source-sink}]
		\label{thm2}
	Let $A=kQ_A/I_A$ be a quiver algebra and let $B=kQ_B/I_B$ be a radical embedding obtained by gluing a source vertex  and a sink vertex from the same block of $A$. Then the one-dimensional Lie ideal $\mathcal{I}$ lies in the center of $\HH^{1}(B)$ and $\HH^1(B)$ is a central extension of $\HH^1(A)$ by $\mathcal{I}$. In addition, if $\mathrm{char}(k)=0$ and if $A$ is a monomial algebra, then there is a Lie algebra isomorphism 
\[\HH^{1}(B)\simeq \HH^{1}(A) \oplus \mathcal{I}.
	\]
\end{theorema}

Let $c_A,c_B$ be the number of blocks of $A, B$, respectively. We also compare the dimensions of $\HH^1(A)$ and $\HH^1(B)$ when gluing of two arbitrary idempotents:

\begin{theorema}[Theorem \ref{hh1-quotient-isomorphism-for-monomial-algebras}]
		\label{thm3}
	Let $A=kQ_A/I_A$ be a quiver algebra and let $B=kQ_B/I_B$ be a radical embedding obtained by gluing two idempotents of $A$. If $\mathrm{char}(k)$ is zero or big enough, then we have $$\dim \HH^{1}(A)=\dim \HH^{1}(B)-1-\dim V_{\AP}+\dim V_{\E} +c_A-c_B.$$
\end{theorema}
 In \cite[Theorem 1]{CS}  the authors give a formula to compute the dimension of $\HH^1(A)$ for a monomial algebra $A$ which allows to give another interpretation for the dimension of $V_{\AP}$ for monomial algebras, see Remark \ref{CiSa} for further details.  Furthermore, in Section \ref{inverse-gluing-sec} we give an interpretation of  S\'{a}nchez-Flores' description of the first Hochschild cohomology for radical square zero algebras in terms of  gluing operations.

Finally, we study the relation between a radical embedding $\phi: B \to A$  and the centers $Z(A), Z(B)$ of $A$ and $B$, respectively.

\begin{theorema}[Proposition \ref{center-from-the-same-block}, Proposition \ref{center-from-different-block}]
		\label{thm4} Let $A$ be a quiver algebra and let $B$ be a radical embedding of
	$A$ obtained by gluing two idempotents of $A$. Then there is an algebra monomorphism
	\begin{itemize}
		\item $Z(A)\hookrightarrow Z(B)$, if we glue from the same block of $A$.
		\item $Z(B)\hookrightarrow Z(A)$,  if we glue from different blocks of $A$.
	\end{itemize}
\end{theorema}

We also provide an explicit combinatorial formula to calculate the difference of the dimensions between $Z(A)$ and $Z(B)$.

We note that the authors of this paper have obtained similar results for monomial algebras in the case of gluing arrows  \cite{LRW2}.

{\bf Outline.} In Section \ref{pre}, we introduce some notation that will be used throughout the paper and provide background on various topics. In Section \ref{HH1sec} we prove Theorem \ref{thm1}, Theorem \ref{thm1.2}, Theorem \ref{thm3} and first part of Theorem \ref{thm2}. In Section \ref{mon-sec} we prove the second part of Theorem \ref{thm2}. In Section \ref{rsz-sec} we apply our main results  to radical square zero algebras. In Section \ref{inverse-gluing-sec} we give an interpretation on S\'{a}nchez-Flores' description  of the first Hochschild cohomology for radical square zero algebras \cite{S-F} by inverse gluing operations. In Section \ref{HHsec} we prove Theorem \ref{thm4}. In Section \ref{eg}
we provide various examples to illustrate our definitions and results.

\section{Preliminaries}
\label{pre}

\medskip

\subsection{Bound quivers}\label{sec:bound-quiver}
\ 
\newline

Le $Q$ be a quiver and $n \in \mathbb{N}$. Set $Q_{n}$ to be the set of paths of length $n$ of $Q$ and let $Q_{\ge n}$ be the set of paths of length greater than or equal to  $n$. Note that $Q_{0}$ is the set of vertices and  $Q_{1}$ is the set of arrows of $Q$. The number of vertices and arrows of $Q$ is denoted by $|Q_0|$ and $|Q_1|$, respectively. We denote by $s(\gamma)$ the source vertex of an (oriented) path $\gamma$ of $Q$ and by $t(\gamma)$ its terminal vertex. The \emph{path algebra} $kQ$ is the $k$-linear span of the set of paths of $Q$, where the multiplication of $\beta \in Q_{i}$ and $\alpha \in Q_{j}$ is provided by the concatenation $\beta\alpha \in Q_{i+j}$ if $t(\alpha)=s(\beta)$ and 0 otherwise. We denote by $l(p)$ the length of a path $p$. A path $p$ of length $l\geq 1$ is an \emph{oriented cycle} (or an \emph{oriented} $l$-\emph{cycle}) if $s(p)=t(p)$. An oriented $1$-cycle is called a \emph{loop}. Two paths $\epsilon,\gamma$ of $Q$ are called \emph{parallel} if $s(\epsilon)=s(\gamma)$ and $t(\epsilon)=t(\gamma)$, denoted by $\epsilon \pll  \gamma$. If $\epsilon$ and $\gamma$ are not parallel, we denote by $\epsilon \notparallel \gamma$. If $X,Y$ are sets of paths of $Q$, we denote by $X \pll  Y$ the set of parallel paths $\epsilon \pll  \gamma$ with $\epsilon \in X$ and $\gamma \in Y$, and denote by $k(X \pll  Y)$ the $k$-vector space with basis $X \pll  Y$. An element in $kQ$ is called \emph{uniform} if it is a linear combination of parallel paths. We will always denote by $\dim V$ the dimension of a $k$-vector space $V$.

Let $Q$ be a finite quiver and let $k$ be a field. Let $J$ be the two-sided ideal of the path algebra $kQ$ generated by all the arrows of $Q$. A two-sided ideal $I$ of $kQ$  is \emph{admissible} if there exists  $m\geq 2$ such that $J^m\subseteq I\subseteq J^2$.  If $I$ is an admissible ideal of $kQ$, the pair $(Q, I)$ is said to be a \emph{bound
quiver} and the quotient algebra $kQ/I$ is called a \emph{quiver algebra} with the Gabriel quiver $Q$. Throughout this paper, all algebras considered are
finite dimensional algebras which are isomorphic to quiver algebras. Any homomorphism between two algebras sends the identity element to the identity element.

We fix a quiver algebra $A=kQ_{A}/I_{A}$. Denote the vertices of $Q_A$ by $e_1,\cdots,e_n$.  A vertex $e_i$ is \emph{isolated} if there is no arrow $\alpha$ such that $s(\alpha)=e_i$ or $t(\alpha)=e_i$. A \emph{source vertex} $e_i$ of $Q_A$ is a vertex such that there is no arrow $\alpha$ with $t(\alpha)=e_i$. A \emph{sink vertex} $e_j$ of $Q_A$ is a vertex such that there is no arrow $\alpha$ with $s(\alpha)=e_j$. By abuse of notation, we denote by $e_1,\cdots,e_n$ the corresponding primitive orthogonal idempotents in the algebra $A$. For a path $p$ in $Q_A$, we use the same notation to denote its image $\overline{p}=p+I_A$ in $A$. If $A=A_1\times \cdots \times A_s$ is a decomposition of $A$ into a product of indecomposable algebras, then $A_i$'s are called \emph{blocks} of $A$.  Note that such a decomposition of $A$ is unique and if $s=1$, then $A$ is an indecomposable algebra. We denote by $c_A$ the number of blocks of $A$ which is also equal to the number of connected components
of the Gabriel quiver $Q_A$ of $A$.

\subsection{Gr{\"{o}}bner basis theory for quiver algebras}
\ 
\newline
Let $A=kQ/I$ be a quiver algebra. We briefly recall the Gr{\"{o}}bner basis (or Gr{\"{o}}bner-Shirshov basis) theory for the ideal $I$. 
Recall that
	$\prec$ is a \emph{well-order} on the $k$-basis $Q_{\geq 0}$ of the path algebra $kQ$ if $\prec$ is a total order on the  $k$-basis $Q_{\geq 0}$  and every nonempty subset of the
	$k$-basis $Q_{\geq 0}$ has a minimal element. First, we fix an admissible well-order $\prec$ on the $k$-basis $Q_{\geq 0}$ of the path algebra $kQ$, that is, a well-order on $Q_{\geq0}$ which is compatible with multiplication. More precisely,

\begin{Def}{\rm(\cite[Section 2.2.]{G})}\label{AO}
    Let $kQ$ be a path algebra with $k$-basis $Q_{\geq 0}$. We call a well-order $\prec$ on $Q_{\geq0}$ \emph{admissible} if the following three conditions are satisfied for $p,q,r,s\in Q_{\geq0}$:
    \begin{itemize}
        \item if $p\prec q$, then $pr\prec qr$ for both $pr\neq 0$ and $qr\neq 0$;
        \item if $p\prec q$, then $sp\prec sq$ for both $sp\neq 0$ and $sq\neq 0$;
        \item if $p=qr$, then $p\succeq q$ and $p\succeq r$.
    \end{itemize}
\end{Def}

The following \emph{left length-lexicographic order} on $Q_{\geq 0}$ is an admissible well-order for any path algebra (cf.~\cite[Example~2.1]{LX}):

Order the vertices $Q_0=\{e_1,\cdots,e_n\}$ and arrows $Q_1=\{\alpha_1,\cdots,\alpha_m\}$ arbitrarily and set the vertices smaller than the arrows. Thus, for example, let $$e_1 \prec \cdots \prec e_n\prec\alpha_1\prec\cdots\prec\alpha_m.$$
 If $p$ and $q$ are paths of length at least $1$, set  $p\prec q$ if
\begin{itemize}
  \item  $l(p)<l(q)$; or
  \item  $l(p)=l(q)=r$, $p=\beta_1\cdots\beta_r$, $q=\beta_1'\cdots\beta_r'$ with $\beta_j,\beta_j'\in Q_1$, and
        at the first index $i$ such that $\beta_i\neq \beta_i'$ we have $\beta_i\prec\beta_i'$.
\end{itemize}

Unless otherwise specified, we will always use the left length-lexicographic order in the present paper. Let $r=\sum_{p\in Q_{\geq 0},\lambda_p\in k}\lambda_pp$ be a $k$-linear combination of paths and let $$\mathrm{Supp}(r)=\{ \textrm{paths } p \textrm{ in }  r \ |\ \lambda_p\neq 0 \}.$$ The \emph{tip} of $r$, denoted by $\mathrm{Tip}(r)$, is the maximal monomial appearing with nonzero coefficient in $r$. In other words, $\mathrm{Tip}(r)=p$ if $p\in \mathrm{Supp}(r)$ and $\tilde{p}\preceq p$ for all $\tilde{p}\in \mathrm{Supp}(r)$. Moreover, we write $\mathrm{CTip}(r)$ for the coefficient of the tip of $r$. For a subset $X$ of $kQ$, we denote by $\mathrm{Tip}(X)=\{\mathrm{Tip}(r)\ |\ r\in X,r\neq0\}$ and set $\mathrm{NonTip}(X):=Q_{\geq0}\setminus \mathrm{Tip}(X)$. 

Let $A=kQ/I$ be a quiver algebra. By \cite{G} there is a $k$-vector space decomposition 
$$kQ=I\oplus \mathrm{Span}_k(\mathrm{NonTip}(I)).$$ So $\B:=\mathrm{NonTip}(I)$ (modulo $I$) gives a ``monomial'' $k$-basis of the quiver algebra $A=kQ/I$. 
Let $b_1, b_2\in kQ$. Then we say that $b_1$ \emph{divides} $b_2$, and we denote $b_1\mid b_2$,  if there are elements $c, d\in kQ$ such that $b_2 = cb_1d$. If $b_1$ does not divide $b_2$ we write $b_1\nmid b_2$. We can give now the definition of a Gr{\"{o}}bner basis:

\begin{Def}{\rm(\cite[Definition 2.4]{G})}\label{GS}
   Using the above notation, we say that a subset $\mathcal{G}$ of uniform elements in $I$
is a \emph{Gr{\"{o}}bner basis} for the ideal $I$ with respect to the order $\prec$ if $$\langle\mathrm{Tip}(\mathcal{G})\rangle=\langle\mathrm{Tip}(I)\rangle,$$ that is, $\mathrm{Tip}(\mathcal{G})$ and $\mathrm{Tip}(I)$ generate the same ideal in $kQ$. 
\end{Def}

Note that in this case $I=\langle\mathcal{G}\rangle$. We will see in the next theorem that there is a criterion in \cite{G}, called the \emph{Termination Theorem}, to judge whether a set of generators of an ideal $I$ in $kQ$ is a Gr{\"{o}}bner basis. Such criterion is based on the overlap relations.

\begin{Def}{\rm(\cite[Definition 2.7]{G})}\label{OR}
    Let $kQ$ be a path algebra, $\prec$ an admissible order on $Q_{\geq0}$ and $f,g\in kQ$. Suppose $b,c\in Q_{\geq0}$, such that
    \begin{itemize}
        \item $\mathrm{Tip}(f)c=b\mathrm{Tip}(g)$,
        \item $\mathrm{Tip}(f)\nmid b$ and $\mathrm{Tip}(g)\nmid c$.
    \end{itemize}
    Then the \emph{overlap relation} of $f$ and $g$ by $b,c$ is $$o(f,g,b,c)=(C\mathrm{Tip}(f))^{-1}\cdot fc-(C\mathrm{Tip}(g))^{-1}\cdot bg.$$
\end{Def}

It is clear that $\mathrm{Tip}(o(f,g,b,c))\prec \mathrm{Tip}(f)c=b\mathrm{Tip}(g)$. We can describe now the Termination Theorem. 

\begin{Thm}{\rm(\cite[Theorem 2.3]{G})}\label{TT}
    Let $kQ$ be a path algebra, $\prec$ an admissible order on $Q_{\geq0}$ and $\mathcal{G}$ a set of uniform elements of $kQ$. Suppose for every overlap relation, we have $$o(g_1,g_2,p,q)\Rightarrow_\mathcal{G}0,$$ that is, $o(g_1,g_2,p,q)$ can be divided by $\mathrm{Tip}(\mathcal{G})$, with $g_1,g_2\in \mathcal{G}$ and $p,q\in Q_{\geq 0}$. Then $\mathcal{G}$ is a Gr{\"{o}}bner basis of the ideal $\langle \mathcal{G}\rangle$ generated by $\mathcal{G}$.
\end{Thm}

 For the definition of divisibility of $o(g_1, g_2, p, q)$ by $\mathrm{Tip}(\mathcal{G})$, see \S 2.3.2 and Definition 2.6 in \cite{G}.
 
In general, an ideal in $kQ$ has many Gr{\"{o}}bner bases, but it has only one that is reduced in the sense
of the following definition.

\begin{Def}{\rm(\cite[Definition 2.5 and Proposition 2.6]{G})}\label{RGB}
    A Gr{\"{o}}bner basis $\mathcal{G}$ for the ideal $I$ is \emph{reduced} if the following three conditions are satisfied:
    \begin{itemize}
        \item $\mathcal{G}$ is tip-reduced: $\mathrm{Tip}(g)\nmid\mathrm{Tip}(h)$, for any $g\neq h\in \mathcal{G}$;         \item $\mathcal{G}$ is monic: $\mathrm{CTip}(g)=1$, for any $g\in \mathcal{G}$; 
        \item $g-\mathrm{Tip}(g)\in \mathrm{Span}_k(\mathrm{NonTip}(I))$, for any $g\in \mathcal{G}$.
    \end{itemize}
\end{Def}

For the existence of a reduced Gröbner basis for any admissible ideal $I$ of $kQ$, we refer the reader to  \cite[Corollary 2.2]{G}. It is easy to see, under a given admissible order, that $I$ has a unique reduced Gr{\"{o}}bner basis $\mathcal{G}$, and in this case $\mathrm{Tip}(\mathcal{G})$ is a minimal generator set of $\langle \mathrm{Tip}(I)\rangle$. Unless otherwise stated, we always assume in the sequel that $\mathcal{G}$ is a reduced Gr{\"{o}}bner basis of $I$.

We also recall Lemma 3.10 in \cite{LX}, which will be useful in Section \ref{HH1sec}. Since \cite[Lemma 3.10]{LX} holds for any Gr{\"{o}}bner basis which is not necessarily reduced, we include a proof here.

\begin{Lem} {\rm(cf. \cite[Lemma 3.10]{LX})}\label{loop}
    Let $A=kQ/I$ be a quiver algebra with $\mathcal{G}$ a Gr{\"{o}}bner basis for $I$. If $\alpha$ is a loop in $Q$, then $\alpha^m\in \Tip(\mathcal{G})$ and $\alpha^{m-1}\in \NTip(I)$ for some integer $m\geq 2$.
\end{Lem}

\begin{proof} First we claim that there exists an integer $m\geq 2$ such that $\alpha^m\in \Tip(\mathcal{G})$. If this is not true, then the elements $\alpha^2,\alpha^3,\cdots$ are contained in $\NTip(I)$, that is, the elements $\alpha^2,\alpha^3,\cdots$ are not contained in $\Tip(I)$. Indeed, if some $\alpha^m$ is contained in $\Tip(I)$ ($m\geq 2$), then clearly some $\alpha^{m'}$ is contained in $\Tip(\mathcal{G})$ ($m\geq m'\geq 2$), a contradiction. So the elements $\alpha^2,\alpha^3,\cdots$ are contained in $\NTip(I)$ and hence $k$-linearly independent, which contradicts the finite-dimensionality of $A$. The claim is proved.

Now take a minimum $m\geq 2$ with $\alpha^m\in \Tip(\mathcal{G})$. Then clearly $\alpha^{m-1}\in \NTip(I)$. 
\end{proof}

\subsection{Hochschild cohomology of quiver algebras.} 
\ 
\newline 
Let $A=kQ_{A}/I_{A}$ be a quiver algebra. The \emph{Hochschild cohomology} $$\HH^{*}(A):=\mathrm{Ext}^{*}_{A^{e}}(A,A)$$
of the $k$-algebra $A$ can be computed using different projective resolutions of $A$ over its \emph{enveloping algebra} $A^{e}:=A \otimes_{k} A^{op}$. The zero-th Hochschild cohomology group $\HH^{0}(A)$ is identified with the center $Z(A)$ of the algebra $A$. In particular, $Z(A)$ is a commutative subalgebra of $A$. The first Hochschild cohomology $\HH^1(A)$ is the quotient of the space of derivations $\Der(A)$ by the space of inner derivations  $\mathrm{Inn}(A)$. It is well-known that $\Der(A)$ is a Lie algebra under the Lie bracket $[f,g]=f\circ g-g\circ f$, where $f, g\in \Der(A)$. In addition,  $\mathrm{Inn}(A)$ is a Lie ideal of $\Der(A)$, therefore $\HH^{1}$($A$) has a Lie algebra structure. If the field $k$ has positive characteristic $p$, then $\HH^1(A)$ is a \emph{restricted Lie algebra}, that is, it is a Lie algebra endowed with a map called $p$-\emph{power map} that satisfies some compatibility properties with respect to the Lie algebra structure. For further background on restricted Lie algebras see for example \cite[Chapter 2]{FS}. The $p$-power map of a derivation $f$ is  defined by composing $f$ with itself $p$ times. The inner derivations form a restricted Lie ideal of the space of derivations, therefore $\HH^1(A)$ is a restricted Lie algebra.

In order to compute the first Hochschild cohomology group, one can use the following truncated projective resolution $\mathcal{P}_{min}$ (which is minimal on the degrees $0$ and $1$) of the $A$-bimodule $A$ given by Bardzell in \cite[Proposition 2.1]{Ba} (see also Chouhy and Solotar \cite{CS}. For a proof using the algebraic Morse theory, see \cite[Lemma 3.6]{LX}.):
\begin{center}
	\begin{tikzcd}
	 A \otimes_{E} k(\mathrm{Tip}(\mathcal{G})) \otimes_{E} A \arrow[r,"d_1"] & A \otimes_{E} kQ_{1} \otimes_{E} A \arrow[r, "d_0"] & A \otimes_{E} kQ_0 \otimes_{E} A \arrow[r, "\mu"] & A \arrow[r] & 0,
	\end{tikzcd}
\end{center}
where $E\simeq kQ_0$ is the subalgebra of $A$ spanned by the vertices, which is separable, and the $A$-bimodule morphisms are given by
\begin{equation*}
\begin{split}
\mu(1 \otimes_{E} e_i \otimes_{E} 1) =& e_i, \\
d_{0}(1 \otimes_{E} \alpha \otimes_{E} 1) =& \alpha \otimes_{E}s(\alpha) \otimes_{E} 1 - 1 \otimes_{E}t(\alpha)\otimes_E \alpha \text{ and} \\
d_{1}(1 \otimes_{E} \mathrm{Tip}(g) \otimes_{E} 1) =&\sum_{p=\alpha_{n}\cdots\alpha_{1}\in \mathrm{Supp}(g)} c_g(p) \sum_{i=1}^{n} \alpha_{n}\cdots \alpha_{i+1} \otimes_{E} \alpha_{i} \otimes_{E} \alpha_{i-1} \cdots \alpha_{1}
\end{split}
\end{equation*}
for all $e_i\in Q_0, \alpha \in Q_{1}$ and $g=\sum_{p\in\Supp(g)}c_g(p)p\in \mathcal{G}$ with $c_g(p)\in k$ and $p=\alpha_n\cdots\alpha_1$ for all $\alpha_i\in Q_1$.
Applying the contravariant functor $\mathrm{Hom}_{A^{e}}(-,A)$ to $\mathcal{P}_{min}$ we obtain the following cochain complex $\mathcal{C}_{min}$ (cf. \cite[Section 2]{STR} in the monomial case):
\begin{center}
		\begin{tikzcd}
		0 \arrow[r] & \mathrm{Hom}_{E^e}(kQ_{0},A) \arrow[r, "{d_{0}}^*"] & \mathrm{Hom}_{E^e}(kQ_{1},A) \arrow[r, "{d_{1}}^*"] & \mathrm{Hom}_{E^e}(k(\Tip(\mathcal{G})),A),
		\end{tikzcd}	
	\end{center}
where the differentials are given by
	\begin{equation*}
	\begin{split}
	({d_{0}}^*f)(\alpha)&=\alpha f(s(\alpha))-f(t(\alpha))\alpha,  \quad \\
	({d_{1}}^*h)(\Tip(g))&=\sum_{p=\alpha_{n}\cdots\alpha_{1}\in \mathrm{Supp}(g)} c_g(p) \sum_{i=1}^{n} \alpha_{n}\cdots \alpha_{i+1} h(\alpha_{i}) \alpha_{i-1} \cdots \alpha_{1},
	\end{split}
	\end{equation*}
where $f\in \mathrm{Hom}_{E^e}(kQ_{0},A)$, $\alpha\in Q_1$, $h\in \mathrm{Hom}_{E^e}(kQ_{1},A)$ and $g=\sum_{p\in\Supp(g)}c_g(p)p\in \mathcal{G}$ with $c_g(p)\in k$ and $p=\alpha_n\cdots\alpha_1$ for all $\alpha_i\in Q_1$.
In particular,   we have $\HH^{1}
(A)\simeq\Ker
({d_{1}}^*)/\Im({d_{0}}^*)$ as 
$k$-vector spaces. Similar to \cite[Proposition 2.8, Corollary 2.9]{STR}, we have that
$\Ker$(${d_{1}}^*$) is
isomorphic to the space  $E^e$-derivations of $A$ 
 and  
$\Im$(${d_{0}}^*$) is 
isomorphic to the space 
of the inner $E^e$-derivations of $A$.

By carrying out the identification $k(X\pll Y)\simeq \mathrm{Hom}_{E^{e}}(kX,kY)$ (which is given by sending $\epsilon \pll  \gamma$ to the map: $\epsilon\mapsto\gamma$ and $\epsilon'\mapsto 0$ for $\epsilon'\neq \epsilon$) in \cite[Lemma 2.3]{STR}, where $X$ and $Y$ are two finite subsets of paths of $Q_A$,
we can rewrite the above cochain complex which gives a more practical way of computing $\HH^{1}$.

\begin{Prop}{\rm(\cite[Proposition 3.7]{LX})}\label{GPC}
Let $A=kQ_A/I_A$ be a quiver algebra. Let $\mathcal{G}$ be a reduced Gr{\"{o}}bner basis of $I_A$, and denote by $\B$ the $k$-basis of $A$ given by $\mathrm{NonTip}(I)$ (modulo $I$). By the above mentioned identifications, the cochain complex $\mathcal{C}_{min}$ is isomorphic to the following complex
\begin{center}
		$\mathcal{C}_{para}:$ \quad\quad \begin{tikzcd}
		0 \arrow[r] & k(Q_{0}\pll \B) \arrow[r, "\delta^{0}"] & k(Q_{1}\pll \B) \arrow[r, "\delta^{1}"] & k(\mathrm{Tip}(\mathcal{G})\pll \B) \arrow[r, "\delta^{2}"] & \cdots,
		\end{tikzcd}	
	\end{center}
    where the differentials are given by
	\begin{equation*}
	\begin{split}
	\delta^{0}:k(Q_{0}\pll \B) &\to k(Q_{1}\pll \B)\\
	e\pll \gamma &\mapsto \sum_{\substack{a\in Q_1e,\\ a\gamma\in \B}}a\pll a\gamma- \sum_{\substack{a\in eQ_1,\\ \gamma a\in \B}}a\pll \gamma a,\\
	\delta^{1}:k(Q_{1}\pll \B) &\to k(\mathrm{Tip}(\mathcal{G}) \pll \B)\\
	a\pll \gamma &\mapsto \sum_{\substack{r\in \mathcal{G},\\ p\in \mathrm{Supp}(r)}}c_r(p)\mathrm{Tip}(r)\pll p^{a\pll\gamma},
	\end{split}
	\end{equation*}
    where $r=\sum_{p\in \mathrm{Supp}(r)}c_r(p)p$ with $c_r(p)\in k$ and $p^{a\pll\gamma}$ denotes the sum of all paths obtained by replacing each appearance of the arrow $a$ in $p$ by the path $\gamma$. We always identify $a\gamma,\gamma a$ and $p^{a\pll\gamma}$ with their images under $\pi_A:kQ_A\to A$, respectively. In particular, we have $\HH^0(A)\simeq \Ker(\delta^0)$ and $\HH^1(A)\simeq \Ker(\delta^1)/\Im(\delta^0)$ as $k$-vector spaces. 
\end{Prop}

The isomorphism $\HH^1(A)\simeq \Ker(\delta^1)/\Im(\delta^0)$ in Proposition \ref{GPC}  is induced by the following map:  send each $f$ in $\mathrm{Hom}_{E^{e}}(kQ_{1},k\B)$ to the element $\underset{a\pll  \gamma\in Q_1\pll  \B}{\sum}\lambda_{a,\gamma}(a\pll  \gamma)$ in $k(Q_{1}\pll  \B)$, where $f(a)=\underset{\gamma\in \B}\sum \lambda_{a,\gamma}\gamma$. Moreover, the inverse of the above isomorphism is induced by sending an element $a \pll  \gamma$ in $k(Q_{1}\pll  \B)$ to $f$ in $\mathrm{Hom}_{E^{e}}(kQ_{1},k\B)$ with $f(a)=\gamma$ and $f(b)=0$ for $a\neq b\in Q_1$.
 
The method of computing $\HH^{1}$ using parallel paths was first given by Strametz for monomial algebras in \cite{STR}. In \cite[Section 2.2]{RSS} and in \cite[Section 3.2]{LX}, this was generalized to arbitrary quiver algebras and called the \emph{generalized parallel paths method} in \cite{LX}. Moreover, Theorem 3.8 in \cite{LX} shows that the second isomorphism in Proposition \ref{GPC} is an isomorphism of Lie algebras.

\begin{Thm} \label{Lie-bracket}
    The bracket $$[a\pll\gamma,b\pll\eta]=b\pll\eta^{a\pll\gamma}-a\pll\gamma^{b\pll\eta}$$ for all $a\pll\gamma,b\pll\eta\in Q_1\pll\B$ induces a Lie algebra structure on $\Ker(\delta^1)/\Im(\delta^0)$ such that $\HH^1(A)$ and $\Ker(\delta^1)/\Im(\delta^0)$ are isomorphic as Lie algebras.
\end{Thm}

For quiver algebras, it is easy to describe the $p$-power map using the chain map from  $\mathcal{C}_{min}$ to  $\mathcal{C}_{para}$ and its inverse chain map. For example, for $p=3$, the $p$-power map of $a\pll  \gamma$  is $(a \pll  \gamma^{ a \pll  \gamma})^{ a \pll  \gamma}$. We note that several  results in this paper, such as Proposition \ref{Kernel-delta-one-injective} and Corollary \ref{quotient-lie-algebra-source-and-sink}, can be readily extended from the context of `Lie' algebras to `restricted Lie' algebras.

\begin{Rem}\label{center-ker-zero} The center $Z(A)$ of $A$ is naturally isomorphic to $\Ker(\delta^{0})$. For an explicit map between $\Ker(\delta^{0})$ and $Z(A)$, see 
the proof of Lemma \ref{Image-zero-part}. 
\end{Rem}

\subsection{Gluing of two idempotents and radical embeddings.} 
\ 
\newline
Let $A=kQ_{A}/I_{A}$ be a quiver algebra. Since each radical embedding reduces to a gluing of two idempotents, from now on we are going to consider $B$ to be a radical embedding which is obtained by gluing only two idempotents of $A$. More precisely, let $e_1,\cdots ,e_n$ be a complete set of primitive orthogonal idempotents in $A$ such that each idempotent $e_i$ corresponds to the vertex $i$ in $Q_A$. Let $B$ be a subalgebra of $A$ obtained by gluing two idempotents $e_{1}$ and $e_{n}$ of $A$. In other words, $B$ is identified as a subalgebra of $A$ generated by $f_1:=e_1+e_n,f_2:=e_2,\cdots,f_{n-1}:=e_{n-1}$ and all arrows in $Q_A$. We have $\dim B=\dim A-1$. Indeed, let $J$ be the ideal of $A$ generated by all arrows in $Q_A$. Then it is easy to see that $\dim A/J=n$ and $\dim B/J=n-1$. Note that the choice of idempotents to glue is arbitrary; however, we prefer to fix the notation such that $f_1 := e_1 + e_n$.

Clearly, there is a bijection between the arrows
of $A$ and the arrows of $B$. For each arrow $\alpha$ in
$Q_{A}$, we denote the corresponding arrow in $Q_{B}$ by
$\alpha'$. We define the quiver morphism $$\varphi:Q_{A}
\to Q_{B}$$
as follows: let $\varphi(e_{i})=f_{i}$ for
$2\le i \le n-1$, let
$\varphi(e_{1})=\varphi(e_{n})=f_{1}$, and let
$\varphi(\alpha)=\alpha'$. By extending the map $\varphi:Q_{A}\to Q_{B}$, we define $\varphi_n: (Q_A)_n \to (Q_B)_n$. More precisely, let $p=a_n\cdots a_1$ be a path in $(Q_A)_n$. Then $\varphi_n(p)=p'=a'_n\cdots a'_1$.

\begin{Rem}\label{rmk:equiv}
Note that $Q_B$ can equivalently be  described as follows: first observe that a quiver  $Q$ is  a $4$-tuple $(Q_0, Q_1, s, t)$. Let $\sim$ be an equivalence relation on the set $Q_0$. Denote by $Q_0/\sim$
 the corresponding set of equivalence classes, and let $\pi : Q_0 \to Q_0/\sim$ be the quotient map. Then the
quiver obtained from $Q$ by identifying vertices according to $\sim$ is defined as the $4$-tuple
$\hat{Q} = (Q_0/\sim, Q_1,\, \pi \circ s,\, \pi \circ t)$. Using this notation, by setting $Q=Q_A$ and  taking $\sim$ to be  the equivalence relation on $Q_0$ that identifies the two vertices $e_1$ and $e_n$ on $Q_0$ we obtain $Q_B=\hat{Q} $.
\end{Rem}

Denote by $Z_{new}$ the set of all paths $\beta'\alpha'$ of length $2$ in $Q_B$ such that the corresponding arrows $\alpha$ and $\beta$ of $Q_A$ are such that $\{ t(\alpha), s(\beta)\} = \{ e_1, e_n\}.$

\begin{Lem}\label{gabriel-quiver} 
Let $A=kQ_{A}/I_{A}$ be a quiver algebra and let $B$ be a subalgebra of $A$ obtained by gluing two idempotents $e_{1}$ and $e_{n}$ of $A$. Then $B\simeq kQ_{B}/I_{B}$, where $Q_B$ is the quiver obtained from $Q_A$ by identifying the vertices $e_1$ and $e_n$, and $I_B$ is the admissible ideal of $kQ_{B}$ generated by the elements in $I_A\cup Z_{new}$. In particular, $Q_B$ is the Gabriel quiver of $B$.
\end{Lem}

\begin{proof}
    Let $B'$ be the algebra $kQ_B/I_B$. Then there is an algebra monomorphism from $B'$ to $A$ by sending $f_1$ to $e_1+e_n$, $f_i$ to $e_i$ for $2\leq i\leq n-1$ and each arrow in $Q_B$ to the same arrow in $Q_A$. It is clear that this map factors through the inclusion $B\hookrightarrow A$, which gives rise to another algebra monomorphism from $B'$ to $B$. Moreover, since $B'$ has dimension $\dim A-1$, it must be isomorphic to $B$.
\end{proof}

The following proposition shows how a Gr{\"{o}}bner basis behaves under gluing of two idempotents. 

\begin{Prop}\label{corresponding-RGB}
 Let $A=kQ_{A}/I_{A}$ be a quiver algebra and let $B$ be a subalgebra of $A$ obtained by gluing two idempotents $e_{1}$ and $e_{n}$ of $A$. Let $\mathcal{G}_A$ be a reduced Gr{\"{o}}bner basis of $I_A$ under some left length-lexicographic order on $(Q_A)_{\geq 0}$. Consider a left length-lexicographic order on $(Q_B)_{\geq 0}$ defined as follows: order the vertices $f_i$ ($1\leq i\leq n-1$) arbitrarily and let $\alpha'\prec\beta'$ if $\alpha\prec\beta$ for $\alpha,\beta\in (Q_A)_{1}$. Identify $\mathcal{G}_A$ in $Q_A$ with $\varphi(\mathcal{G}_A)$ in $Q_B$, and similarly for $\mathrm{Tip}(\mathcal{G}_A)$. 
  Then $$\mathcal{G}_B:=\mathcal{G}_A\cup Z_{new}$$ is a reduced Gr{\"{o}}bner basis of $I_B$ under the above left length-lexicographic order on $(Q_B)_{\geq 0}$. In particular, $$\mathrm{Tip}(\mathcal{G}_B)=\mathrm{Tip}(\mathcal{G}_A)\cup Z_{new}.$$
\end{Prop}

\begin{proof}
    First we show that $\mathcal{G}_B:=\mathcal{G}_A\cup Z_{new}$ is a Gr{\"{o}}bner basis of $I_B$. Since $\mathcal{G}_A$ is a reduced Gr{\"{o}}bner basis, for each $g\in \mathcal{G}_B$, we have that $\mathrm{CTip}(g)=1$. By Theorem \ref{TT} and Lemma \ref{gabriel-quiver}, it suffices to show that $$o(g_1,g_2,p,q)=g_1q-pg_2\Rightarrow_{\mathcal{G}_B}0$$ for $g_1,g_2\in \mathcal{G}_B$ and $p,q\in (Q_B)_{\geq 0}$. The proof is divided into four cases.
    
    {\it Case 1: Let $g_1,g_2\in Z_{new}$.} Suppose $g_1=a_2'a_1'$ and $g_2=b_2'b_1'$ with $a_i', b_i'\in (Q_B)_1$ for $i=1,2$. Then $\mathrm{Tip}(g_1)q=p\mathrm{Tip}(g_2)$ is equivalent to $g_1q=pg_2$. It follows that $o(g_1,g_2,p,q)=g_1q-pg_2=0$.    
    
    {\it Case 2: Let $g_1\in \mathcal{G}_A$ and $g_2\in Z_{new}$.} Then the condition $\mathrm{Tip}(g_1)q=p\mathrm{Tip}(g_2)$ implies that \begin{equation*}
        \begin{split}
            o(g_1,g_2,p,q)& =g_1q-pg_2\\ &=(g_1-\mathrm{Tip}(g_1))q-p(g_2-\mathrm{Tip}(g_2))\\ &=(g_1-\mathrm{Tip}(g_1))q\\ &=(\sum_i\lambda_ip_i)q,
        \end{split}
    \end{equation*} where $p_i\in \mathrm{NonTip}(I_A)$ and $\lambda_i\in k$. The last two equalities follow from the facts that $g_2\in Z_{new}$ whence $g_2=\mathrm{Tip}(g_2)$, and since $g_1-\mathrm{Tip}(g_1)\in \mathrm{Span}_k(\mathrm{NonTip}(I_A))$, we can assume that $g_1-\mathrm{Tip}(g_1)=\sum_i\lambda_ip_i$ for $p_i\in \mathrm{NonTip}(I_A)$ and $\lambda_i\in k$. We claim that $q\in (Q_B)_1$, that is, the length $l(q)$ of $q$ equals 1. Indeed, if $l(q)=0$, then $\mathrm{Tip}(g_1)=pg_2$ should have a preimage  in $Q_A$, which is absurd since $g_2\in Z_{new}$. Therefore, we have $l(q)\geq 1$. Moreover, $l(q)<2$, otherwise $pg_2=\mathrm{Tip}(g_1)q$ and $l(g_2)=2$ yield that $\mathrm{Tip}(g_1)\mid p$, a contradiction.

    Assume that $g_2=a_2'a_1'$ with $a_1',a_2'\in (Q_B)_1$ and $t(a_1)\neq s(a_2)$. As a consequence, we have $q=a_1'$ and $\mathrm{Tip}(g_1)=pa_2'$ since $pg_2=\mathrm{Tip}(g_1)q$. It follows that $\Tip(g_1)$ starts at $s(a_2)\in (Q_A)_0$, and the same holds for all $p_i$ since $g_1$ is uniform. Therefore each $p_ia_1'$ has a subpath in $Z_{new}$ and we have $o(g_1,g_2,p,q)=(\sum_i\lambda_ip_i)q=(\sum_i\lambda_ip_i)a_1'\Rightarrow_{Z_{new}}0$.
    
    {\it Case 3: Let $g_1\in Z_{new}$ and $g_2\in \mathcal{G}_A$.} The proof is similar to that of Case 2.
    
    {\it Case 4: Let $g_1, g_2\in \mathcal{G}_A$.} If $l(p)=0$, then $\mathrm{Tip}(g_1)q=\mathrm{Tip}(g_2)$ which yields  $\mathrm{Tip}(g_1)\mid \mathrm{Tip}(g_2)$. Since $\mathcal{G}_A$ is reduced, we have $g_1=g_2$ and $l(q)=0$. Consequently, $o(g_1,g_2,p,q)=0$. If $l(p)>0$ such that $p$ has a subpath in $Z_{new}$, then the conditions $\mathrm{Tip}(g_1)q=p\mathrm{Tip}(g_2)$ and $\mathrm{Tip}(g_1)\nmid p$  imply that $p$ is a proper subpath of $\mathrm{Tip}(g_1)$. Hence $\mathrm{Tip}(g_1)$ has a subpath in $Z_{new}$, a contradiction. Similarly, if $l(q)=0$ or $l(q)>0$ such that $q$ has a subpath in $Z_{new}$, it will lead to a contradiction. So we may assume that $l(p)>0$, $l(q)>0$ and both $p$ and $q$ do not contain a subpath in $Z_{new}$. Then, under our assumption on the admissible order on $(Q_B)_{\geq 0}$, the overlap relation  $o(g_1,g_2,p,q)$ in $kQ_B$ becomes an overlap relation in $kQ_A$. Since $o(g_1,g_2,p,q)\Rightarrow_{\mathcal{G}_A}0$, then $o(g_1,g_2,p,q)\Rightarrow_{\mathcal{G}_B}0$.
    
    This proves that $\mathcal{G}_B$ is a Gr{\"{o}}bner basis of $I_B$. Finally, it is obvious that the Gr{\"{o}}bner basis $\mathcal{G}_B$ is reduced.
\end{proof}

\section{First Hochschild cohomology}
\label{HH1sec}

In this section we assume that $A$ is a finite dimensional algebra isomorphic to a quiver algebra $kQ_A/I_A$. Let $B=kQ_B/I_B$ be a radical embedding obtained by gluing two idempotents $e_{1}$ and $e_{n}$ of $A$. Throughout this section, we will refer to this
as a radical embedding of $A$. We exclude the case in which $e_1$ or $e_n$ is an isolated vertex.  We denote the vertices of $Q_{B}$ by $f_1,\cdots,f_{n-1}$, where $f_1$ is obtained by gluing $e_1$ and $e_n$.  For the rest of this section, we always assume that $A$ and $B$ are as in Proposition \ref{corresponding-RGB} so that $I_A$ has a reduced Gr{\"{o}}bner basis $\mathcal{G}_A$ and $I_B$ has a reduced Gr{\"{o}}bner basis $\mathcal{G}_B=\mathcal{G}_A\cup Z_{new}$ under some appropriate left length-lexicographic orders. Moreover, $A$ has a `monomial' $k$-basis $\B_A$ given by $\mathrm{NonTip}(I_A)$ (modulo $I_A$) and $B$ has a `monomial' $k$-basis $\B_B$ given by $\mathrm{NonTip}(I_B)$ (modulo $I_B$).

We briefly outline the main results of this section.  Firstly, we will compare $\Im(\delta^0_A)$ and $\Im(\delta^0_B)$. Then we will study the Lie algebra structures of $\Ker(\delta^1_A)$ and $\Ker(\delta^1_B)$. Lastly, we will compare the dimensions and the Lie structures of $\HH^1(A)$ and $\HH^1(B)$.

We will use  the cochain complex \(\mathcal{C}_{\text{para}}\) from the previous section in order to understand the behaviour of the first Hochschild cohomology under idempotent gluings. We  start by considering how idempotent gluings behave with respect to parallelism of arrows and paths. Recall from Section \ref{pre} that the quiver morphism $\varphi:Q_{A} \to Q_{B}$ sends a vertex $e_{i}$ to $f_{i}$ for $2\le i \le n-1$ and $e_{1}, e_n$ to $f_{1}$. In addition, $\varphi$ sends an arrow $\alpha$ in $Q_{A}$ to an arrow $\alpha'$ in $Q_B$.

\begin{Lem} \label{parallel-arrows}
Let $B$ be a radical embedding of $A$. Let $\alpha,\beta \in (Q_{A})_{1}$. If $\alpha \pll  \beta$, then $\alpha' \pll  \beta'$.
\end{Lem}

\begin{proof}
The proof follows from the definition of gluing of two idempotents.
\end{proof}

\begin{Lem}\label{bijection-parallel-arrows}
Let $B$ be a radical embedding obtained by gluing a source and a sink of $A$. Let $\alpha,\beta \in (Q_{A})_{1}$. Then $\alpha \pll  \beta$ if and only if $\alpha' \pll  \beta'$.
\end{Lem}

\begin{proof}
The sufficiency is obvious by Lemma \ref{parallel-arrows}, it is enough to show the necessity.  Assume  we are gluing a source, say $e_{1}$, and a sink, say $e_{n}$. Write $\sim$ for the equivalence relation on $(Q_A)_0$ that identifies the two vertices $e_1$ and $e_n$, as in Remark \ref{rmk:equiv}.  If $\alpha' \pll  \beta'$ but $\alpha  \notparallel \beta$, then
$s(\alpha) \sim s(\beta)$, $t(\alpha) \sim t(\beta)$, and $s(\alpha) \neq s(\beta)$
or $t(\alpha) \neq t(\beta)$. In the first case $\{s(\alpha),s(\beta)\}=\{e_1,e_n\}$,
but that is impossible because  $e_n$ is a sink; similarly, in the
second case $\{t(\alpha),t(\beta)\}=\{e_1,e_n\}$, and again this is impossible because
$e_1$ is a source.
\end{proof}

We now partially extend the above results to parallel paths. 

\begin{Prop}\label{parallel-paths-in-monomial-algebras}
Let $A=kQ_A/I_A$ be a quiver algebra and let $B=kQ_B/I_B$ be a radical embedding. Then the following hold:

$(1)$ The map $\varphi: Q_{A} \to Q_{B}$ induces a surjective map, also denoted by $\varphi: \B_A\rightarrow \B_B$, such that $\varphi^{-1}(p')=\{p\}$ for $p'\neq f_1$ and $\varphi^{-1}(f_1)=\{e_1,e_n\}$, where we denote $\varphi(p)$ by $p'$ for $p\in \B_A$.

$(2)$ Let $p,q \in \B_A$. If $p \pll  q$ in $Q_A$, then $p' \pll  q'$ in $Q_B$.

$(3)$ The map $\varphi: \B_A\rightarrow \B_B$ induces $k$-linear maps $$\varphi_0: k((Q_A)_0\pll  \B_A) \to k((Q_B)_0\pll  \B_B),$$ $$\varphi_1: k((Q_A)_1\pll  \B_A) \to k((Q_B)_1\pll  \B_B),$$ $$\varphi_2: k(\mathrm{Tip}(\mathcal{G}_A)\pll  \B_A) \to k(\mathrm{Tip}(\mathcal{G}_B)\pll  \B_B).$$
\end{Prop}

\begin{proof} We identify $\B_A$ with $\mathrm{NonTip}(I_A)$  (modulo $I_A$) and observe that 
$$\mathrm{NonTip}(I_A):=(Q_A)_{\geq 0}\setminus \Tip(I_A)$$ consists of monomial elements. The same holds for $\B_B$.

The quiver morphism $\varphi: Q_{A} \to Q_{B}$ induces a $k$-linear map $kQ_{A} \to kQ_{B}$ between path algebras by sending a path $p=a_m\cdots a_1$ ($a_i\in (Q_A)_1$ for $1\leq i\leq m$) in $Q_A$ to a path $p':=a'_m\cdots a'_1$ in $Q_B$. Clearly, the condition $p\in \B_A$ is equivalent to $p\notin \langle \Tip(I_A) \rangle=\langle \Tip(\mathcal{G}_A) \rangle\subseteq kQ_A$. We deduce that $p'\notin \langle \Tip(I_B) \rangle=\langle \Tip(\mathcal{G}_B) \rangle=\langle \Tip(I_A)\cup Z_{new}\rangle\subseteq kQ_B$ since the elements in the set $Z_{new}$ are the newly formed relations in $I_B$. Therefore $p'\in \B_B$. Now Statement (1) follows from the fact that $\dim B=\dim A-1$, and Statement (2) follows from (1) and the definition of gluing of two idempotents, which in turn implies Statement (3).
\end{proof}

We have the following  \emph{non-commutative} diagram:

\begin{center}
	\begin{tikzcd}
	0 \arrow[r] & k((Q_A)_{0}\pll  \B_A) \arrow[d, "\varphi_{0}"] \arrow[r, "\delta_A^{0}"] & k((Q_A)_{1}\pll  \B_A) \arrow[d, "\varphi_{1}"] \arrow[r, "\delta_A^{1}"] & k(\Tip(\mathcal{G}_A)\pll  \B_A) \arrow[d, "\varphi_{2}"]\\
    0 \arrow[r] & k((Q_B)_{0}\pll  \B_B) \arrow[r, "\delta_B^{0}"] & k((Q_B)_{1}\pll  \B_B) \arrow[r, "\delta_B^{1}"] & k(\Tip(\mathcal{G}_B)\pll  \B_B)~.
	\end{tikzcd}	
	\quad \quad $(*)$
\end{center}
Note that the top and the bottom  complexes are truncations of the complexes $\mathcal{C}_{para}$ of $A$ and of $B$, respectively.
Although both squares in the diagram $(*)$ are not commutative in general, there are  close connections between the coboundaries (resp. the cocycles) of the top complex and the coboundaries (resp. the cocycles) of the bottom complex in the diagram $(*)$.

In order to compare $\Im(\delta^0_A)$ and $\Im(\delta^0_B)$ we need some definitions and a lemma.
With Proposition \ref{GPC} in mind, we introduce the following notation:

\begin{Notation} \label{kerdeltaA0}
Denote by $\delta_{A,0}^0$ and $\delta_{A,\geq 1}^0$  the restrictions of $\delta_A^0$ to the subspaces $k((Q_A)_0\pll (Q_A)_0)$ and $k((Q_A)_0\pll (\B_A)_{\geq 1}$), respectively. We use the same notation for $\delta_{B,0}^0$ and $\delta_{B,\geq 1}^0$.
\end{Notation}

\begin{Lem}\label{Image-zero-part} Let $A=kQ_A/I_A$ be a quiver algebra. Then
\[\dim\Im(\delta^0_{A,0})=n_A-c_A,\]
where $n_A=|(Q_A)_0|$ is the number of vertices of $Q_A$ and $c_A$ is the number of connected components
of $Q_A$.
\end{Lem}

\begin{proof}
It is enough to assume that $A$ is indecomposable. Indeed, if it holds for each block $A_i$ of $A$, then
\begin {equation*}
\dim \Im(\delta^0_{A,0})= \sum_{A_i}(|(Q_{A_i})_0|-1)=|(Q_A)_0|-c_A.
\end {equation*}
Hence assume $A$ is indecomposable.
Note that: 
\[\dim k((Q_A)_0\pll (Q_A)_0)=|(Q_A)_0|= \dim  \Im(\delta^0_{A,0})+\dim  \Ker(\delta^0_{A,0}).\]
Consequently, it is enough to show that $\dim \Ker(\delta^0_{A,0})=1$.

Let $Z(A)$ be the
center of $A$. There is a bijective linear map
$\Ker(\delta_A^0)\to Z(A)$ which sends  $\sum_{i\in (Q_A)_0} \lambda_i\cdot e_i \pll p_i\ $ to $ \sum_{i\in (Q_A)_0}\lambda_i p_i$, which restricts  to a bijection $\Ker(\delta_{A,0}^0)\xrightarrow{\sim} Z(A)\cap k(Q_A)_0$  given by $\sum_{i\in (Q_A)_0} \lambda_i\cdot e_i\pll e_i\mapsto \sum_{i\in (Q_A)_0} \lambda_i e_i$.
A linear combination $z=\sum_{i\in (Q_A)_0} a_i e_i$ commutes with an arrow
$\alpha$ exactly when $a_{s(\alpha)}=a_{t(\alpha)}$. Since $Q$ is connected, we have that $\{s(\alpha),t(\alpha)\mid \alpha\in (Q_A)_1\}=(Q_A)_0$. Hence, $a_i=a_j$ for all $i,j\in (Q_A)_0$. Therefore, $Z(A)\cap k(Q_A)_0$ has dimension $1$.
\end{proof}

Let $p$ be a path between  $e_1$ and $e_n$ in $\B_A$. Then $p'$ is an oriented cycle at $f_1$ in $Q_B$. If $p$ is a path from $e_1$ to $e_n$, then we have  \[\delta_{B}^{0}(f_1 \pll  p')=\sum_{\substack{s(a)=e_n,\\a\in (Q_{A})_1,\\ap\in \B_A}}a'\pll  a'p'-\sum_{\substack{t(b)=e_1,\\b\in (Q_{A})_1,\\pb\in \B_A}}b'\pll  p'b'. \]
%$$\delta_{B}^{0}(f_1 \pll  p')=\sum_{s(a)=e_n,a\in (Q_{A})_1,ap\in \B_A}a'\pll  a'p'-\sum_{t(b)=e_1,b\in (Q_{A})_1,pb\in \B_A}b'\pll  p'b'.$$
Note that we have omitted some zero terms in the above sum, for example, if $d\in (Q_{A})_1$ is an arrow starting at $e_1$, then $d'\pll  d'p'$ appears as a term in the above sum, however, it is zero since $d'p'$ lies in $I_B$.
If $p$ is a path from $e_n$ to $e_1$, then we have
\[\delta_{B}^{0}(f_1 \pll  p')=\sum_{\substack{s(a)=e_1,\\a\in (Q_{A})_1,\\ap\in \B_A}}a'\pll  a'p'-\sum_{\substack{t(b)=e_n,\\b\in (Q_{A})_1,\\pb\in \B_A}}b'\pll  p'b'. \]
%$$\delta_{B}^{0}(f_1 \pll  p')=\sum_{s(a)=e_1,a\in (Q_{A})_1,ap\in \B_A}a'\pll  a'p'-\sum_{t(b)=e_n,b\in (Q_{A})_1,pb\in \B_A}b'\pll  p'b'.$$
As in the previous case, we have omitted some zero terms in the above sum. Moreover, in both cases, $\delta_{B}^{0}(f_1 \pll  p')$ is zero if and only if $ap,pa\in I_A$ for all $a\in (Q_{A})_1$. This observation leads to the following definition:

\begin{Def}\label{special-path}
Let $A=kQ_A/I_A$ be a quiver algebra and let $B=kQ_B/I_B$ be a radical embedding. Let $p$ be a path between $e_1$ and $e_n$ in $\B_A$. We call $p$ an \emph{extendable path} between $e_1$ and $e_n$ in $Q_A$ if $\delta_{B}^{0}(f_1 \pll  p')\neq 0$, or equivalently, if there exists some $a\in (Q_{A})_1$ such that $ap\notin I_A$ or $pa\notin I_A$. 

We denote by 
\begin{itemize}
	\item $\E$ the set of extendable paths between $e_1$ and $e_n$ in $Q_A$,
	\item $V_{\E}$ the $k$-subspace of $\Im(\delta^0_B)$ generated by the elements $\delta_B^0(f_1 \pll  p')$ for $p \in$ $\E$,
	\item $\e$ the dimension of $V_{\E}$.
\end{itemize}
\end{Def}

In the literature, extendable paths are also called \emph{non-maximal paths}.  It is worth noting that in \cite{LRW2}, we use the term \emph{special path} rather than extendable path.

\begin{Lem}\label{summand-disjoint}
Let $p$ be an extendable path between $e_1$ and $e_n$ in $\B_A$ and let $q$ be a path in $\B_A\setminus\E$ such that $q'$ is an oriented cycle in $Q_B$. Then the set of the summands of $\delta_B^0(f_1\pll p')$ and the set of the summands of $\delta_B^0(f_i\pll q')$ $(1\leq i\leq n-1)$ are disjoint.
\end{Lem}

\begin{proof}
Without loss of generality, we assume that $p$ is an extendable path from $e_1$ to $e_n$. Then $$\delta_{B}^{0}(f_1 \pll  p')=\sum_{\substack{s(\alpha)=e_n,\\\alpha\in (Q_{A})_1,\\ \alpha p\in \B_A}}\alpha'\pll  \alpha'p'-\sum_{\substack{t(\beta)=e_1,\\ \beta\in (Q_{A})_1,\\ p\beta\in \B_A}}\beta'\pll  p'\beta',$$ $$\delta_{B}^{0}(f_i \pll  q')=\sum_{\substack{s(a')=f_i,\\ a'\in (Q_{B})_1,\\ a'q'\in \B_B}}a'\pll  a'q'-\sum_{\substack{t(b')=f_i,\\ b'\in (Q_{B})_1,\\ q'b'\in \B_B}}b'\pll  q'b'.$$ Note that $\alpha'\pll  \alpha'p'\neq a'\pll  a'q'$, otherwise, $\alpha'=a'\in (Q_B)_1$ and $\alpha'p'=a'q'\in \B_B$ which imply that $\alpha=a\in (Q_A)_1$ and $\alpha p=aq\in \B_A$ by the bijection between $(\B_A)_{\geq 1}$ and $(\B_B)_{\geq 1}$. Moreover, the equality $\alpha p=aq\in \B_A$ implies that $p\pll q$. Hence $q$ is also a path in $\B_A$ from $e_1$ to $e_n$ and $aq\notin I_A$ for an arrow $a$. This means that $q\in \E$ , a contradiction. In addition, we have $\alpha'\pll  \alpha'p'\neq b'\pll  q'b'$, otherwise $\alpha=b \in (Q_A)_1$ and $\alpha p=qb$ $\in \B_A$ which implies $e_1=s(p)=s(\alpha p)=s(q\alpha)=s(\alpha)=e_n$,  a contradiction. We can similarly show that $\beta'\pll  p'\beta'\neq a'\pll  a'q'$ and $\beta'\pll  p'\beta'\neq b'\pll  q'b'$.
\end{proof}

\begin{Rem} \label{dimension-zsp}
\begin{itemize}
\item[$(1)$] The dimension of $V_{\E}$ is less than or equal to the number of extendable paths, that is, $\e\le |\E|$. This follows from the fact that the summands of $\delta_B^0(f_1\pll p')$ and of $\delta_B^0(f_1\pll q')$ may cancel each other out for $p,q\in \E, \  p\neq q$ (cf. Example \ref{eg-not-equal}). However, for monomial algebras, this lemma holds for any paths (between $e_1$ and $e_n$) $p,q\in \B_A$ with $p\neq q$, and therefore in this case we have $\e=|\E|$. The reason is that for monomial algebras we can simplify the above proof as follows: $\alpha=b$ and $\alpha p=qb$ will cause $p$ to be of the form $x_m\cdots x_2\alpha$, where $x_2, \cdots, x_m\in (Q_{A})_1$. It follows that $s(p)=s(\alpha)=e_n$, which is a contradiction.
\item[$(2)$] If $e_1$ and $e_n$ belong to different blocks of $A$ or $A$ is a radical square zero algebra, then clearly $\E=\emptyset$.
\item[$(3)$]In general, the natural number $\e$ could be arbitrarily large, see Example \ref{ex-image-kernel}.
\end{itemize}
\end{Rem}

We can now compare the dimensions of $\Im(\delta^0_A)$ and $\Im(\delta^0_B)$:

\begin{Prop}\label{Iamge-delta-zero}
Let $A=kQ_A/I_A$ be a quiver algebra and let $B=kQ_B/I_B$ be a radical embedding. Then $$\dim \Im(\delta^0_A) = \dim \Im(\delta^0_B)+1+c_B-c_A-\e.$$ In particular, if we glue $e_1$ and $e_n$ from the same block of $A$, then $$\dim \Im(\delta^0_A)=\dim \Im(\delta^0_B)+1-\e;$$ if we glue $e_1$ and $e_n$ from different blocks of $A$, then $$\dim \Im(\delta^0_A)=\dim \Im(\delta^0_B).$$
\end{Prop}

\begin{proof} As usual, the vertices of $Q_{A}$ are  $e_1,\cdots,e_n$ and the vertices of $Q_{B}$ are $f_1,\cdots,f_{n-1}$, where $f_1$ is obtained by gluing $e_1$ and $e_n$. We begin with describing a generating set of $\Im(\delta_{A}^{0})$ (resp. of $\Im(\delta_{B}^{0})$).

Let $e_i \pll  p \in k((Q_A)_{0}\pll  \B_A)$. We consider two cases, depending on whether \( p = e_i \) or \( p \neq e_i \).

($a_1$) If $p=e_i$ ($1\leq i\leq n$), then we have $$\delta_{A}^{0}(e_i \pll  e_i)=\sum_{\substack{s(a)=e_i,\\ a\in (Q_{A})_1}}a\pll  a-\sum_{\substack{t(b)=e_i,\\ b\in (Q_{A})_1}}b\pll b.$$
By Lemma \ref{Image-zero-part}, the subspace $\Im$($\delta^0_{(A)_0}$) of $\Im(\delta^0_A)$ generated by the elements of the form $\delta^0_A(e_i\pll  e_i)$ has dimension $n_A-c_A$.

($a_2$) If $p\neq e_i$, then $p$ is an oriented cycle at $e_i$ and $$\delta_{A}^{0}(e_i \pll  p)=\sum_{\substack{s(a)=e_i,\\ a\in (Q_{A})_1,\\ ap\in \B_A}}a\pll  ap-\sum_{\substack{t(b)=e_i,\\ b\in (Q_{A})_1,\\ pb\in \B_A}}b\pll pb.$$
It is clear that $$\Im(\delta^0_A) = \Im(\delta^0_{A,0})\ \oplus\ \Im(\delta^0_{A,\geq 1}).$$

Similarly, we let $f_i \pll  q \in k((Q_B)_{0}\pll  \B_B)$ and consider four cases.

($b_1$) If $q=f_i$ ($1\leq i\leq n-1$), then we have $$\delta_{B}^{0}(f_i \pll  f_i)=\sum_{\substack{s(a')=f_i,\\ a'\in (Q_{B})_1}}a'\pll  a'-\sum_{\substack{t(b')=f_i,\\ b'\in (Q_{B})_1}}b'\pll b'.$$
By Lemma \ref{Image-zero-part}, the subspace $\Im$($\delta^0_{B,0}$) of $\Im(\delta^0_B)$ generated by the elements of the form $\delta^0_B(f_i\pll  f_i)$ has dimension $n_B-c_B$.

($b_2$) If $q$ is an oriented cycle at $f_i$ and $i\neq 1$, then by Proposition \ref{parallel-paths-in-monomial-algebras} we have $q=p'$ for some oriented cycle $p\in \B_A$ at $e_i$ ($2\leq i\leq n-1$). Therefore $$\delta_{B}^{0}(f_i \pll  p')=\sum_{\substack{s(a')=f_i,\\ a'\in (Q_{B})_1}}a'\pll  a'p'-\sum_{\substack{t(b')=f_i,\\ b'\in (Q_{B})_1}}b'\pll  p'b'=\varphi_1(\delta_{A}^{0}(e_i \pll p)).$$

($b_3$) If $q$ is an oriented cycle at $f_1$ such that $q=p'$, for some oriented cycle $p\in \B_A$ at $e_1$, then $$\delta_{B}^{0}(f_1 \pll  p')=\sum_{\substack{s(a')=f_1,\\ a'\in (Q_{B})_1,\\ a'p'\in \B_B}}a'\pll  a'p'-\sum_{\substack{t(b')=f_1,\\ b'\in (Q_{B})_1,\\ p'b'\in \B_B}}b'\pll  p'b'=\varphi_1(\delta_{A}^{0}(e_1 \pll p)).$$
If $q$ is an oriented cycle at $f_1$ such that $q=p'$, for some oriented cycle $p\in \B_A$ at $e_n$, then $$\delta_{B}^{0}(f_1 \pll p')=\sum_{\substack{s(a')=f_1,\\ a'\in (Q_{B})_1,\\ a'p'\in \B_B}}a'\pll  a'p'-\sum_{\substack{t(b')=f_1,\\ b'\in (Q_{B})_1,\\ p'b'\in \B_B}}b'\pll  p'b'=\varphi_1(\delta_{A}^{0}(e_n \pll  p)).$$

($b_4$) If $q$ is an oriented cycle at $f_1$ with $q=p'$ for some path $p$ between $e_1$ and $e_n$ in $\B_A$, then we assume that $p$ is an extendable path since otherwise $\delta_{B}^{0}(f_1 \pll  p')$ is zero. Note that $q$ is of the form $f_1\stackrel{a'}{\rightarrow}\cdots\stackrel{b'}{\rightarrow}f_1$ and might be a loop at $f_1$. If $p$ is a path from $e_1$ to $e_n$, then we have $$\delta_{B}^{0}(f_1 \pll  p')=\sum_{\substack{s(a)=e_n,\\ a\in (Q_{A})_1,\\ ap\in \B_A}}a'\pll  a'p'-\sum_{\substack{t(b)=e_1,\\ b\in (Q_{A})_1,\\ pb\in \B_A}}b'\pll  p'b'.$$ If $p$ is a path from $e_n$ to $e_1$, then we have $$\delta_{B}^{0}(f_1 \pll  p')=\sum_{\substack{s(a)=e_1,\\ a\in (Q_{A})_1,\\ ap\in \B_A}}a'\pll  a'p'-\sum_{\substack{t(b)=e_n,\\ b\in (Q_{A})_1,\\ pb\in \B_A}}b'\pll  p'b'.$$ In both cases, $\delta_{B}^{0}(f_1 \pll  p')$ is nonzero since $p$ is an extendable path.

In addition, we have
$$\Im(\delta^0_B) = \Im(\delta^0_{B,0})\oplus \Im(\delta^0_{B,\geq 1}).$$
We claim that 
  $$\Im(\delta^0_{B,\geq 1})=\varphi_1({\Im(\delta^0_{A,\geq 1})})\oplus V_{\E}.$$
It suffices to show that the set of the summands of $\delta_B^0(f_1\pll p')$ and the set of the summands of $\varphi_1({\Im(\delta^0_{A,\geq 1})})$ are disjoint for $p\in \E$. Since any summand of an element in $\varphi_1({\Im(\delta^0_{A,\geq 1})})$ is of the form $\varphi_1(\delta^0_A(e_i \pll q))=\delta_B^0(f_i \pll q')$, where $q$ is an oriented cycle at $e_i \ (1\leq i\leq n)$ (here we identify $f_n$ with $f_1$), the statement follows from Lemma \ref{summand-disjoint}. Note also that the map $\varphi_1: \Im(\delta^0_{A,\geq 1}) \to \Im(\delta^0_{B,\geq 1})$ is clearly injective. Therefore we have
\begin{equation}\label{eq:im1}
\dim 	\Im(\delta^0_{A,\geq 1})=\dim 	\Im(\delta^0_{B,\geq 1})-\e.
\end{equation}
By ($a_1$) and ($b_1$) and since $n_A=n_B+1$, we get
 \begin{equation}\label{eq:im2}
\dim \Im(\delta^0_{A,0})=\dim \Im(\delta^0_{B,0})+1+c_B-c_A. 
 \end{equation}
Then we have
 \begin{equation}
 	\begin{split}
\dim \Im(\delta^0_A)&= \dim \Im(\delta^0_{A,\geq 1})+\dim \Im(\delta^0_{A,0})\\
&=\dim \Im(\delta^0_{B,\geq 1})-\e+\dim \Im(\delta^0_{B,0})+1+c_B-c_A\\
&=\dim \Im(\delta^0_B)+1+c_B-c_A-\e,
 	\end{split}
 \end{equation} 
where the second equality follows from Equations (\ref{eq:im1}) and (\ref{eq:im2}). In particular, if we glue $e_1$ and $e_n$ from the same block of $A$, then we have $c_B=c_A$. If $e_1$ and $e_n$ are from two different blocks of $A$, then $c_B=c_A-1$ and $\e=0$.
\end{proof}

We obtain a corollary that will be useful for stable equivalences induced by idempotent gluings.

\begin{Cor}\label{Iamge-delta-zero-source-and-sink}
Let $A=kQ_A/I_A$ be a quiver algebra and let $B=kQ_B/I_B$ be a radical embedding. Assume one of the following two conditions holds:
\begin{itemize}
		\item[$(i)$] $e_1$ is a source and $e_n$ is a sink;
		\item[$(ii)$] $A$ is a radical square zero algebra.
\end{itemize}

Then we have $$\dim \Im(\delta^0_A)=\dim \Im(\delta^0_B)+1+c_B-c_A.$$
\noindent In particular, if we glue $e_1$ and $e_n$ from the same block of $A$, then $$\dim \Im(\delta^0_A)=\dim \Im(\delta^0_B)+1;$$ if $e_1$ and $e_n$ are from two different blocks of $A$, then we have $$\dim \Im(\delta^0_A)=\dim \Im(\delta^0_B).$$
\end{Cor}

\begin{proof} It is clear that under the condition $(i)$ or $(ii)$ there are no extendable paths between $e_1$ and $e_n$. Therefore $\e=0$. If we glue $e_1$ and $e_n$ from the same block of $A$, then $c_B=c_A$; if $e_1$ and $e_n$ are from two different blocks of $A$, then $c_B=c_A-1$. Thus, the result  follows from Proposition \ref{Iamge-delta-zero}.
\end{proof}

We will often use the following assumption in the sequel:
\begin{Assumption}\label{assump}
Let $A=kQ_A/I_A$ be a quiver algebra and let $B=kQ_B/I_B$ be a radical embedding obtained by gluing two idempotents $e_1$ and $e_n$ of $A$. For each loop $\alpha$ at $e_1$ or at $e_n$ with $\alpha^m\in \Tip(\mathcal{G}_A)$,  we have that  $\mathrm{char}(k)\nmid m$.
\end{Assumption}

Clearly, Assumption \ref{assump} holds if the characteristic of the field $k$ is zero or big enough. We now proceed to compare the Lie structures of $\Ker(\delta^1_A)$ and $\Ker(\delta^1_B)$:

\begin{Prop}\label{Kernel-delta-one-injective}
Let $A=kQ_A/I_A$ be a quiver algebra and let $B=kQ_B/I_B$ be a radical embedding. If Assumption \ref{assump} is satisfied, then there exists an injective (restricted) Lie algebra homomorphism $\Ker(\delta^1_A)\hookrightarrow$ $\Ker(\delta^1_B)$ induced from $\varphi_1: k((Q_A)_1\pll  \B_A) \to k((Q_B)_1\pll  \B_B)$, which we still denote by $\varphi_1$.
\end{Prop}

\begin{proof} First we notice that $I_A=\langle \mathcal{G}_A\rangle$ and $I_B=\langle \mathcal{G}_B\rangle$, and by Proposition \ref{corresponding-RGB} we can write $\mathcal{G}_B=\mathcal{G}_A\cup Z_{new}$, where $Z_{new}=\{b'c'\mid b',c'\in (Q_B)_1, t(c')=f_1=s(b'), bc\notin \B_A\}$.
Having the diagram  $(*)$  in mind,
let $\alpha \pll  p \in k((Q_A)_{1}\pll  \B_A)$ and let $\varphi_1(\alpha \pll  p)=\alpha' \pll  p'$ be the corresponding element in $k((Q_B)_1\pll  \B_B)$. On the one hand, we have
\[\varphi_2\left(\delta_{A}^{1}(\alpha \pll  p)\right)=\varphi_2\left(\sum_{\substack{r \in \mathcal{G}_A,\\ q\in \Supp(r)}}c_r(q)\cdot\Tip(r)\pll q^{\alpha \pll p}\right)=\sum_{\substack{r \in \mathcal{G}_A,\\ q\in \Supp(r)}}c_r(q)\cdot\Tip(r)'\pll q'^{\alpha' \pll p'}; \]
On the other hand, we have
\begin{align*}
	\delta_{B}^{1}\left(\varphi_1(\alpha \pll  p)\right)&=\delta_{B}^{1}(\alpha' \pll  p')\\ &= \sum_{\substack{r' \in \mathcal{G}_B,\\q'\in \Supp(r')}}c_{r'}(q')\cdot\Tip(r')\pll q'^{\alpha' \pll p'} \\ &=\sum_{\substack{r \in \mathcal{G}_A,\\ q\in \Supp(r)}}c_r(q)\cdot\Tip(r)'\pll q'^{\alpha' \pll p'}+\sum_{r' \in Z_{new}}r'\pll r'^{\alpha' \pll p'}\\ &=\varphi_2\left(\delta_{A}^{1}(\alpha \pll  p)\right)+\sum_{r' \in Z_{new}}r'\pll r'^{\alpha' \pll p'}.
\end{align*}
We claim that $\varphi_1$ induces a well-defined map from $\Ker(\delta_A^1)$ to $\Ker(\delta_B^1)$ under Assumption \ref{assump}. Since each element in $\Ker(\delta_A^1)$ is a $k$-linear combination of  elements $\alpha\pll p\in (Q_A)_1\pll \B_A$, it suffices to consider the following four cases.

($c_1$) If $\alpha$ is a loop at $e_i$, for $2\leq i\leq n-1$, and $p=e_i$ or $p$ is an oriented cycle at $e_i$, then $\sum_{r' \in Z_{new}}r'\pll r'^{\alpha' \pll p'}=0$. Indeed, $\alpha'$ does not appear in any $r'\in Z_{new}$. Therefore $\varphi_2(\delta_{A}^{1}(\alpha \pll p))=\delta_{B}^{1}(\varphi_1(\alpha \pll  p))$.

($c_2$) If $\alpha$ is a loop at $e_1$ (resp. $e_n$) and $p=e_1$ (resp. $p=e_n$), then Lemma \ref{loop} shows that there exists a relation $r$ in $\mathcal{G}_A$ such that $\Tip(r)=\alpha^m$ and $\alpha^{m-1}\in \mathrm{NonTip}(I_A)$ for some integer $m\geq 2$. By Assumption \ref{assump} we have that char$(k)\nmid m$. Let us assume that $\alpha$ is a loop at $p=e_1$. Then $\delta_{A}^{1}(\alpha\pll e_1)$ contains a nonzero summand $m\Tip(r)\pll \alpha^{m-1}=m\alpha^m\pll \alpha^{m-1}$. In addition, by definition of $\delta_A^1$ we have that it does not exist an element $a\pll q\in (Q_A)_1\pll \B_A$, $a\pll q\neq \alpha\pll e_1$, such that $\delta_A^1(a\pll q)$ has a summand of the form $\alpha^m\pll \alpha^{m-1}$. Hence $\alpha \pll  e_1$ cannot appear as a summand of an element of $\Ker(\delta_A^1)$. Note that $\alpha'$ appears in some $r'\in Z_{new}$. Therefore $\delta_{B}^{1}(\varphi_1(\alpha \pll  e_1))=\delta_{B}^{1}(\alpha' \pll  f_1)$ contains a summand $r'\pll r'^{\alpha'\pll f_1}$, which by the same argument as above, cannot be cancelled in $\Im(\delta_B^1)$. Therefore,  $\varphi_1(\alpha \pll  e_1)$ cannot appear as a summand of an element of $ \Ker(\delta_{B}^{1})$. A similar result holds if  $\alpha$ is a loop at $p=e_n$.

($c_3$) If $\alpha$ is a loop at $e_1$ (resp. $e_n$) and $p$ is an oriented cycle at $e_1$ (resp. $e_n$), then once we replace $\alpha'$ in any $r'\in Z_{new}$ by $p'$, $r'$ becomes a path in $Q_B$ that still contains some relation in $Z_{new}$. Hence $\sum_{r' \in Z_{new}}r'\pll  {r'}^{\alpha' \pll  p'}=0$. Therefore $\varphi_2(\delta_{A}^{1}(\alpha \pll  p))=\delta_{B}^{1}(\varphi_1(\alpha \pll  p))$.

($c_4$) If $\alpha$ is a an arrow which is not a loop such that $\alpha'$ appears in some $r'\in Z_{new}$ and if $p\in \B_A$ is a  path parallel to $\alpha$, then, by the same argument as in $(c_3)$, the element obtained from $r'$ by replacing $\alpha'$ by $p'$ is not in $\B_B$. Hence  $\sum_{r' \in Z_{new}}r'\pll  {r'}^{\alpha' \pll  p'}=0$ and consequently $\varphi_2(\delta_{A}^{1}(\alpha \pll  p))=\delta_{B}^{1}(\varphi_1(\alpha \pll  p))$.

The above discussion shows that, for each element $x=\sum_i \lambda_i\cdot a_i\pll p_i\in \Ker(\delta_A^1)$ with $\lambda_i\in k$,  if Assumption \ref{assump} is satisfied, then the summand $a_i\pll p_i$ of $x$ could come from the cases $(c_1), (c_3)$ and $(c_4)$. As a consequence, we obtain that $\delta_B^1(\varphi_1(x))=\varphi_2(\delta_A^1(x))=0$. So $\varphi_1(x)\in \Ker(\delta_B^1)$. In other words, if Assumption \ref{assump} is satisfied, then there is a $k$-linear map $\varphi_1: \Ker(\delta_{A}^{1})\longrightarrow \Ker(\delta_{B}^{1})$ induced from  $\varphi_1: k((Q_A)_1\pll  \B_A) \to k((Q_B)_1\pll  \B_B)$. It is also clear that $\varphi_1: \Ker(\delta_{A}^{1})\longrightarrow \Ker(\delta_{B}^{1})$ is injective and preserves the Lie bracket, since $\varphi_1: k((Q_A)_1\pll  \B_A) \to k((Q_B)_1\pll  \B_B)$ preserves parallel paths.
\end{proof}

Since the characteristic condition is only needed in $(c_2)$ to ensure that $\varphi_1:\Ker(\delta_A^1)\to \Ker(\delta_B^1)$ is well-defined, we do not need Assumption \ref{assump} in Proposition \ref{Kernel-delta-one-injective} in some cases. For instance, we have

\begin{Cor}\label{no-Assump1}
	Let $B$ be a radical embedding of $A$  where the base field $k$ has arbitrary characteristic. Then $\Ker(\delta^1_A)\hookrightarrow$ $\Ker(\delta^1_B)$ is an injective Lie algebra homomorphism under one of the following conditions:
	\begin{enumerate}
\item There is no loop both at $e_1$ and at $e_n$. In particular if $e_1$ (resp. $e_n$) is a source vertex and $e_n$ (resp. $e_1$) is a sink vertex.
	\item If $A$ (and hence also $B$) is a radical square zero algebra excluding the case that $e_1$ and $e_n$ are glued from different blocks of $A$ with one of these blocks isomorphic to $k[x]/(x^2)$.
	\end{enumerate}
\end{Cor}

\begin{proof}
The proof of (1) is obvious so we proceed to prove (2). We divide the proof in two cases:
	
(a) In the excluded case, we show that $\chr(k)\neq 2$ is necessary to ensure $\varphi_1:\Ker(\delta_A^1)\to \Ker(\delta_B^1)$ is well-defined. Otherwise, assume $\chr(k)=2$ and we glue $e_1$ and $e_n$ from different blocks, say $e_1\in A_1$ and $e_n\in A_2$, of $A$ such that $A_1$ or $A_2$ is isomorphic to $k[x]/(x^2)$. Without loss of generality, we assume that $A_1$ is isomorphic to $k[x]/(x^2)$, that is,  the quiver of $A_1$ is given by a loop $\alpha$ and a vertex $e_1$. In this case we have that $\delta_A^1(\alpha\pll e_1)=2\alpha^2\pll \alpha=0$.  Note that there  exists an  arrow $a$ in $(Q_A)_1$ such that either $t(a)=e_n$ or $s(a)=e_n$. Thus, either $\alpha' a'$ or $a'\alpha'$ lies in $Z_{new}$. This give to a nonzero summand $\alpha' a'\pll a'$ or $a'\alpha'\pll a'$ of $\delta_B^1(\varphi_1(\alpha\pll e_1))$, that is, $\varphi_1(\alpha\pll e_1) \notin \Ker(\delta_B^1)$. So $\varphi_1:\Ker(\delta_A^1)\to \Ker(\delta_B^1)$ is not well-defined.

(b) If we are not in the excluded case, then we show that $\varphi_1:\Ker(\delta_A^1)\to \Ker(\delta_B^1)$ is well-defined in any $\chr(k)$. In fact, now either `$e_1$ and $e_n$ are glued from the same block' or `we glue $e_1$ and $e_n$ from different blocks, say $e_1\in A_1$ and $e_n\in A_2$, of $A$ such that neither $A_1$ nor $A_2$ is isomorphic to $k[x]/(x^2)$'. Following the proof of Proposition~\ref{Kernel-delta-one-injective}, it suffices to show that no element $x\in \Ker(\delta_A^1)$ has $\alpha \pll s(\alpha)$ as a summand whenever $s(\alpha)\in\{e_1,e_n\}$, where $\alpha$ is a loop at $e_1$ or $e_n$. Note that for such loop $\alpha$, there exists either an arrow $\beta \neq \alpha$ starting at $s(\alpha)$ or an arrow $\gamma \neq \alpha$ ending at $s(\alpha)$. Since $A$ is a radical square zero algebra, we have $\beta\alpha \in \mathcal{G}_A$ or $\alpha\gamma \in \mathcal{G}_A$ for such $\beta$ and $\gamma$. Therefore, either $\beta\alpha \pll \beta$ or $\alpha\gamma \pll \gamma$ is a nonzero summand of $\delta_A^1(\alpha \pll s(\alpha))$. Moreover, such summand cannot be cancelled by the summands of $\delta_A^1(a\pll p)$ for any $a\pll p\neq \alpha \pll s(\alpha)$.  Consequently, each summand of $x\in \Ker(\delta_A^1)$ is different from $\alpha\pll s(\alpha)$. 
\end{proof}

In order to describe the elements in $\Ker(\delta_{B}^{1})$ which are in the complement of the subspace $\varphi_1(\Ker(\delta_{A}^{1}))$, we introduce the following definition.

\begin{Def}\label{elements-in-Z-spp}
Let $A=kQ_A/I_A$ be a quiver algebra and let $B=kQ_B/I_B$ be a radical embedding. Let $\alpha$ be an arrow of $Q_A$ and $p$ be a path in $\B_A$. We call $(\alpha,p)$  an \emph{almost parallel pair} with respect to the gluing of $e_1$ and $e_n$ if the following two conditions are satisfied:

$(1)$ $\alpha \notparallel p$ in $Q_A$;

$(2)$ $\alpha' \pll  p'$ in $Q_B$.
 
We denote by 
	\begin{itemize}
		\item $\AP$ the set of all almost parallel pairs $(\alpha,p)$ in $(Q_A)_1\pll \B_A$ with respect to the gluing of $e_1$ and $e_n$.
		\item $k(\AP)$ the $k$-subspace of $k((Q_B)_1\pll \B_B)$ generated by the elements $\alpha' \pll  p'$, where $(\alpha,p) \in \AP$.
		\item $V_{\AP}$ the intersection of $k(\AP)$ and $\Ker(\delta_B^1)$.
		\item $\ap$ the dimension of the $k$-subspace $V_{\AP}$ of $\Ker(\delta_B^1)$.
	\end{itemize}
\end{Def}

Observe that $V_{\E}$ is a subspace of $V_{\AP}$ and therefore we always have $\ap\ge\e$. Note that every nonzero element of $V_{\AP}$ is a linear combination of parallel paths corresponding to almost parallel pairs (cf. Example \ref{linear-combination-zspp}). Moreover, conditions (1) and (2) imply that $\alpha$ is starting from $e_1$, or ending at $e_1$, or starting from $e_n$, or ending at $e_n$. Note also that the notion of almost parallel pairs leads to various possible configurations of the pairs $(\alpha,p)\in (Q_A)_1\pll \B_A$, see Example \ref{example-type-special-pair}.

It is worth noting that in \cite{LRW2} the authors use the term `special pair' rather than almost parallel pair.

\begin{Rem}\label{spp-rad-square-zero}
Let $A$ be a radical square zero algebra. Then there are no extendable paths in $Q_A$. However,  almost parallel pairs could appear when we glue two vertices in $Q_A$ either from the same block (see Examples \ref{rad-1} and \ref{ex-converse-coro}) or from two different blocks of $A$.
 
Indeed, if we glue $e_1$ and $e_n$ from two different blocks of $A$, then 
\begin{align*} \AP=&\{(\alpha,e_n),(\alpha,\beta),(\beta,e_1),(\beta,\alpha)\mid \text{$\alpha$ (resp. $\beta$) is a loop at $e_1$ (resp. $e_n$)} \},\\ V_{\AP}:=&k(\AP)\cap\Ker(\delta_B^1)=\langle\alpha'\pll \beta',\beta'\pll \alpha'\mid \text{$\alpha$ (resp. $\beta$) is a loop at $e_1$ (resp. $e_n$)} \rangle. 
\end{align*} 
It follows that in this case if there is a loop at either $e_1$ or $e_n$, then the set $\AP$ of almost parallel pairs is not empty. In addition, note that if at least one of $e_1$ or $e_n$ has no loop, then $V_{\AP}=0$. 
\end{Rem}

\begin{Prop}\label{Kernel-delta-one-difference}
Let $A=kQ_A/I_A$ be a quiver algebra and let $B=kQ_B/I_B$ be a radical embedding. If Assumption \ref{assump} is satisfied, then we have a decomposition  $$\Ker(\delta_B^1)=\varphi_1(\Ker(\delta_{A}^{1}))\oplus V_{\AP},$$ as $k$-vector spaces and therefore $$\dim \Ker(\delta^1_B) = \dim \Ker(\delta^1_A)+\ap.$$
\end{Prop}

\begin{proof}
By Proposition \ref{Kernel-delta-one-injective}, we only need to describe the elements $\theta$ in $\Ker(\delta_{B}^{1})$ which are in the complement of the subspace $\varphi_1(\Ker(\delta_{A}^{1}))$, under Assumption \ref{assump}. According to the proof of Proposition \ref{Kernel-delta-one-injective}, we may assume that $\theta$ is a linear combination of  elements of the form $\alpha' \pll  p'$ such that $(\alpha,p)$ is a almost parallel pair with respect to the gluing of $e_1$ and $e_n$. Clearly in this case $\theta\in V_{\AP}$.  Therefore, we have the following decomposition: $\Ker(\delta_B^1)=\varphi_1(\Ker(\delta_{A}^{1}))\oplus V_{\AP}$.
\end{proof}

In several results concerning gluing a source and a sink, we have to exclude a very special case which we call `Exceptional case \ref{e-case}'.

\begin{Ec}\label{e-case}
	Let $\mathrm{char}(k)=2$, by gluing we obtain a block of $B$ of the form $k[x]/(x^{2})$, in other words, $A$ has one block which has a Gabriel quiver of type $A_2$ and we perform the gluing in this block.
\end{Ec}

\begin{Cor}\label{Kernel-delta-one-source-and-sink}
Let $A=kQ_A/I_A$ be a quiver algebra and let $B=kQ_B/I_B$ be a radical embedding obtained by gluing a source vertex $e_1$ and a sink vertex $e_n$ of $A$. Then we have $\Ker(\delta^1_A) \simeq \Ker(\delta^1_B)$ as Lie algebras, except in the Exceptional case \ref{e-case}.
\end{Cor}

\begin{proof} By Lemma \ref{bijection-parallel-arrows}, the only possible almost parallel pair with respect to the gluing of $e_1$ and $e_n$ has the form $(\alpha,e_1)$ or $(\alpha,e_n)$ such that $\alpha$ is an arrow from $e_1$ to $e_n$. Therefore $k(\AP)$ is generated by the elements of the form $\alpha'\pll f_1$. Suppose now that $\alpha'\pll f_1 \in\Ker(\delta_{B}^{1})$. Then we consider two cases. If $Q_A$ contains a connected component $e_1\stackrel{\alpha}\longrightarrow e_n$ so that $B$ has a block isomorphic to $k[x]/(x^{2})$, then $\delta_{B}^{1}(\alpha' \pll  f_1)=2r'\pll  \alpha'=0$ (where $r'=\alpha'\alpha'$) implies that $\mathrm{char}(k)=2$. If $Q_A$ is not the above case, then either there is an arrow $\beta'\neq \alpha'$ starting from $f_1$ or there is an arrow $\gamma'\neq \alpha'$ ending at $f_1$ in $Q_B$. Therefore $\delta_{B}^{1}(\alpha' \pll  f_1)$ will contain a summand $\beta'\alpha'\pll  \beta'$ or a summand $\alpha'\gamma'\pll  \gamma'$, which clearly cannot be cancelled in $\Im(\delta_{B}^{1})$, so $\alpha' \pll  f_1\notin \Ker(\delta_{B}^{1})$. It follows that $\alpha'\pll  f_1\in \Ker(\delta_{B}^{1})$ if and only if $B$ has a block isomorphic $k[x]/(x^{2})$ and  char$(k)=2$. Summarising the above discussion, we get $\ap=0$ when gluing a source and a sink and excluding the Exceptional case \ref{e-case}. The statement follows from Proposition \ref{Kernel-delta-one-difference}, Proposition \ref{Kernel-delta-one-injective} and Corollary \ref{no-Assump1} (1).	
\end{proof}

\begin{Rem}\label{exceptional-case}
For the Exceptional case \ref{e-case},  suppose that $A$ has only one block which has a Gabriel quiver of type $A_2$.
In this case $\mathrm{char}(k)=2$ and $B\simeq k[x]/(x^{2})$, $A=kQ_A$ where $Q_A$ is given by the quiver $1\stackrel{\alpha}\longrightarrow 2$. By a direct computation, we have the following:  $\Im(\delta^0_A)=\Ker(\delta^1_A)$ is 1-dimensional with $k$-basis $\{\alpha\pll \alpha\}$, $\Im(\delta^0_B)=0$ and $\Ker(\delta^1_B)$ is 2-dimensional with $k$-basis $\{\alpha'\pll f_1, \alpha'\pll \alpha'\}$. Note that $\mathcal{AP}^2_1=\{(\alpha,e_1),(\alpha,e_2)\}$ and $V_{\mathcal{AP}^2_1}=\langle \alpha'\pll  f_1 \rangle $. Therefore, $\dim \Ker(\delta^1_B)$ = $\dim \Ker(\delta^1_A)+$ $\mathfrak{ap}^2_1$.
\end{Rem}

We can finally compare the dimensions of $\HH^1(A)$ and of $\HH^1(B)$.

\begin{Thm}\label{hh1-quotient-isomorphism-for-monomial-algebras}
Let $A=kQ_A/I_A$ be a quiver algebra and let $B=kQ_B/I_B$ be a radical embedding. If Assumption \ref{assump} is satisfied, then we have $$\dim \HH^{1}(A)=\dim \HH^{1}(B)-1-\ap+\e +c_A-c_B.$$
In particular, if we glue $e_1$ and $e_n$ from the same block of $A$, then $$\dim \HH^{1}(A)=\dim \HH^{1}(B)-1-\ap+\e;$$ if $e_1$ and $e_n$ are from two different blocks of $A$, then $\HH^{1}(A)$ is a Lie subalgebra of $\HH^{1}(B)$ and $$\dim \HH^{1}(A)=\dim \HH^{1}(B)-\ap.$$
\end{Thm}

\begin{proof} 
Since $\HH^{1} \simeq \Ker(\delta^1) / \Im(\delta^0)$, the statement follows from Propositions \ref{Iamge-delta-zero} and \ref{Kernel-delta-one-difference}.
\end{proof}

\begin{Rem}\label{CiSa}
    In \cite[Theorem 1]{CS}  the authors give a formula to compute the dimension of $\HH^1(A)$ for a monomial algebra $A$ using an exact sequence in \cite[Page 98]{C3}.
    They introduce the following notions: an element $a\pll p$ in $(Q_A)_1\pll \B$ is \emph{admissible} if $a\pll p\in \Ker(\delta^1_A)$. An element $a\pll p$ in $(Q_A)_1\pll \B$ is \emph{glued} if  $p$ is a vertex or $a$ is the first or the last arrow of $p$. The image of $\delta^1$ restricted to the subspace spanned by glued elements is denoted $\Im(R_g)$. An element $a\pll p$ is called \emph{effective} if it is neither glued nor admissible. We denote by  $ ((Q_A)_1\pll \B_A)_e$ the set of effective elements. Then:
    \[
    \mathrm{dim} \HH^1(A)= |(Q_A)_1\pll\B_A|-|((Q_A)_1\pll\B_A)_e|-\dim\Im(R_g) -(|(Q_A)_0\pll\B_A| -\dim Z(A)).
    \]
   This gives another interpretation for the dimension of $V_{\AP}$ for monomial algebras, that is, $\ap=K_A-K_B$, where $K_A:=|(Q_A)_1\pll\B_A|-|((Q_A)_1\pll\B_A)_e|-\dim\Im(R_g)$ and so does for $B$.
\end{Rem}

\begin{Notation}
    We set $$Y:=\varphi_1(\Im(\delta_A^0))\oplus V_{\E} \subseteq \Ker
(\delta_B^1),$$ where $V_{\E}$ is the subspace of $\Im(\delta_{(B)_{\ge1}}^0)$ generated by the elements $\delta_{B}^{0}(f_1 \pll  p')$ for $p \in \E$.
\end{Notation}

We have the following result, which can be seen as a structural interpretation of the dimension formula in Theorem \ref{hh1-quotient-isomorphism-for-monomial-algebras}.

\begin{Thm}\label{Lie-strucutre-of-hh1} Under the conditions of Theorem \ref{hh1-quotient-isomorphism-for-monomial-algebras}, we have the following exact commutative diagram:
\begin{center}
	\begin{tikzcd}
	& 0 & 0 & 0 \\
	0 & {V_{\E}} & {V_{\AP}} & {\mathrm{Coker}(\varphi)} & 0 \\
	0 & Y & {\Ker(\delta_B^1)} & {\frac{\Ker(\delta_B^1)}{Y}} & 0 \\
	0 & {\Im(\delta_A^0)} & {\Ker(\delta_A^1)} & {\HH^1(A)} & 0, \\
	& 0 & 0 & 0
	\arrow[from=2-1, to=2-2]
	\arrow["{\tilde{\iota}_B}", from=2-2, to=2-3]
	\arrow[from=2-3, to=2-4]
	\arrow[from=2-4, to=2-5]
	\arrow[from=3-1, to=3-2]
	\arrow["{\iota_B}", from=3-2, to=3-3]
	\arrow[from=3-3, to=3-4]
	\arrow[from=3-4, to=3-5]
	\arrow[from=4-1, to=4-2]
	\arrow["{\iota_A}", from=4-2, to=4-3]
	\arrow[from=4-3, to=4-4]
	\arrow[from=4-4, to=4-5]
	\arrow["{\varphi_1|_{\Im(\delta_A^0)}}", from=4-2, to=3-2]
	\arrow["{\pi^0}", from=3-2, to=2-2]
	\arrow["{\pi^1}", from=3-3, to=2-3]
	\arrow["{\varphi_1}", from=4-3, to=3-3]
	\arrow["\pi", from=3-4, to=2-4]
	\arrow["\varphi", from=4-4, to=3-4]
	\arrow[from=5-2, to=4-2]
	\arrow[from=5-3, to=4-3]
	\arrow[from=2-2, to=1-2]
	\arrow[from=2-3, to=1-3]
	\arrow[from=2-4, to=1-4]
	\arrow[from=5-4, to=4-4]
	\end{tikzcd}
\quad \quad $(**)$
\end{center}
where $\pi^0,\pi^1$ are canonical projections, $\iota_A$ and $\iota_B$ are canonical injections, $\varphi$ is an
injective map induced from $\varphi_1$ and $\pi$ is a surjective map induced from $\pi^1$, $Y:=\varphi_1(\Im(\delta_A^0))\oplus V_{\E}$ is a subspace of $\Ker
(\delta_B^1)$. In particular, $\tilde{\iota}_B: V_{\E}\rightarrow V_{\AP}$ is an injective map. In addition,  
\begin{itemize}
    \item $Y$ is equal to $\Im(\delta_B^0)$ in the case that $e_1$ and $e_n$ are from different blocks of $A$.
    \item $Y$ contains $\Im(\delta_B^0)$ as a codimension one subspace in case that $e_1$ and $e_n$ are from the same block of $A$.
\end{itemize}
\end{Thm}

\begin{proof} By Proposition \ref{Kernel-delta-one-injective}, there exists an injective Lie algebra homomorphism $\varphi_1: \Ker(\delta^1_A)\hookrightarrow \Ker(\delta^1_B)$, which is induced from the canonical map $\varphi_1: k((Q_A)_1\pll  \B_A) \to k((Q_B)_1\pll  \B_B)$. Moreover, by Proposition \ref{Kernel-delta-one-difference}, we have the decomposition $\Ker(\delta_B^1)=\varphi_1(\Ker(\delta_{A}^{1}))\oplus V_{\AP}$.

Therefore by Proposition \ref{Iamge-delta-zero} and by the fact that $\delta_B^0(f_1\pll f_1)=\varphi_1(\delta_A^0(e_1\pll e_1))+\varphi_1(\delta_A^0(e_n\pll e_n))$, we have that $\varphi_1: \Ker(\delta^1_A)\hookrightarrow \Ker(\delta^1_B)$ restricts to an injective map $$\varphi_1|_{\Im(\delta_A^0)}:\Im(\delta_A^0)=\Im(\delta^0_{A,0})\oplus \Im(\delta^0_{A,\geq 1}) \hookrightarrow X\oplus \Im(\delta^0_{B,\geq 1})\subseteq \Ker(\delta_B^1),$$ where $X$ is the subspace of $\Ker(\delta_B^1)$ generated by the elements $\varphi_1(\delta_A^0(e_1\pll e_1))$, $\varphi_1(\delta_A^0(e_n\pll e_n))$ and $\delta_B^0(f_i\pll f_i)$ $(2\le i\le n-1)$.

Note that $X\oplus \mathrm{\Im}(\delta^0_{B,\geq 1})=Y$. In addition, the dimension of $X$ is equal to $\dim \Im(\delta^0_{A,0})$. It follows from Lemma \ref{Image-zero-part} that if $e_1$ and $e_n$ are from two different blocks of $A$, then $X=\Im(\delta^0_{B,0})$. By the same reasoning, if $e_1$ and $e_n$ are from the same block of $A$, then $\Im(\delta^0_{B,0})\subseteq X$ has codimension one in $X$. Therefore $\Im$($\delta_B^0$)=$\Im$($\delta^0_{B,0}$)$\oplus \Im$($\delta^0_{B,\geq 1}$) is equal to $Y$ if $e_1$ and $e_n$ are from two
different blocks of $A$, and $\Im(\delta_B^0)$ has codimension one in $Y$ if $e_1$ and $e_n$ are from the same block of $A$. 

If $p$ is an extendable path from $e_1$ to $e_n$, then each summand $a'\pll a'p'$ (or $b'\pll p'b'$) of $\delta_{B}^{0}(f_1 \pll  p')$, where $a$ is an arrow starting from $e_n$ such that $ap\in \B_A$ (or where $b$ is an arrow ending at $e_1$ such that $pb\in \B_A$), is induced from an almost parallel pair $(a,ap)$ (or $(b,pb)$). In the case that $p$ is an extendable path from $e_n$ to $e_1$, we have the similar conclusion. Therefore the canonical injective map $Y\hookrightarrow \Ker(\delta^1_B)$ restricts to an injective map $V_{\E}\hookrightarrow V_{\AP}$.

Hence we obtain the exact commutative diagram $(**)$.
\end{proof}

\begin{Lem}\label{lem:Lieideal}
Under the conditions of Theorem \ref{Lie-strucutre-of-hh1}, the space $Y$ is a Lie ideal of $\Ker(\delta_B^1)$ if and only if $[\varphi_1(\Im(\delta_A^0)),V_{\AP}] \subseteq Y$.
\end{Lem}
\begin{proof}
        By the definition of $Y$ we have that
        \begin{equation*}
            \begin{split}
                [Y,\Ker(\delta_B^1)]&=[\varphi_1(\Im(\delta_A^0)),\Ker(\delta_B^1)]+[V_{\E},\Ker(\delta_B^1)]\\&=
                [\varphi_1(\Im(\delta_A^0)),\varphi_1(\Ker(\delta_A^1))]+[\varphi_1(\Im(\delta_A^0)),V_{\AP}]+[V_{\E},\Ker(\delta_B^1)]\\&
                =\varphi_1([\Im(\delta_A^0),\Ker(\delta_A^1)])+[\varphi_1(\Im(\delta_A^0)),V_{\AP}]+[V_{\E},\Ker(\delta_B^1)]\\ &\subseteq \varphi_1(\Im(\delta_A^0))+[\varphi_1(\Im(\delta_A^0)),V_{\AP}]+[\Im(\delta_B^0),\Ker(\delta_B^1)]\\ &\subseteq \varphi_1(\Im(\delta_A^0))+[\varphi_1(\Im(\delta_A^0)),V_{\AP}]+\Im(\delta_B^0)\\ &\subseteq Y+[\varphi_1(\Im(\delta_A^0)),V_{\AP}]+Y,
            \end{split}
        \end{equation*}
        where the second equality follows from Proposition \ref{Kernel-delta-one-difference}, the third equality follows from the fact that $\varphi_1$ is a Lie algebra homomorphism, and the last three equalities follow  from the definition of $Y$ and the fact that $\Im(\delta^0)$ is a Lie ideal of $\Ker(\delta^1)$.
    \end{proof}

\begin{Thm}\label{Lie-strucutre2-of-hh1}
Let $A=kQ_A/I_A$ be a quiver algebra and let $B=kQ_B/I_B$ be a radical embedding. Suppose that Assumption \ref{assump} is satisfied. If $V_{\AP}=V_{\E}$, then  
\begin{itemize}
    \item $Y$ is a Lie ideal of $\Ker(\delta_B^1)$ and
    \item there is a Lie algebra epimorphism from $\HH^1(B)$ to $
\Ker(\delta_B^1)/Y\simeq \HH^1(A)$ with kernel $\mathcal{I}:=Y/\Im(\delta_B^0)$, where $\mathcal{I}$ is zero if $e_1$ and $e_n$ are from two different blocks of $A$ and $\dim\mathcal{I}=1$ if $e_1$ and $e_n$ are
from the same block of $A$. In particular, we have a short exact sequence of Lie algebras:
\begin{equation*}
\begin{tikzcd}[column sep=small]
0\arrow[r] &\mathcal{I} \arrow[r] &\HH^1(B)\arrow[r] &\HH^1(A)\arrow[r] &0.
\end{tikzcd}
\end{equation*}
\end{itemize}

\end{Thm}

\begin{proof}
    For the first part, note that if $e_1$ and $e_2$ are from two different blocks then by Theorem \ref{Lie-strucutre-of-hh1}  we have  $Y= \Im(\delta^0_B)$. Hence $Y$ is a Lie ideal of $\Ker(\delta^1_B)$. If $e_1$ and $e_2$ are from the same block, then the statement follows from Lemma \ref{lem:Lieideal} and from the fact that $V_{\E}$ is a $k$-subspace of $\Im(\delta^0_B)$. The second part of the proof follows from Theorem \ref{Lie-strucutre-of-hh1}. 
\end{proof}

\begin{Rem}\label{rmk:basisofI}
  Note that the one-dimensional ideal $\mathcal{I}:=Y/\Im(\delta_B^0)$ of $\HH^1(B)$ in Theorem \ref{Lie-strucutre2-of-hh1} is spanned by $\varphi_1(\delta_A^0(e_1\pll e_1))=\varphi_1(\sum_{\alpha\in (Q_A)_1e_1}\alpha\pll \alpha-\sum_{\beta \in e_1(Q_A)_1}\beta\pll \beta)$ modulo $\Im(\delta_B^0)$ (cf. the proof of Proposition \ref{Iamge-delta-zero}). In case we glue a source $e_1$ and a sink $e_n$, the ideal $\mathcal{I}$ is spanned by $\sum_{\alpha\in (Q_A)_1e_1}\alpha'\pll \alpha'$ modulo $\Im(\delta_B^0)$.
 \end{Rem}
\begin{Cor}	\label{cor:iso}
Let $A=kQ_A/I_A$ be a quiver algebra and let $B=kQ_B/I_B$ be a radical embedding. Suppose that Assumption \ref{assump} is satisfied. If $V_{\AP}=V_{\E}$, then
	\[
	\HH^1(A)/\mathrm{rad}(\HH^1(A))\simeq \HH^1(B)/\mathrm{rad}(\HH^1(B)).
	\]
\end{Cor}

\begin{proof}
Since by Theorem \ref{Lie-strucutre2-of-hh1} the ideal $\mathcal{I}$ is at most one-dimensional, then $\mathcal{I}$ is solvable. Since the radical contains every solvable ideal, then $\mathcal{I}\subseteq \mathrm{rad}(\HH^1(B))$. Hence by quotienting by the radical we obtain the desired isomorphism.
\end{proof}

\begin{Cor}\label{quotient-lie-algebra-source-and-sink}
	Let $A=kQ_A/I_A$ be a quiver algebra and let $B=kQ_B/I_B$ be a radical embedding obtained by gluing a source vertex $e_1$ and a sink vertex $e_n$ of $A$. Then we have $$\dim \HH^{1}(A)=\dim \HH^{1}(B)+c_A-c_B-1,$$ except in the Exceptional case \ref{e-case}. In particular, if we glue $e_1$ and $e_n$ from two different blocks of $A$, then there is a (restricted) Lie algebra isomorphism $$\HH^{1}(A)\simeq \HH^{1}(B);$$ if $e_1$ and $e_n$ are from the same block of $A$, then $$\HH^{1}(A)\simeq \HH^{1}(B)/\mathcal{I}$$ as (restricted) Lie algebras, where $\mathcal{I}$ is a one-dimensional (restricted) Lie ideal of $\HH^{1}(B)$.  
\end{Cor}

\begin{proof} First we notice that by Corollary \ref{no-Assump1} (1), we do not need Assumption \ref{assump}. By Corollary \ref{Kernel-delta-one-source-and-sink}, we have $V_{\AP}=0$ since we glue a source and a sink. The statement follows  from Theorem \ref{hh1-quotient-isomorphism-for-monomial-algebras} and Theorem \ref{Lie-strucutre2-of-hh1}.
\end{proof}

Recall that an exact sequence of Lie algebra homomorphisms
\begin{tikzcd}[column sep=small]
0\arrow[r] &\mathfrak{a}\arrow[r] &\mathfrak{h}\arrow[r] &\mathfrak{g}\arrow[r] &0
\end{tikzcd}
is called a \textit{central extension} of $\mathfrak{g}$ by $\mathfrak{a}$ if $[\mathfrak{a},\mathfrak{h}] = 0$, where we identify $\mathfrak{a}$ with the corresponding Lie ideal of $\mathfrak{h}$. 

\begin{Thm}\label{thm:centralext}
Let $A=kQ_A/I_A$ be a quiver algebra and let $B=kQ_B/I_B$ be a radical embedding obtained by gluing a source vertex $e_1$ and a sink vertex $e_n$ from the same block of $A$.  Then $\HH^1(B)$ is a central extension of $\HH^1(A)$ by $\mathcal{I}$, except in the Exceptional case \ref{e-case}.
\end{Thm}
\begin{proof}
By Theorem \ref{Lie-strucutre2-of-hh1} we have a short exact sequence of Lie algebras: 
\begin{equation*}
\begin{tikzcd}[column sep=small]
0\arrow[r] &\mathcal{I} \arrow[r] &\HH^1(B)\arrow[r] &\HH^1(A)\arrow[r] &0.
\end{tikzcd}
\end{equation*}
It remains to show that $\mathcal{I}$ is an ideal in the center of $\HH^1(B)$. Recall by Remark \ref{rmk:basisofI} that the one dimensional ideal  $\mathcal{I}$  has a basis given by $\sum_{\alpha\in (Q_A)_1e_1}\alpha'\pll \alpha'$, where $s(\alpha)=e_1$. Note that an element in  $\HH^1(B)$ is a linear combination of elements $\beta'\pll p'$, where $\beta'$ is an arrow in $Q_B$ and $p'$ is a path in $\B_B$. We show that $
[\sum_{\alpha\in (Q_A)_1e_1}\alpha'\pll \alpha', \beta'\pll p']=0$ for every $\beta'\pll p'$. 

For each arrow $\alpha\in (Q_A)_1$ with $s(\alpha)=e_1$, first observe that $p'$ contains $\alpha'$ if and only if $p'=p'_n\cdots p'_2\alpha'$ and $p'_i\neq \alpha'$ for $i=2, \cdots, n$. Now for such arrow $\alpha$, let us consider two cases. If $s(\beta')\neq f_1$, then $s(p')\neq f_1$, so $\beta\neq \alpha'$ and $p'$ does not contain arrow $\alpha'$. Therefore $[\sum_{\alpha\in (Q_A)_1e_1}\alpha'\pll \alpha', \beta'\pll p']=0$. On the other hand, if $s(\beta')= f_1$, then $\beta'=\alpha'_j$ and $p'=p'_n\cdots p'_2\alpha'_i$, where $s(\alpha_j)=s(\alpha_i)=e_1$. In this case, $p'_i\neq \alpha'$ for $i=2, \cdots, n$ and $\alpha\in (Q_A)_1e_1$. Hence  $[\sum_{\alpha\in (Q_A)_1e_1}\alpha'\pll \alpha', \beta'\pll p']= \alpha_j\pll p'-\alpha_j\pll p'= 0$.
\end{proof}

It is worth noting that in Theorem \ref{thm:centralext} we could have non-trivial central extensions. 
\begin{Ex1}
Consider $A$ to be the algebra with $Q_A$ given by the $n$-Kronecker quiver and let char$(k$) divide $n$. Note that in this case $Q_B$ is a quiver with one vertex and $n$-loops. Then it is easy to observe, see   Remark \ref{reduced-case}, that $\HH^1(A)\cong \pgl_n(k)$ and $\HH^1(B)\cong \gl_n(k)$ so the central extension is non-trivial.
\end{Ex1}

\section{Monomial algebras}
Recall that a finite dimensional $k$-algebra $\Lambda$ is called \emph{monomial} if it is a quotient $kQ/I$ of a path algebra, where the two-sided ideal $I$ of $kQ$ is generated by a set $Z$ of paths of length $\ge$ 2. We assume that $Z$ is \emph{minimal}, that is, no proper subpath of a path in $Z$ is again in $Z$. Clearly $Z$ is a reduced Gr{\"{o}}bner basis of $I$ under any left length-lexicographic order on $Q_{\geq 0}$. Let $\B=\B_{\Lambda}$ be the set of paths of $Q$ which do not contain any element of $Z$ as a subpath. It is clear that the (classes modulo $I$ of) elements of $\B$ form a basis of $\Lambda$. We shall denote by $\B_{n}$ the subset $Q_{n} \cap \B$ of $\B$ formed by the paths of length $n$.

For the quiver $Q$, the parallelism relation is an equivalence relation on the set of arrows $Q_1$; for $\alpha\in Q_1$, $[\alpha]$ denotes the equivalence class of $\alpha$. We denote $\bar{Q}_1$ the set of equivalence classes of parallel arrows. The quiver which has $Q_0$ as vertices and
$\bar{Q}_1$ as set of arrows, will be denoted by $\bar{Q}$. We denote by $\chi(\bar{Q})$ the \emph{first Betti number} of $\bar{Q}$ which is equal to $|\bar{Q}_1|-|\bar{Q}_0|+c_{\bar{Q}}$, where $c_{\bar{Q}}$ is the number of connected components of $\bar{Q}$.
\subsection{A direct sum decomposition of $\HH^1$}\label{mon-sec}
In this subsection, we will show that our Theorem \ref{thm:centralext}  can be strengthened to Corollary \ref{product-decomp-source-sink} for monomial algebras. More precisely, we will show that when we glue a source and a sink in a monomial algebra $A$, the central extension is actually a trivial extension, that is, we have a direct sum of Lie algebras.

\begin{comment}
	We also need recall some notations and results in \cite{STR}. Recall that a finite dimensional $k$-algebra $\Lambda$ is called \emph{monomial} if it is a quotient $kQ/I$ of a path algebra where the two-sided ideal $I$ of $kQ$ is generated by a set $Z$ of paths of length $\ge$ 2. We assume that $Z$ is \emph{minimal}, that is, no proper subpath of a path in $Z$ is again in $Z$. Clearly $Z$ is a reduced Gr{\"{o}}bner basis of $I$ under any left length-lexicographic order on $Q_{\geq 0}$. Let $\B=\B_{\Lambda}$ be the set of paths of $Q$ which do not contain any element of $Z$ as a subpath. It is clear that the (classes modulo $I$ of) elements of $\B$ form a basis of $\Lambda$. We shall denote by $\B_{n}$ the subset $Q_{n} \cap \B$ of $\B$ formed by the paths of length $n$.
\end{comment}

According to \cite[Section 4]{STR}, the Lie algebra $\mathrm{HH^1}
(\Lambda)$ of a monomial algebra $\Lambda=kQ/\langle Z \rangle$ has a 
natural grading. Indeed, if $a\pll \gamma \in Q_1\pll  
\B_n$ and $b\pll \epsilon \in Q_1\pll  \B_m$, then 
the Lie bracket defined in Theorem \ref{Lie-bracket} shows that $\left[ a 
\pll  \gamma,b \pll  \epsilon \right] \in k(Q_1\pll  
\B_{n+m-1})$. Thus, we have a grading on the Lie algebra 
$k(Q_1\pll  \B)=\oplus_{i\in \mathbb{N}} k(Q_1\pll  \B_i)$ by considering that the elements of $k(Q_1\pll  
\B_i)$ have degree $i-1$ for all $i\in \mathbb{N}$. It is clear 
that the Lie subalgebra $\Ker(\delta^{1})$ of $k(Q_1\pll  
\B)$ preserves this grading and that $\Im(\delta^{0})$ 
is a graded ideal, which induces a grading on the Lie algebra 
$\mathrm{HH^1}(\Lambda)\simeq \Ker(\delta^{1})/\Im
(\delta^{0})$. More precisely, if we set $$L_{-1}:=k(Q_1\pll  Q_0)\cap 
\Ker(\delta^{1}),$$ $$L_0:=(k(Q_1\pll Q_1)\cap \Ker
(\delta^{1}))/\langle \delta^0(e\pll e)\ |\ e\in Q_0 \rangle \ and$$  
$$L_i:=(k(Q_1\pll \B_{i+1})\cap \Ker
(\delta^{1}))/\langle \delta^0(e\pll p)\ |\ e\pll p \in Q_1\pll  
\B_i \rangle$$
for all $i\geq 1$, then $\mathrm{HH^1}(\Lambda)=\bigoplus_{i\geq -1}L_i$ and $\left[ L_i,L_j \right] \subseteq L_{i+j}$ for all $i,j\geq -1$, where $L_{-2}=0$.

\begin{Rem}\label{graded-HH-one}
Note that if the characteristic of the field $k$ is equal to 0, then $L_{-1}=0$  by Proposition 4.2 in \cite{STR}. It follows that $\bigoplus_{i\geq 1}L_i$ is a solvable Lie ideal of $\mathrm{HH^1}(\Lambda)=\bigoplus_{i\geq 0}L_i$ since $\mathrm{HH^1}(\Lambda)$ is finite dimensional. It is also obvious that $L_0$ is a Lie subalgebra of $\mathrm{HH^1}(\Lambda)$. 
\end{Rem}

In order to ensure each  $L_0^{[\alpha]}$ (in the Lie algebra decomposition $(\dag)$ of $L_0$ below) to be a Lie ideal, we need to use the following variation of \cite[Proposition 4.7]{STR}. In the following, when we say an undirected path in a quiver $Q$, we mean a walk consisting of concatenating arrows and formal inverses of arrows.  In particular, an oriented cycle is a closed directed walk and an undirected cycle is  a closed walk which is not directed. 
\begin{Def}\label{def:simuc}
A \emph{simple undirected (resp. oriented) cycle} is an undirected  (resp. oriented)  cycle in which only the first and last vertices are equal.
\end{Def} 

\begin{Lem}\label{basis-L_0} There is a basis $\B_{L_0}$ of $L_0$ whose elements are given by the union of the following sets:
	\begin{itemize}
		\item [$(i)$] all the elements $a\pll b \in L_0$ such that $a\neq b$;
		\item [$(ii)$] for every class $[\alpha]$ of parallel arrows  we fix a total order on the set of arrows so that $\alpha_i\prec\alpha_j$ for $i<j$. We write $[\alpha]=\{\alpha_1,\alpha_2,\cdots,\alpha_n \} \in \bar{Q}_1$. We pick all the elements $\alpha_i \pll  \alpha_i-\alpha_n\pll  \alpha_n \in L_0$ such that $i< n$; 
		\item [$(iii)$] using the same notation  in $(ii)$, for each (oriented or undirected) cycle in $\bar{Q}$, choose one class of parallel arrows $[\alpha]=\{\alpha_1,\alpha_2,\cdots,\alpha_n \}$ in this cycle and take the largest element in the parallelism class, that is, $\alpha_n\pll  \alpha_n$. Note that there are $\chi(\bar{Q})$ linearly independent elements in $(iii)$.
				\end{itemize}
\end{Lem}

For each class of parallel arrows $[\alpha] \in \bar{Q}_1$ we denote by $L_0^{[\alpha]}$ the Lie ideal of $L_0$ generated by the elements of the form $\alpha_i\pll \alpha_j$ and $\alpha_i \pll  \alpha_i-\alpha_n\pll \alpha_n$ in $ \B_{L_0}$, where $[\alpha]=\{\alpha_1,\alpha_2,\cdots,\alpha_n\}$ and $1\leq i,j\leq n$. Obviously the Lie algebra $L_0$ is the direct sum of these Lie algebras:
$$L_0=\oplus_{[\alpha]\in \bar{Q}} L_0^{[\alpha]},\quad \quad (\dag)$$
where this decomposition depends on the basis $\B_{L_0}$ and $L_0^{[\alpha]}$ may be equal to zero for some $[\alpha]$.

\begin{Rem}\label{minimal-generator-of-ideal-I} Let $A$ be a monomial algebra and let $B$ be a radical embedding obtained by gluing two idempotents in $A$. Then $B$ is also a monomial algebra, hence both $\HH^1(A)$ and $\HH^1(B)$ have a canonical grading. Note that the one-dimensional ideal $\mathcal{I}:=Y/\Im(\delta_B^0)$ of $\HH^1(B)$ in Theorem \ref{Lie-strucutre2-of-hh1} is generated by $\varphi_1(\delta_A^0(e_1\pll e_1))=\varphi_1(\sum_{[\alpha]\in (\bar{Q}_A)_1e_1}\mathcal{I}_{[\alpha]}-\sum_{[\alpha]\in e_1(\bar{Q}_A)_1}\mathcal{I}_{[\alpha]})$, where $\mathcal{I}_{[\alpha]}:=\sum_{i=1}^{m}\alpha_i\pll \alpha_i$ for $[\alpha]=\{\alpha_1,\alpha_2,\cdots,\alpha_m \}$.
\end{Rem}

We can rewrite the generator $\varphi_1(\delta_A^0(e_1\pll e_1))$ of $\mathcal{I}$ after introducing the following definition.

\begin{Def}\label{def-Delta}
Let $Q^c_A$ and $Q^d_A$ be two sub-quivers of $Q_A$ as follows:
  \begin{itemize}
	\item the arrows of $Q^c_A$ consist of all the arrows in $(Q_A)_1$ such that themselves or their formal inverses lie in a walk from $e_1$ to $e_n$ in $Q_A$ that passes through $e_1$ only once;
	\item the arrows of $Q^d_A$ are the arrows of $Q_A$ which are not in $Q^c_A$.
  \end{itemize} 
We also define the corresponding sub-quivers $Q^c_B$ and $Q^d_B$ via the map $\varphi$ in Section \ref{pre}.
	
Denote by $\Delta:=(\bar{Q}^c_A)_1e_1$ the subset of $(\bar{Q}_A)_1e_1$ consisting of the equivalence classes of parallel arrows $[\alpha]$ starting from $e_1$ in $\bar{Q}^c_A$.
\end{Def}

For a concrete example for $\Delta$, see Example \ref{eg-I-summand-of-HH}.

\begin{Lem}\label{fact-Delta}
	Let $A=kQ_A/I_A$ be a monomial algebra and let $B=kQ_B/I_B$ be a radical embedding obtained by gluing a source vertex $e_1$ and a sink vertex $e_n$ of $A$. If $e_1$ and $e_n$ are from the same block of $A$, then the one-dimensional Lie ideal $\mathcal{I}$ in Corollary \ref{quotient-lie-algebra-source-and-sink} is generated by $\varphi_1(\sum_{[\alpha]\in \Delta}\mathcal{I}_{[\alpha]})$ (modulo an element in $\Im(\delta_B^0)$).
\end{Lem}

\begin{proof}
	Since $e_1$ is a source vertex, then Remark \ref{minimal-generator-of-ideal-I} yields that 
	\begin{equation*}
	\begin{split}
	\varphi_1(\delta_A^0(e_1\pll e_1))&=\varphi_1(\sum_{[\alpha]\in (\bar{Q}_A)_1e_1}\mathcal{I}_{[\alpha]})= \varphi_1(\sum_{[\alpha]\in (\bar{Q}^c_A)_1e_1}\mathcal{I}_{[\alpha]})+  \varphi_1(\sum_{[\alpha]\in (\bar{Q}^d_A)_1e_1}\mathcal{I}_{[\alpha]})\\
	&=  \varphi_1(\sum_{[\alpha]\in\Delta}\mathcal{I}_{[\alpha]})+  \varphi_1(\sum_{[\alpha]\in (\bar{Q}^d_A)_1e_1}\mathcal{I}_{[\alpha]}).
	\end{split}
	\end{equation*}
	Note that $Q_A^c$ and $Q_A^d$ can be obtained from $Q_A$ by splitting in $e_1$ since Definition \ref{def-Delta} shows that $Q_A^c$ and $Q_A^d$ are disjoint and they only share the vertex $e_1$ when $Q^d_A$ is not empty. By combining this with the fact that $e_1$ is a source vertex, we deduce that $\delta_A^0(\sum_{i=1}^{n}e_i\pll e_i)=0$ if and only if $\delta_A^0(\sum_{e_i\in Q_A^c}e_i\pll e_i)=0$ and $\delta_A^0(\sum_{e_i\in Q_A^d}e_i\pll e_i)=0$, whence $$\varphi_1(\sum_{[\alpha]\in (\bar{Q}^d_A)_1e_1}\mathcal{I}_{[\alpha]})=-\varphi_1(\sum_{\substack{e_i \in (Q_A^d)_0,\\ e_i\neq e_1}} \delta^0_A(e_i \pll  e_i))= -\sum_{\substack{f_i \in (Q_B^d)_0,\\ f_i\neq f_1}} \delta^0_B(f_i \pll  f_i) \in \Im(\delta_B^0).$$
\end{proof}
Now we can give the main result in this subsection.

\begin{Cor}\label{product-decomp-source-sink}
	Let $A=kQ_A/I_A$ be a monomial algebra and let $B=kQ_B/I_B$ be a radical embedding obtained by gluing a source vertex $e_1$ and a sink vertex $e_n$ of $A$. If $e_1$ and $e_n$ are from the same block of $A$ and $\mathrm{char}(k)=0$, then $$\HH^{1}(B)\simeq \HH^{1}(A) \oplus \mathcal{I}\simeq \HH^{1}(A) \oplus k$$ as Lie algebras. 
\end{Cor}

\begin{proof}
We claim that it is enough to show that we have a decomposition as vector spaces $L_0=\mathcal{I}\oplus G$, where $G$ is a Lie subalgebra of $L_0$. Indeed, if this is the case, by the grading on $\HH^{1}(B)$ we have a decomposition as vector spaces: $$\HH^{1}(B)=L_0\oplus \bigoplus_{i\geq 1}L_i=(\mathcal{I}\oplus G)\oplus \bigoplus_{i\geq 1}L_i=\mathcal{I}\oplus(G\oplus \bigoplus_{i\geq 1}L_i)=:\mathcal{I}\oplus L.$$ Note that  $L$ is a Lie subalgebra of $\HH^{1}(B)$ since $G$ is a Lie subalgebra  and $\bigoplus_{i\geq 1}L_i$ is a Lie ideal of $\HH^{1}(B)$. In addition, by Theorem \ref{thm:centralext} the ideal $\mathcal{I}$ is in the center of $\HH^1(B)$, hence $L$ is a Lie ideal of $\HH^1(B)$.  Therefore we have a direct sum decomposition as Lie algebras: $\HH^1(B)=\mathcal{I}\oplus L$. Since $\HH^{1}(A)\simeq \HH^{1}(B)/\mathcal{I}$ as Lie algebras, then $L\simeq \HH^{1}(A)$ as Lie algebras.  Therefore,   $\HH^{1}(B)=\mathcal{I}\oplus \HH^{1}(A)$ as Lie algebras.

We show that $\mathcal{I}$ is a direct summand of $L_0$ as vector spaces. From now on, we fix the `minimal' generator $\varphi_1(\sum_{[\alpha]\in \Delta}\mathcal{I}_{[\alpha]})$ of the Lie ideal $\mathcal{I}$ which is given by Lemma \ref{fact-Delta}. We sketch the proof in the case that $\Delta$ only contains two equivalence classes of parallel arrows (cf. Definition \ref{def-Delta}), namely $\Delta=\{[\alpha],[\beta] \}$, where $[\alpha]=\{\alpha_1,\cdots,\alpha_m \}$ and $[\beta]=\{\beta_1,\cdots,\beta_t \}$. Then $\mathcal{I}=\langle \varphi_1(\mathcal{I}_{[\alpha]}+\mathcal{I}_{[\beta]}) \rangle=\langle \sum_{i=1}^{m}\alpha'_i\pll  \alpha'_i+\sum_{j=1}^{t}\beta'_j\pll  \beta'_j \rangle$. Since $L_0^{[\alpha]}=\langle \alpha'_i\pll \alpha'_j,\alpha'_l\pll \alpha'_l \ |\ 1\leq l\leq m,1\leq i\neq j\leq m,\alpha_i\pll \alpha_j \in \Ker(\delta_A^1) \rangle$ and $L_0^{[\beta]}=\langle \beta'_i\pll \beta'_j,\beta'_l\pll \beta'_l \ |\ 1\leq l\leq t,1\leq i\neq j\leq t,\beta_i\pll \beta_j \in \Ker(\delta_A^1) \rangle$, it is easy to see that $\mathcal{I}$ is a summand of $\oplus_{[\alpha] \in \Delta}L_0^{[\alpha]}=L_0^{[\alpha]}\oplus L_0^{[\beta]}$, hence it is a summand of $L_0$. In fact, we have vector space decompositions:
	\begin{equation*}
		\begin{split}
		    L_0^{[\alpha]} & =\langle \sum_{i=1}^{m}\alpha'_i\pll  \alpha'_i\rangle \oplus \langle \alpha'_i\pll \alpha'_j,\alpha'_l\pll \alpha'_l-\alpha'_m\pll \alpha'_m \ |\ 1\leq l\leq m-1,1\leq i\neq j\leq m,\alpha_i\pll \alpha_j \in \Ker(\delta_A^1) \rangle \\ & =:\langle \varphi_1(I_{[\alpha]}) \rangle \oplus J_1,
		\end{split}
	\end{equation*}	
	\begin{equation*}
		\begin{split}
		L_0^{[\beta]} &=\langle \sum_{i=1}^{t}\beta'_i\pll  \beta'_i\rangle \oplus \langle \beta'_i\pll \beta'_j,\beta'_l\pll \beta'_l-\beta'_t\pll \beta'_t \ |\ 1\leq l\leq t-1,1\leq i\neq j\leq t,\beta_i\pll \beta_j \in \Ker(\delta_A^1) \rangle\\&=:\langle \varphi_1(I_{[\beta]}) \rangle \oplus J_2.
		\end{split}
	\end{equation*}
	
As a consequence,  there are vector space decompositions
    \begin{equation*}
    	\begin{split}
    	L_0^{[\alpha]}\oplus L_0^{[\beta]} & =(\langle \varphi_1(\mathcal{I}_{[\alpha]}) \rangle \oplus \langle \varphi_1(\mathcal{I}_{[\beta]}) \rangle) \oplus J_1 \oplus J_2= (\mathcal{I}\oplus \langle \varphi_1(\mathcal{I}_{[\alpha]}) \rangle) \oplus J_1 \oplus J_2 \\ & = \mathcal{I}\oplus (\langle \varphi_1(\mathcal{I}_{[\alpha]}) \rangle \oplus J_1 \oplus J_2)=:\mathcal{I}\oplus J.
    	\end{split}
    \end{equation*}
   It follows from the definition of the Lie bracket in Theorem \ref{Lie-bracket} that  $J$ is a subalgebra of $L_0$. As a consequence, $L_0=\oplus_{[\alpha] \in (\bar{Q}_B)_1}L_0^{[\alpha]}=\mathcal{I}\oplus G$ as vector spaces, where $G:=J\oplus \bigoplus_{\alpha\in(\bar{Q}_B)_1\backslash\Delta}L_0^{[\alpha]}$. Moreover, by combining the direct sum decomposition $(\dag)$ with $J$ being a subalgebra of $L_0$, we have that $G$ is also a Lie subalgebra of $L_0$.
\end{proof}

\subsection{Radical square zero algebras}\label{rsz-sec}
\ 

We now apply our main results in Section \ref{HH1sec} to a subclass of monomial algebras: radical square zero algebras. An application of these results can be found in Subsection \ref{inverse-gluing-sec}. Throughout this subsection, we let $A=kQ_A/I_A$ be a radical square zero algebra and let $B=kQ_B/I_B$ be a radical embedding obtained by gluing two idempotents $e_1$ and $e_n$ of $A$.

\begin{Cor}\label{hh1-radical-square-zero-case}
Let $A$ be a radical square zero algebra and let $B$ be a radical embedding obtained by gluing two idempotents $e_1$ and $e_n$ of $A$. If $\mathrm{char}(k)\neq 2$, then we have $$\dim \HH^{1}(A)=\dim \HH^{1}(B)-1-\ap-c_B+c_A.$$
In particular, if we glue $e_1$ and $e_n$ from the same block of $A$, then
	$$\dim \HH^{1}(A)=\dim \HH^{1}(B)-1-\ap;$$ if we glue $e_1$ and $e_n$ from two different blocks of $A$, then
$$\dim \HH^{1}(A)=\dim \HH^{1}(B)-\ap$$ and  $\HH^{1}(A)\simeq \HH^{1}(B)$ as Lie algebras if we exclude the case that there are loops both at $e_1$ and $e_n$. 
\end{Cor}

\begin{proof} For radical square zero algebras, there are no extendable paths and Assumption \ref{assump} is equivalent to the condition that $\mathrm{char}(k)\neq 2$. The dimension formulas follow immediately from Theorem \ref{hh1-quotient-isomorphism-for-monomial-algebras} and Theorem \ref{Lie-strucutre-of-hh1}. Moreover, if we glue $e_1$ and $e_n$ from two different blocks of $A$ and exclude the case that there are loops at $e_1$ and at $e_n$ simultaneously, then $V_{\AP}=0$ by Remark \ref{spp-rad-square-zero}. Since $V_{\E}$ is a subspace of $V_{\AP}$, then $V_{\E}=0$. And by Theorem \ref{Lie-strucutre2-of-hh1} we have $\HH^{1}(A)\simeq \HH^{1}(B)$ as Lie algebras.
\end{proof}

Moreover, it is easy to see that if one of the following conditions holds, then the results in Corollary \ref{hh1-radical-square-zero-case} still hold in the case $\mathrm{char}(k)=2$ by Corollary \ref{no-Assump1} (2):
	\begin{itemize}
		\item[$(i)$] glue $e_1$ and $e_n$ from the same block of $A$;
		\item[$(ii)$] glue $e_1\in A_1$ and $e_n\in A_2$ from the different blocks of $A$ such that both $A_1$ and $A_2$ are not isomorphic to $k[x]/(x^2)$.
\end{itemize}

\begin{Rem} \label{exceptional-case-different-blocks} Let $A$ and $B$ as above and let $A_1$ and $A_2$ be two different blocks of $A$. Suppose $e_1\in A_1$ and $e_n\in A_2$.
\begin{itemize}
  \item[$(1)$] If there are loops at $e_1$ or at $e_n$, then in general $\HH^{1}(A)$ is not isomorphic to $\HH^{1}(B)$ and the difference between the dimensions of $\HH^1(A)$ and $\HH^1(B)$ can be arbitrarily large, see Example \ref{diffblockarblarg}.
  \item[$(2)$] If $\mathrm{char}(k)=2$ and exactly one of $A_1, A_2$ is isomorphic to $k[x]/(x^2)$, then $$\dim \HH^{1}(A)=\dim \HH^{1}(B)+1.$$
  \item[$(3)$] If $\chr(k)=2$ and $A_1\simeq k[x]/(x^2)\simeq A_2$, then $\HH^1(A)\not\simeq \HH^1(B)$ as Lie algebras although $\HH^1(A)\simeq \HH^1(B)$ as $k$-vector spaces, see Example \ref{Diffblockinc}.
\end{itemize}
\end{Rem}

\begin{Cor}\label{hh1-radical-square-zero-case-structure}
Let $A$ be a radical square zero algebra and let $B$ be a radical embedding obtained by gluing two idempotents $e_1$ and $e_n$  from the same block of $A$. If $V_{\AP}=0$ and $\mathrm{char}(k)=0$, then we have a Lie algebra isomorphism $$\HH^1(B)\simeq \HH^1(A)\oplus k.$$
\end{Cor}

\begin{proof} We use the notation in Theorem \ref{Lie-strucutre-of-hh1}. Since $V_{\AP}=0$, we have a Lie algebra epimorphism from $\HH^1(B)$ to $\Ker(\delta_B^1)/Y\simeq \HH^1(A)$ with one-dimensional kernel $I:=Y/\Im(\delta_B^0)$, where $\Ker(\delta_B^1)=\varphi_1(\Ker(\delta_A^1))$ and $Y$ is a Lie ideal of $\Ker(\delta_B^1)$. Also note that this epimorphism and equality do not depend on Assumption \ref{assump} since we glue $e_1$ and $e_n$ from the same block, cf. Corollary \ref{no-Assump1} (2). Since $V_{\AP}=0$, then $\dim \HH^{1}(A)=\dim \HH^{1}(B)-1$. Note that by gluing $e_1$ and $e_n$ from the same
block we have $\chi(\bar{Q}_B)=\chi(\bar{Q}_A)+1$. By Theorem 2.9 in \cite{S-F} (see also Theorem \ref{S-F-1}) there is an injective Lie algebra homomorphism:
$$\HH^1(A)\simeq \oplus_{\alpha \in S} \sl_{|
\alpha|}(k) \oplus k^{\chi(\bar{Q}_A)}\to  \oplus_{\alpha \in S} \sl_{|\alpha|}(k) \oplus k^{\chi(\bar{Q}_A)}\oplus k \simeq\HH^1(B).$$
Therefore it gives rise to the following Lie algebra isomorphisms: $$\HH^1(B)\simeq \HH^1(A)\oplus I\simeq \HH^1(A)\oplus k.$$
\end{proof}

\begin{Rem}\label{reduce-cycle}
	Let $A$ and $B$ as above, and suppose that $e_1$ and $e_n$ are in the same block of $A$. If we exclude the Exceptional case \ref{e-case},
then it is straightforward to check that $V_{\AP}=0$ under each of the following conditions:
		\begin{itemize}
			\item[$(i)$] $e_1$ is a source and $e_n$ is a sink;
			\item[$(ii)$] Both $e_1$ and $e_n$ are sinks such that $$\{s(\alpha)\ |\ t(\alpha)=e_1,\alpha\in (Q_A)_1\} \cap \{s(\beta)\ |\ t(\beta)=e_n,\beta\in (Q_A)_1\}=\emptyset;$$
			\item[$(iii)$] Both $e_1$ and $e_n$ are sources such that $$\{t(\alpha)\ |\ s(\alpha)=e_1,\alpha\in (Q_A)_1\} \cap \{t(\beta)\ |\ s(\beta)=e_n,\beta\in (Q_A)_1\}=\emptyset.$$
	\end{itemize}
\end{Rem}

\begin{Rem}\label{reduced-case}
Let $A$ be a radical square zero algebra having Gabriel quiver $Q$. By direct computations, we can determine the Lie algebra structure of $\HH^1(A)$ in the following well-known cases, which are recalled here for completeness:
\begin{itemize}
	\item[$(1)$] $\HH^1(A)\simeq \gl_n(k)$ if $Q$ is the quiver with one vertex and $n$ loops , except in the case $n=1$ and $\mathrm{char}(k)=2$ (for this exceptional case, see Remark \ref{exceptional-case}); The isomorphism sends $\alpha_i\pll \alpha_j$ to $E_{ji}$, where $E_{ij}$ is the matrix that has  $1$ in position $(i,j)$ and $0$ elsewhere. Note that if the characteristic of the field $k$ does not divide $n$, then $\gl_n(k) \simeq \sl_n(k)\oplus k$ as Lie algebras.
	
	\item[$(2)$] $\HH^1(A)\simeq \pgl_n(k)$ if $Q$ is the $n$-Kronecker quiver, with the convention that $1$-Kronecker quiver is the Dynkin quiver $A_2$, where $\pgl_n(k)$ is the quotient of $\gl_n(k)$ by its center $k\cdot \mathrm{Id}$. Let $e$ be the source vertex of the $n$-Kronecker quiver. Then the above isomorphism can be obtained by observing that $\Ker(\delta^1_A)\simeq \gl_n(k)$ via the isomorphism in (1). In addition, this isomorphism sends the unique generator $\sum_{s(\alpha_i)=e} \alpha_i \pll  \alpha_i$ of $\Im(\delta^0_A)$ to $\mathrm{Id}$. If the characteristic of the field $k$  does not divide $n$, then $\pgl_n(k) \simeq \sl_n(k)$.
\end{itemize}
\end{Rem}

\subsection{S\'{a}nchez-Flores' decomposition via inverse gluing} \label{inverse-gluing-sec} 
\

In this section, we provide an interpretation of S\'{a}nchez-Flores' description of the Lie algebra structure of the first Hochschild cohomology for radical square zero algebras \cite{S-F} using inverse gluing operations.

Given a quiver $Q$, denote by $S$ a complete set of representatives of the non-trivial classes on the set of arrows $Q_1$, that is, equivalence classes having at least two arrows, and for $\alpha\in S$, $|\alpha|$ denotes the number of arrows in the equivalence class $[\alpha]$ of $\alpha$. S\'{a}nchez-Flores' description of the first Hochschild cohomology for radical square zero algebras can be stated as follows.

\begin{Thm} {\rm(\cite[Theorem 2.9]{S-F})} \label{S-F-1}
Let $k$ be a field  of characteristic zero and let $A$ be an indecomposable radical square
zero algebra having Gabriel quiver $Q$. Then we have an isomorphism of  Lie algebras:
\[
    \HH^1(A)\simeq \bigoplus_{\alpha \in S} \sl_{|\alpha|}(k) \oplus k^{\chi(\bar{Q})}.
\]
\end{Thm}

Note that intuitively we can say that $\chi(\bar{Q})$ counts the number of holes in $\bar{Q}$. From this point of view we could give an interpretation of the above result by inverse gluing operations. To be more intuitive we will demonstrate our method by Example \ref{typical-example} that includes all possible cases. 

\begin{Rem}\label{gen-assump-SF-decomp}
	Note that we can slightly generalize Theorem \ref{S-F-1} to the case that $\chr(k) \nmid |\alpha|$ for each $\alpha\in S$. 
	\begin{itemize}
		\item Indeed, Example \ref{typical-example} shows that our inverse gluing technique to prove Theorem \ref{S-F-1} is based on the Lie algebra decomposition $\gl_{|\alpha|}(k) \simeq \sl_{|\alpha|}(k)\oplus k$ for $\alpha\in (Q_A)_1$. This decomposition holds if $\chr(k) \nmid |\alpha|$, see Exercise 7 in Chapter 1 in \cite{Hu}.
		\item  By reviewing Example \ref{typical-example} we can see that this decomposition only occurs in Step 1 (reduce $n$-loops) and Step 4 (reduce $m$-Kronecker). Hence it is enough to require that $\chr(k)$ does not divide $|a|$ and $|b|$ for each representative $a$ lying in an $n$-loop and each representative $b$ lying in a $m$-Kronecker with $n,m\ge 2$. 
			\end{itemize} 
\end{Rem}

\begin{Ex1}\label{typical-example}
	Let $k$ be a field with $\chr(k)\neq 2$ and let $A$ be a radical square zero algebra having Gabriel quiver $Q_A$. Note that in this case $\chi(\bar{Q}_A)=4$ and $S=\{[\alpha_1],[\beta_1]\}$.
	\begin{center}
		\begin{tikzcd}
		& j\bullet \arrow["\alpha_2", loop, distance=2em, in=55, out=125] \arrow["\alpha_1"', loop, distance=2em, in=215, out=145] \arrow[r, "\eta_1"] & \bullet i \arrow[r, "\xi_4"] \arrow[d, "\xi_1"]                                                                              & \bullet h \arrow[r, "\eta_2"]                   & \bullet g                     \\
		Q_A: & a\bullet \arrow[d, "\gamma_2", shift left=2]                                                                                                 & \bullet d \arrow[l, "\beta_1" description, shift right=2] \arrow[l, "\beta_2" description, shift left=2] \arrow[ld, "\beta_3"] & \bullet e \arrow[u, "\xi_3"] \arrow[l, "\xi_2"] & \bullet f \arrow[l, "\eta_3"] \\
		& b\bullet \arrow[u, "\gamma_1", shift left]                                                                                                   &                                                                                                                              &                                                 &
		\end{tikzcd}
	\end{center}

\textbf{Step 1 (separate and reduce loops):} We separate the loops at the vertex $j$ of $Q_A$ to get $Q_B$. The  algebra $B$ has two blocks, say $B_1$ and $B_2$.
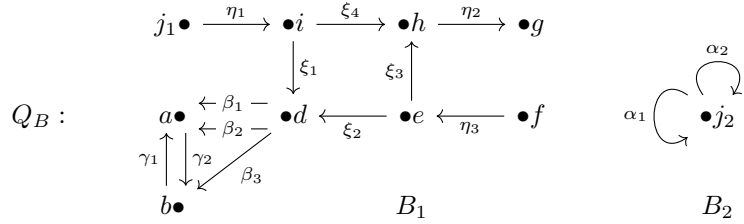
\begin{center}
	\begin{tikzcd}
	& j_1\bullet \arrow[r, "\eta_1"]               & \bullet i \arrow[r, "\xi_4"] \arrow[d, "\xi_1"]                                                                              & \bullet h \arrow[r, "\eta_2"]                   & \bullet g                     &  &                                                                                                                           \\
	Q_B: & a\bullet \arrow[d, "\gamma_2", shift left=2] & \bullet d \arrow[l, "\beta_1" description, shift right=2] \arrow[l, "\beta_2" description, shift left=2] \arrow[ld, "\beta_3"] & \bullet e \arrow[u, "\xi_3"] \arrow[l, "\xi_2"] & \bullet f \arrow[l, "\eta_3"] &  & \bullet j_2 \arrow["\alpha_1"', loop, distance=2em, in=215, out=145] \arrow["\alpha_2", loop, distance=2em, in=55, out=125] \\
	& b\bullet \arrow[u, "\gamma_1", shift left]   &                                                                                                                              & B_1                                             &                               &  & B_2
	\end{tikzcd}
\end{center}

The inverse operation is given by gluing two vertices (one of which has no loops) from two different blocks, that is, we glue $j_1\in Q_{B_1}$ and $j_2\in Q_{B_2}$. By Corollary \ref{hh1-radical-square-zero-case}, this operation does not change the dimension and the Lie structure of $\HH^1(A)$, that is, $$\HH^1(A)\simeq \HH^1(B) \simeq \HH^1(B_1)\oplus \HH^1(B_2).$$ By Remark \ref{reduced-case} (1) we obtain $\HH^1(B_2) \simeq \gl_2(k)\simeq \sl_2(k)\oplus k$, where the summand $k$ contributes $1$ to the value of $\chi(\bar{Q}_A)$. After this step, we have reduced $Q_A$ to the no loop quiver $Q_{B_1}$.

\textbf{Step 2 (reduce simple oriented $l$-cycles ($l\geq 2$)):} We reduce the  simple oriented cycle $p:=\gamma_2\gamma_1$ in $Q_{B_1}$. Choose the vertex $b$ in $p$ and split it into a source vertex $b_1$ and a sink vertex $b_2$:
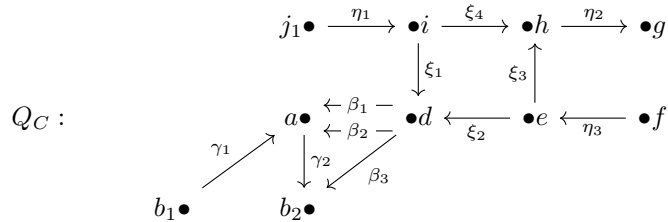
\begin{center}
	\begin{tikzcd}
	&                                   & j_1\bullet \arrow[r, "\eta_1"]             & \bullet i \arrow[r, "\xi_4"] \arrow[d, "\xi_1"]                                                                              & \bullet h \arrow[r, "\eta_2"]                   & \bullet g                     \\
	Q_C: &                                   & a\bullet \arrow[d, "\gamma_2", shift left] & \bullet d \arrow[l, "\beta_1" description, shift right=2] \arrow[l, "\beta_2" description, shift left=2] \arrow[ld, "\beta_3"] & \bullet e \arrow[u, "\xi_3"] \arrow[l, "\xi_2"] & \bullet f \arrow[l, "\eta_3"] \\
	& b_1\bullet \arrow[ru, "\gamma_1"] & b_2\bullet                                 &                                                                                                                              &                                                 &
	\end{tikzcd}
\end{center}

The inverse operation is given by gluing $b_1$ and $b_2$ from the same block. By Remark \ref{reduce-cycle} (i), by reducing $p$ from $Q_{B_1}$ we get one summand isomorphic to $k$ (cf. Corollary \ref{hh1-radical-square-zero-case-structure}), which contributes $1$ to the value of $\chi(\bar{Q}_B)=\chi(\bar{Q}_A)$. So $$\HH^1(B_1)\simeq \HH^1(C) \oplus k$$ and we have reduced $Q_{B_1}$ to the no  simple oriented cycle quiver $Q_C$.

\textbf{Step 3 (reduce simple undirected $l$-cycles ($l\geq 3$)):}  We first deal with the simple  undirected $3$-cycle $q_1:=\beta_3^{-1}\gamma_2\beta_1$  in $Q_C$. We can split $b_2$ into two sinks, say $b_3$ and $b_4$, and denote the corresponding quiver and algebra by $Q_D$ and $D$, respectively.
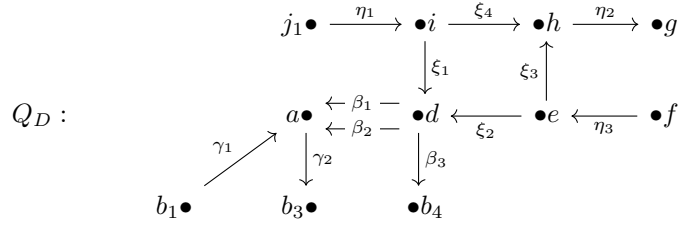
\begin{center}
	\begin{tikzcd}
	&                                   & j_1\bullet \arrow[r, "\eta_1"]             & \bullet i \arrow[r, "\xi_4"] \arrow[d, "\xi_1"]                                                                                          & \bullet h \arrow[r, "\eta_2"]                   & \bullet g                     \\
	Q_D: &                                   & a\bullet \arrow[d, "\gamma_2", shift left] & \bullet d \arrow[l, "\beta_1" description, shift right=2] \arrow[l, "\beta_2" description, shift left=2] \arrow[d, "\beta_3", shift right] & \bullet e \arrow[u, "\xi_3"] \arrow[l, "\xi_2"] & \bullet f \arrow[l, "\eta_3"] \\
	& b_1\bullet \arrow[ru, "\gamma_1"] & b_3\bullet                                 & \bullet b_4                                                                                                                              &                                                 &
	\end{tikzcd}
\end{center}

The inverse operation is given by gluing $b_3$ and $b_4$ from the same block. By Corollary \ref{hh1-radical-square-zero-case-structure} and Remark \ref{reduce-cycle}, by reducing $q_1$ from $Q_C$ we get a summand isomorphic to $k$, which again contributes $1$ to the value of $\chi(\bar{Q}_A)$. Therefore $$\HH^1(C)\simeq \HH^1(D) \oplus k.$$

We then reduce another simple  undirected cycle $q_2:=\xi_1^{-1}\xi_2\xi_3^{-1}\xi_4$. Choose the vertex $i$ in $q_2$ and split $i$ into a sink vertex $i_1$ and a source vertex $i_2$ to get $Q_E$, denote the corresponding algebra by $E$.
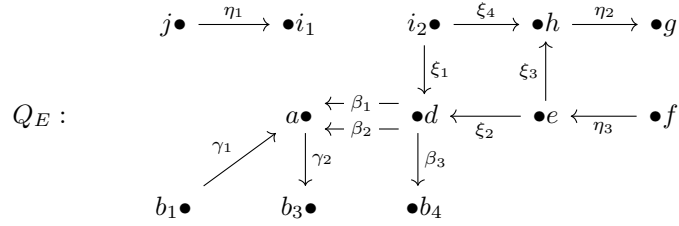
\begin{center}
	\begin{tikzcd}
	& j\bullet \arrow[r, "\eta_1"]      & \bullet i_1                                & i_2\bullet \arrow[r, "\xi_4"] \arrow[d, "\xi_1"]                                                                                         & \bullet h \arrow[r, "\eta_2"]                   & \bullet g                     \\
	Q_E: &                                   & a\bullet \arrow[d, "\gamma_2", shift left] & \bullet d \arrow[l, "\beta_1" description, shift right=2] \arrow[l, "\beta_2" description, shift left=2] \arrow[d, "\beta_3", shift right] & \bullet e \arrow[u, "\xi_3"] \arrow[l, "\xi_2"] & \bullet f \arrow[l, "\eta_3"] \\
	& b_1\bullet \arrow[ru, "\gamma_1"] & b_3\bullet                                 & \bullet b_4                                                                                                                              &                                                 &
	\end{tikzcd}
\end{center}

The inverse operation is given by gluing $i_1$ and $i_2$ from two different blocks. By Corollary \ref{hh1-radical-square-zero-case}, this operation does not change the dimension and the Lie structure of $\HH^1(D)$, that is, $$\HH^1(D) \simeq \HH^1(E).$$

Note that the above reduction produces a new simple  undirected cycle  $q'_2:=\xi_1^{-1}\xi_2\xi_3^{-1}\xi_4$ in $Q_E$. However, we can reduce $q_2'$ in $Q_E$ by splitting $i_2$ into two sources, say $i_3$ and $i_4$ (the corresponding quiver is $Q_F$).
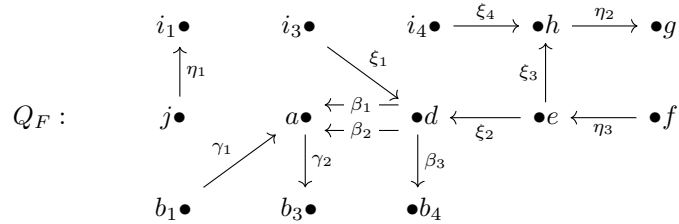
\begin{center}
	\begin{tikzcd}
	& i_1\bullet                                 & i_3\bullet \arrow[rd, "\xi_1"]             & i_4\bullet \arrow[r, "\xi_4"]                                                                                                            & \bullet h \arrow[r, "\eta_2"]                   & \bullet g                     \\
	Q_F: & j\bullet \arrow[u, "\eta_1"', shift right] & a\bullet \arrow[d, "\gamma_2", shift left] & \bullet d \arrow[l, "\beta_1" description, shift right=2] \arrow[l, "\beta_2" description, shift left=2] \arrow[d, "\beta_3", shift right] & \bullet e \arrow[u, "\xi_3"] \arrow[l, "\xi_2"] & \bullet f \arrow[l, "\eta_3"] \\
	& b_1\bullet \arrow[ru, "\gamma_1"]          & b_3\bullet                                 & \bullet b_4                                                                                                                              &                                                 &
	\end{tikzcd}
\end{center}

The inverse operation is given by gluing two sources from the same block. Again by Corollary \ref{hh1-radical-square-zero-case-structure} and Remark \ref{reduce-cycle}, we get that $$\HH^1(E)\simeq \HH^1(F) \oplus k,$$ where the summand $k$ also contributes $1$ to the value of $\chi(\bar{Q}_A)$. We have reduced to a quiver $Q_F$ that has neither oriented cycles nor simple  undirected cycles.

\textbf{Step 4 (Split into several $m$-Kronecker quivers):} Since $Q_F$ contains no simple cycles (whether oriented or undirected), we can split $Q_F$ into several quivers. Note that each of these quivers is a $m$-Kronecker quiver for some $m\geq 1$.
\begin{center}
	\begin{tikzcd}
	&                                                &                                                & i_4\bullet \arrow[r, "\xi_4"']                                                    & \bullet h_3                     &                                 \\
	& i_1\bullet                                     & i_3\bullet \arrow[r, "\xi_1"']                 & \bullet d_1                                                                       & \bullet h_2 \arrow[r, "\eta_2"] & \bullet g                       \\
	Q_G: & j\bullet \arrow[u, "\eta_1"', shift right]     & a\bullet                                       & \bullet d_2 \arrow[l, "\beta_1"', shift right] \arrow[l, "\beta_2", shift left=2] & h_1\bullet                      & \bullet e_3 \arrow[l, "\xi_3"]  \\
	& a_1\bullet                                     & \bullet a_2 \arrow[d, "\gamma_2", shift right] & \bullet d_3 \arrow[d, "\beta_3", shift right]                                     & f\bullet                        & \bullet e_2 \arrow[l, "\eta_3"] \\
	& b_1\bullet \arrow[u, "\gamma_1"', shift right] & \bullet b_3                                    & \bullet b_4                                                                       & d_4\bullet                      & \bullet e_1 \arrow[l, "\xi_2"']
	\end{tikzcd}
\end{center}

The inverse of the above operations are given by repeatedly applying three types of operations: gluing a source and a sink from different blocks, gluing two sources from different blocks, gluing two sinks from different blocks. By Corollary \ref{hh1-radical-square-zero-case}, these operations do not change the dimension and the Lie structure of $\HH^1(F)$. Therefore, $$\HH^1(F)\simeq \HH^1(G).$$
By Remark \ref{reduced-case} (2) the $\HH^1$ of a $m$-Kronecker algebra is $\sl_m(k)$, consequently $\HH^1(G) \simeq \sl_2(k)$.

We conclude that $\HH^1(B_1) \simeq \HH^1(G) \oplus k^3$, therefore $$\HH^1(A)\simeq \HH^1(B_1) \oplus \HH^1(B_2) \simeq {\sl_2(k)}^2\oplus k^4.$$ 
\end{Ex1}

\section{Center}
\label{HHsec}

Inspired by the method used in Section \ref{HH1sec}, we now study the behaviour of the centers of quiver algebras under the gluing of idempotents. Throughout this section we will denote by $Z(A)$ the center of an algebra $A$.

Recall from the proof of Lemma \ref{Image-zero-part} that there is an obvious bijection between $\Ker(\delta_{A}^{0})$ and $Z(A)$ given by 
	\begin{align*}
		\theta:\quad\Ker(\delta_{A}^{0})\quad&\xrightarrow{\sim}\quad Z(A)\\ \sum_{i=1}^{n} \lambda_i\cdot e_i \pll p_i\ &\mapsto\quad \sum_{i=1}^{n}\lambda_i p_i,
	\end{align*} where $\lambda_i\in k$ and $e_i\pll p_i\in (Q_A)_0\pll \B_{\geq 0}$. Thus, the decomposition $\Ker(\delta_A^0)=\Ker(\delta_{A,0}^0)\oplus\Ker(\delta_{A,\ge 1}^0)$ on $\Ker(\delta_A^0)$ induces a decomposition $Z(A)=Z(A)_0\oplus Z(A)_{\geq 1}$ on $Z(A)$, where $Z(A)_0:=Z(A)\cap k(Q_A)_0$ and $Z(A)_{\geq 1}:=Z(A)\cap k(\B_A)_{\geq 1}$. A similar argument holds for $B$. Therefore, in order to compare $Z(A)$ with $Z(B)$ we may reduce to comparing $\Ker(\delta_{A,0}^0)$ with $\Ker(\delta_{B,0}^0)$ and $\Ker(\delta_{A,\ge 1}^0)$ with $\Ker(\delta_{B,\ge 1}^0)$.

\begin{Lem}\label{lem:center-0}
	Let $A=kQ_A/I_A$ be a quiver algebra, and $c_A$ be the number of connected components of $Q_A$. Then $\dim Z(A)_0=c_A$ and $Z(A)_0$ is generated by the identity elements in each block of $A$.
\end{Lem}
\begin{proof}
		Note that the bijection
		$\theta:\Ker(\delta_A^0)\to Z(A)$ restricts to a bijection $\Ker(\delta_{A,0}^0)\xrightarrow{\sim} Z(A)_0$, which is given by $\sum_{i\in (Q_A)_0} \lambda_i\cdot e_i\pll e_i\mapsto \sum_{i\in (Q_A)_0} \lambda_i e_i$. We have $\dim \Ker(\delta_{A,0}^0)=c_A$ by  combining Lemma \ref{Image-zero-part} with the isomorphism $\Ker(\delta_{A,0}^0)\simeq k((Q_A)_0\pll(Q_A)_0)/\Im(\delta_{A,0}^0)$. It is easy to see that $Z(A)_0$ is generated by the identity elements in each block of $A$. 
	\end{proof}

In order to compare $\Ker(\delta_{A,\ge 1}^0)$ with $\Ker(\delta_{B,\ge 1}^0)$, we need introduce the following notion.

\begin{Def}\label{non-special-path}
	Let $A=kQ_A/I_A$ be a quiver algebra and let $B=kQ_B/I_B$ be a radical embedding obtained by gluing two idempotents $e_1$ and $e_n$ of $A$. Let  $p$ be a path between  $e_1$ and $e_n$ in $\B_A$. We call  $p$  a \emph{maximal} between $e_1$ and $e_n$ in $Q_A$ if $\delta_{B}^{0}(f_1 \pll  p')= 0$, or equivalently, if $ap\in I_A$ and $pb\in I_A$ for arbitrary $a,b \in (Q_{A})_1$.

Denote by \begin{itemize}
		\item $\M$ the set of maximal paths between $e_1$ and $e_n$ in $Q_A$,
		\item $V_{\M}:=k((Q_B)_0\pll \varphi(\M))$ the $k$-subspace of $k((Q_B)_0\pll \B_B)$ generated by the elements $f_1 \pll  p'$ for $p\in \M$, where $\varphi:\B_A\to \B_B$ is defined in Proposition \ref{parallel-paths-in-monomial-algebras}.
	\end{itemize}
\end{Def}

Note that the notion of maximal path is exactly the opposite notion of extendable path. It is clear that there are no maximal paths between $e_1$ and $e_n$ when we glue these two idempotents from different blocks. Moreover, the dimension of $V_{\M}$ equals the number of maximal paths between $e_1$ and $e_n$, that is, $\dim V_{\M}=|\M|$.
 
\begin{Lem}\label{kernel-delta}
Let $A=kQ_A/I_A$ be a quiver algebra and let $B=kQ_B/I_B$ be a radical embedding obtained by gluing two idempotents $e_1$ and $e_n$ of $A$. Then there is a decomposition as $k$-vector spaces
\[\Ker(\delta^0_{B,\geq1})=\varphi_0(\Ker(\delta^0_{A,\geq1})) \oplus V_{\M}\oplus \Ker(\delta_B^0|_{k((Q_B)_0\pll\varphi(\E))}).\]
In particular, if we glue $e_1$ and $e_n$ from the same block, then \[\dim \Ker(\delta^0_{B,\geq1})=\dim \Ker(\delta^0_{A,\geq1})+|\M|+|\E|-\e;\]
if we glue $e_1$ and $e_n$ from different blocks, then \[\dim \Ker(\delta^0_{B,\geq1})=\dim \Ker(\delta^0_{A,\geq1}).\]
\end{Lem}

\begin{proof}
Recall from Proposition \ref{parallel-paths-in-monomial-algebras} that there is a $k$-linear map $\varphi_0: k((Q_A)_0\pll  \B_A) \to k((Q_B)_0\pll  \B_B)$. A direct computation shows that $\delta_B^0(\varphi_0(e_i \pll  p))=\varphi_1(\delta_A^0(e_i \pll  p))$ for $1\leq i\leq n$ and $p\in \B_A\backslash \{e_1,e_n\}$. It follows that $\varphi_0$ induces an injective $k$-linear map from $\Ker(\delta^0_{A,\geq 1})$ to $\Ker(\delta^0_{B,\geq 1})$.

Let $x\in\Ker(\delta^0_{B,\geq 1})$ be in the complement of the subspace $\varphi_0(\Ker(\delta^0_{A,\geq 1}))$. Then we assume that $x$ is a linear combination of the elements of the form $f_1\pll  p'$ such that $p$ is a path between $e_1$ and $e_n$. If $p$ is a maximal path, then $f_1\pll  p'\in V_{\M}\subseteq \Ker(\delta^0_{B,\geq 1})$. Otherwise, $p\in \E$. Note that there may exist another extendable path $q\neq p$ such that the summands of $\delta_B^0(f_1\pll p')$ and the summands of $\delta_B^0(f_1\pll q')$ can be cancelled by each other. Consequently, $x$ can be a linear combination of the elements in $V_{\M}$ and in $\Ker(\delta_B^0|_{k((Q_B)_0\pll\varphi(\E))})$. Therefore, we have the following formula  \[\Ker(\delta^0_{B,\geq1})=\varphi_0(\Ker(\delta^0_{A,\geq1})) \oplus V_{\M}\oplus \Ker(\delta_B^0|_{k((Q_B)_0\pll\varphi(\E))}).\]
\end{proof}

\begin{Rem}
    Note that in general the generator of the $k$-vector space $\Ker(\delta_B^0|_{k((Q_B)_0\pll\varphi(\E))})$ could be the $k$-linear combination of elements of the form $f_1\pll p'$ for $p\in\E$ (see Example \ref{eg-not-equal}). However,  in monomial case we have $\Ker(\delta_B^0|_{k((Q_B)_0\pll\varphi(\E))})=0$, so Lemma \ref{kernel-delta} can be simplified.  This is because in this case Lemma \ref{summand-disjoint} holds for any paths (between $e_1$ and $e_n$) $p,q\in \B_A$ with $p\neq q$ (cf. Remark \ref{dimension-zsp} (1)). Moreover, by the same reason, we have $\e=|\E|$, that is, the dimension of $V_{\E}$ is equal to the number of extendable paths. In addition, the number of almost parallel pairs is greater than or equal to the number of extendable paths. Note that in general these statements are not true (cf. Example \ref{eg-not-equal}). 
\end{Rem}

First, we deal with the case that the algebra $A$ is indecomposable.

\begin{Prop}\label{center-from-the-same-block}
Let $A$ be an indecomposable quiver $k$-algebra and let $B$ be a radical embedding of
	$A$ obtained by gluing two idempotents $e_{1}$ and $e_{n}$
	of $A$. Then there is an algebra monomorphism $Z(A)\hookrightarrow Z(B)$. Moreover, $$\dim Z(B)=\dim Z(A)+|\M|+|\E|-\e.$$
\end{Prop}

\begin{proof} By Lemma \ref{lem:center-0}, we have $Z(A)_0=k\cdot 1_A$ and $Z(B)_0=k\cdot 1_B$. Thus, $1_A=\sum_{i=1}^n e_i=\sum_{i=1}^{n-1} f_i=1_B$ which shows that $Z(A)_0=Z(B)_0$. The second statement follows immediately from the formula $\dim Z(A)_{\geq 1}=\dim Z(B)_{\geq 1}+|\M|+|\E|-\e$ in Lemma \ref{kernel-delta}. Moreover, the proof of Lemma \ref{kernel-delta} shows that $\varphi_0$ gives an inclusion from $Z(A)_{\geq 1}=\Ker(\delta_{A,\geq 1}^{0})$ to $Z(B)_{\geq 1}=\Ker(\delta_{B,\geq 1}^{0})$. We thereby obtain an algebra monomorphism $Z(A)\hookrightarrow Z(B)$ induced from $\varphi_0$.
\end{proof}

\begin{Cor}\label{center-source-sink-same-block}
	Let $A$ be an indecomposable quiver $k$-algebra and let $B$ be a radical embedding of
	$A$ obtained by gluing a source vertex $e_{1}$ and a sink vertex $e_{n}$. Then $\varphi_0:Z(A)\hookrightarrow Z(B)$ is an isomorphism if and only if there is no path from $e_1$ to $e_n$.
\end{Cor}

\begin{proof} Note that in this case, $\E=\emptyset$ and $p$ is a maximal path between $e_1$ and $e_n$ if and only if $p$ is a path from $e_1$ to $e_n$. Thus the result follows from Proposition \ref{center-from-the-same-block}.
\end{proof}

\begin{Cor}\label{center-radical-square-zero-same-block}
	Let $A$ be a radical square zero indecomposable finite
	dimensional algebra and let $B$ be a radical embedding of
	$A$ obtained by gluing two idempotents $e_{1}$ and $e_{n}$
	of $A$. Then $\varphi_0:Z(A)\hookrightarrow Z(B)$ is isomorphism if and only if there are no arrows between $e_1$ and $e_n$ in $Q_A$.
\end{Cor}

\begin{proof} For radical square zero algebras, the set $\M$ consists of all arrows between $e_1$ and $e_n$ in $Q_A$ and there is no extendable path between $e_1$ and $e_n$ in $Q_A$.
\end{proof}

Note that Cibils has shown in \cite{C} that the dimension of the center of an indecomposable radical square zero algebra is given by $|Q_1\pll  Q_0|+1$. Indeed, by the proof of Proposition \ref{center-from-the-same-block}, we know that the basis of the center of an indecomposable radical square zero algebra is provided by the set of loops together with the unit element of the algebra.

Next we deal with the case that the algebra $A$ is not indecomposable. Without loss of generality we assume that $A$ has two blocks, say $A_1$ and $A_2$, and assume that $B$ is an algebra obtained from $A$ by gluing $e_1 \in A_1$ and $e_n \in A_2$.

\begin{Prop}\label{center-from-different-block}
Let $A$ be a quiver algebra with two blocks $A_1$ and $A_2$. Let $B$ be a radical embedding of $A$ obtained by gluing idempotents $e_1 \in A_1$ and $e_n \in A_2$. Then the radical embedding $B \to A$ restricts to a radical embedding $Z(B)\to Z(A)$. In particular, \[\dim Z(A)=\dim Z(B)+1.\]
\end{Prop}

\begin{proof}It follows from Lemma \ref{kernel-delta} that $\Ker(\delta_{A,\geq 1}^0)\simeq \Ker(\delta_{B,\geq 1}^0)$, which further implies that $Z(A)_{\geq 1}=Z(B)_{\geq 1}$. By Lemma \ref{lem:center-0}, we have that $Z(A)_0=\langle 1_{A_1},1_{A_2} \rangle$ and $Z(B)_0=\langle 1_B=1_{A_1}+1_{A_2} \rangle$. Therefore, the radical embedding from $B$ to $A$ restricts to an radical embedding from $Z(B)$ to $Z(A)$ which sends $1_B$ to $1_{A_1}+1_{A_2}$ and each element in $Z(B)_{\geq 1}$ to the same element in $Z(A)_{\geq 1}$. In particular, we have $\dim Z(A)=\dim Z(B)+1$.
\end{proof}

\section{Examples}
\label{eg}

We give some examples concerning the main  results in this paper. The first one shows that our gluing technique is useful for computing the $\HH^1$ of non-monomial algebras.

\begin{Ex1}
Let the algebra $B$ be obtained from $A$ by gluing source $e_{1}$ and sink $e_{4}$:
  \[\begin{tikzcd}
  Q_A:\hspace{-2em} & e_1\bullet \arrow[r, "\gamma"] \arrow[d, "\alpha"'] & \bullet e_3 \arrow[d, "\eta"] &  &    Q_B:\hspace{-2em}  & f_1\bullet \arrow[d, "\alpha'"', shift right=2] \arrow[r, "\gamma'", shift left=2] & \bullet f_3 \arrow[l, "\eta'", shift left=2] \\
  & e_2\bullet \arrow[r, "\beta"']                      & \bullet e_4                   &  &  & f_2\bullet \arrow[u, "\beta'"', shift right=2]                                     &             
  \end{tikzcd}\]
Where $Z_A=\{\beta\alpha-\eta\gamma\}, Z_{new}=\{\alpha'\beta',\gamma'\beta',\alpha'\eta',\gamma'\eta'\}$ and $Z_B=Z_A\cup Z_{new}$. We fix the order on $(Q_A)_1$ by $\eta\succ \gamma\succ \beta\succ \alpha$, then $\Tip(Z_A)=\{\eta\gamma\}$. It follows that $\mathcal{G}_A=\{\eta\gamma-\beta\alpha\}$ and $\mathcal{G}_B=\mathcal{G}_A\cup Z_{new}$. A direct computation based on Theorem \ref{GPC} shows that  \[\delta_A^1(\alpha\pll\alpha)=\eta\gamma\pll(\eta\gamma)^{\alpha\pll\alpha}-\eta\gamma\pll(\beta\alpha)^{\alpha\pll\alpha}=-\eta\gamma\pll\beta\alpha=\delta_A^1(\beta\pll\beta),
\delta_A^1(\gamma\pll\gamma)=\eta\gamma\pll\eta\gamma=\delta_A^1(\eta\pll\eta).\] Similarly we can compute $\delta_A^0$, $\delta_B^0$ and $\delta_B^1$. Note that in this case $V_{\E}=V_{\AP}=0$. Observe that $\beta\alpha=\eta\gamma$ in $\B_A$, we obtain that $$\Im(\delta_A^0)\simeq \langle \beta\pll\beta-\alpha\pll\alpha,\eta\pll\eta-\gamma\pll\gamma,\alpha\pll\alpha+\gamma\pll\gamma\rangle\simeq \Ker(\delta_A^1)\simeq \Ker(\delta_B^1),$$ $$\Im(\delta_B^0)\simeq \langle\beta'\pll\beta'-\alpha'\pll\alpha',\eta'\pll\eta'-\gamma'\pll\gamma' \rangle.$$ Therefore, $$\HH^1(A)\simeq\Ker(\delta_A^1)/\Im(\delta_A^0)=0,$$ $$\HH^1(B)\simeq\Ker(\delta_B^1)/\Im(\delta_B^0)\simeq \langle \alpha'\pll\alpha'+\gamma'\pll\gamma'\rangle\simeq \HH^1(A)\oplus k.$$ It is clear that $\beta\alpha=\eta\gamma$ is a maximal path between $e_1$ and $e_4$ in $Q_A$, hence $$Z(A)\simeq \Ker(\delta_A^0)\simeq\langle 1_A=\sum_{i=1}^4 e_i\pll e_i\rangle \hookrightarrow \langle 1_B=\sum_{i=1}^3f_i\pll f_i, f_1\pll \beta'\alpha' \rangle \simeq \Ker(\delta_B^0) \simeq Z(B).$$
\end{Ex1}

The second example shows a particular instance of Corollary \ref{hh1-radical-square-zero-case} in which $B$ is not obtained from $A$ by gluing a source and a sink:

\begin{Ex1}\label{rad-1} Assume $\chr(k)\neq 2$. The algebra $B$ is obtained from $A$ by gluing $e_{1}$ and $e_{3}$:

\[\begin{tikzcd}
	Q_A:\hspace{-2em} & {e_1 \bullet} & {\bullet e_2} & {\bullet e_3} & & Q_B:\hspace{-2em} & {f_2 \bullet} & {\bullet f_1}
	\arrow["{\alpha_1}"', from=1-3, to=1-2]
	\arrow["{\alpha_2}", from=1-3, to=1-4]
	\arrow["{\alpha'_1}", shift left=2, from=1-7, to=1-8]
	\arrow["{\alpha'_2}"', shift right=2, from=1-7, to=1-8]
\end{tikzcd}\]
Where $Z_A=Z_{B}=\emptyset$. Note that $A$ is hereditary, and the underlying graph of $Q_A$ is a tree, therefore $\HH^1(A)=0$. We have that $V_{\mathcal{E}_1^3}=0$ since there is no extendable path between $e_1$ and $e_3$. Moreover, observe that the set almost parallel pairs is given by \[\mathcal{AP}_1^3=\{(\alpha_1,\alpha_2),(\alpha_2,\alpha_1)\}\}.\] Thus, we obtain that \[V_{\mathcal{AP}_1^3}:=k(\mathcal{AP}_1^3)\cap \Ker(\delta_B^1)=\langle \alpha'_1\pll\alpha'_2, \alpha'_2\pll\alpha'_1 \rangle,\text{ and } \mathfrak{ap}_1^3=\dim V_{\mathcal{AP}_1^3}=2.\]  By Corollary  \ref{hh1-radical-square-zero-case} we have that $\dim \HH^1(A)=\dim \HH^1(B)-1-\mathfrak{ap}_1^3$, we can deduce that the dimension of $\HH^1(B)$ is $3$. Indeed, a direct computation shows $\HH^1(B)\cong \mathfrak{sl}_2(k)$ having a $k$-basis given by $\{\alpha'_1\pll\alpha'_1,\alpha'_1\pll\alpha'_2, \alpha'_2\pll\alpha'_1 \}$. 
\end{Ex1}
The next two examples show that Assumption \ref{assump}  in Proposition \ref{Kernel-delta-one-injective} is necessary.

\begin{Ex1}\label{example-characteristic-condition}Assume that char$(k)=2$, and that $B$ is obtained from $A$ by gluing $e_{1}$ and $e_{2}$:
	\begin{center}
		\begin{tikzcd}
Q_A:\hspace{-2em} & e_1\bullet \arrow["\alpha"', loop, distance=2em, in=125, out=55] \arrow[r, "\beta", shift left=2] & \bullet e_2 \arrow[l, "\gamma", shift left=2] &  & Q_B: & f_1\bullet \arrow["\alpha'"', loop, distance=2em, in=215, out=145] \arrow["\beta'"', loop, distance=2em, in=125, out=55] \arrow["\gamma'", loop, distance=2em, in=325, out=35]
\end{tikzcd}
	\end{center}
Where $Z_A=\{r_1=\alpha^2-\gamma\beta,r_2=\beta\alpha\gamma,\beta\gamma\}$, $Z_{new}=\{r_3=\alpha'\beta',r_4=\gamma'\alpha',(\gamma')^2,(\beta')^2\}$ and $Z_B=Z_A\cup Z_{new}$. We fix the order on $(Q_A)_1$ by $\gamma\prec \beta\prec \alpha$. Then $\mathcal{G}_A=Z_A$ and $\mathcal{G}_B=Z_B$. A direct computation shows that $\delta_A^1(\alpha\pll e_1)=2\alpha^2\pll\alpha+r_2\pll \beta\gamma=2\alpha^2\pll\alpha=0$ since char$(k)=2$. However, $\delta_B^1(\alpha'\pll f_1)=2(\alpha')^2\pll \alpha'+r_3'\pll \beta'+r_4'\pll \gamma'\ne0$,  which means that although $\alpha\pll e_1\in \Ker(\delta_A^1)$, $\varphi_1(\alpha\pll e_1)=\alpha'\pll f_1\notin \Ker(\delta_B^1)$. Hence $\varphi_1$ does not induce an injective $k$-linear map from $\Ker(\delta_A^1)$ to $\Ker(\delta_B^1)$.
\end{Ex1}

\begin{Ex1}
\label{Diffblockinc}
Let $A$ be given by two blocks $A_1$ and $A_2$ such that $A_1$ is isomorphic to $k[x]/(x^2)$ and $A_2$ is isomorphic to $k[y]/(y^2)$. Let $B$ be obtained by gluing the units of $A_1$ and $A_2$. Then $\Ker(\delta_B^1)=\HH^1(B)\simeq \mathfrak{gl}_2(k)$ and has $k$-basis
given by $\{ x \pll  x,\ x\pll  y,\ y\pll  x,\  y\pll  y\}$. However, there are two cases for $A$.

$(1)$ If $\chr(k)\neq 2$, then $\Ker(\delta_A^1)=\HH^1(A)\simeq k\oplus k$ has
$k$-basis given by $x\pll  x$ and $y\pll  y$, and there is an injective Lie algebra homomorphism $\Ker(\delta_A^1)\hookrightarrow \Ker(\delta_B^1)$.

$(2)$ If $\chr(k)=2$, then $\Ker(\delta_A^1)=\HH^1(A)$ has
$k$-basis given by $\{ x \pll  x,\ x\pll  e_1,\ y\pll  y,\  y\pll  e_2\}$. Note that in this case $\HH^1(A)$ is a direct product of two 2-dimensional  Lie algebras, hence we cannot get an injective Lie algebra homomorphism from $\Ker(\delta_A^1)$ to $\Ker(\delta_B^1)$.
\end{Ex1}

For a finite dimensional $k$-algebra $A$, it is known that every derivation of $A$ preserves its Jacoboson radical in $\chr(k)=0$, which is a special case of Assumption \ref{assump}. However, the following example shows that  Assumption  \ref{assump}  does not imply that every derivation preserves the Jacobson radical on $A$.

\begin{Ex1}
 Let $\mathrm{char}(k)=2$,  and let $Q_A$ given by:
% https://q.uiver.app/#q=WzAsMyxbMCwwLCIxIl0sWzEsMCwiMiJdLFsyLDAsIjMiXSxbMSwxLCJcXGJldGEiXSxbMCwxLCJcXGFscGhhXzIiLDJdLFsxLDIsIlxcZ2FtbWEiLDJdLFswLDAsIlxcYWxwaGFfMSJdXQ==
\[\begin{tikzcd}
1 & 2 & 3
\arrow["{\alpha_1}", from=1-1, to=1-1, loop, in=55, out=125, distance=2em]
\arrow["{\alpha_2}"', from=1-1, to=1-2]
\arrow["\beta", from=1-2, to=1-2, loop, in=55, out=125, distance=2em]
\arrow["\gamma"', from=1-2, to=1-3]
\end{tikzcd}\]
Assume $I_A$ generated by the relations $\alpha_1^3=0, \beta^2=0$. The quiver $Q_B$ is obtained by gluing $e_1$ and $e_3$. Consider the derivation $d$ represented by $\beta||e_2\in \mathrm{Ker}(\delta^1_A)$. Hence $d(\mathrm{rad } (A))$ is not contained in  $\mathrm{rad } (A)$, however, Assumption 1 holds.
\end{Ex1}

In the following example, we compute explicitly the extendable paths and the $k$-space $V_{\E}$ (resp. the almost parallel pairs and the $k$-space $V_{\AP}$) appeared in Definition \ref{special-path} and Proposition \ref{Iamge-delta-zero} (resp. in Definition \ref{elements-in-Z-spp} and Proposition \ref{Kernel-delta-one-difference}).

\begin{Ex1}\label{ex-image-kernel}
Let $B$ be obtained from $A$ by gluing $e_{1}$ and $e_{4}$:	
	\begin{center}
		\begin{tikzcd}
		&                             &                                                                                                                   &                             &              &        &  f_{3} \bullet & {} \arrow[rd, "\ddots", phantom, shift right=-1]                                                                                                                                              &    \\
		Q_{A}:\hspace{-2em} & e_{2}\bullet \arrow[r, "a"] & e_{1}\bullet \arrow[r, "\alpha_{1}", bend left] \arrow[r, "\alpha_{n}"', bend right] \arrow[r, "\vdots", phantom] & e_{4}\bullet \arrow[r, "b"] & e_{3}\bullet & Q_{B}:\hspace{-2em} & f_{2}\bullet \arrow[r, "a'"] & \bullet f_{1} \arrow["\alpha'_{1}"', loop, distance=2em, in=125, out=55] \arrow["\alpha'_{n}"', loop, distance=2em, in=35, out=325] \arrow[lu, "b'"', shift left] & {}
		\end{tikzcd}  	
	\end{center}
Where $Z_A=\emptyset$, $Z_{B}=Z_{new}=\{\alpha'_{i}\alpha'_{j}\ |\ 1\le i,j \le n\}$. Since $\alpha_ia \notin I_A$ for $1\le i \le n$, $\alpha_i$ is an extendable path from $e_1$ to $e_4$ for $1\le i \le n$, we have $\mathcal{E}_1^4=\{\alpha_i\ |\ 1\le i \le n \}$ and
\begin{equation*}
\begin{split}
 V_{\mathcal{E}_1^4} & =\langle \delta_B^0(f_1 \pll  \alpha_i')\ |\ 1\le i\le n \rangle \\ & =\langle b'\pll b'\alpha_i'-a'\pll \alpha_i'a'\ |\ 1\le i \le n \rangle.
\end{split}
\end{equation*}
Hence $\mathfrak{e}^4_1=n=\dim V_{\mathcal{E}_1^4}$. Since $a'\pll \alpha'_{i}a', b'\pll b'\alpha'_{i},\alpha_i'\pll f_1$, $a\notparallel\alpha_{i}a, b\notparallel b\alpha_{i},\alpha_i\notparallel e_1,\alpha_i\notparallel e_4$, we know that $(a,\alpha_{i}a), (b,b\alpha_{i}), (\alpha_i,e_1), (\alpha_i,e_n)$ are almost parallel pairs with respect to the gluing of $e_{1}$ and $e_{4}$ for $1\le i \le n$, and $\mathcal{AP}^4_1=\{(a,\alpha_i a), (b,b\alpha_i), (\alpha_i,e_1), (\alpha_i,e_n)\ |\ 1\le i\le n \}$. As a result we get $$k(\mathcal{AP}^4_1)=\langle a'\pll  \alpha_i'a',b'\pll  b'\alpha_i', \alpha_i'\pll f_1\ |\ 1\le i\le n \rangle,$$
\begin{equation*}
	\begin{split}
	V_{\AP} & = k(\mathcal{AP}^4_1) \cap \Ker(\delta_B^1) \\ & =\langle a'\pll  \alpha_i'a',b'\pll  b'\alpha_i'\ |\ 1\le i\le n\rangle.
	\end{split}
\end{equation*}
Hence $\mathfrak{ap}^4_1=\dim V_{\mathcal{AP}^4_1}=2n$. A direct computation shows that $\Im(\delta_{A}^{0}),\Im(\delta_{B}^{0})$ are 3-dimensional and $(n+2)$-dimensional, respectively, since $$\Im(\delta_{A}^{0})=\langle a\pll a, b\pll b,\sum\limits_{i=1}^n\alpha_{i}\pll \alpha_{i}\rangle,$$
$$\Im(\delta_{B}^{0})=\langle a'\pll a', b'\pll  b',b'\pll  b'\alpha'_{i}-a'\pll  \alpha'_{i}a'\ |\ 1\le i\le n\rangle.$$  Therefore, $$\dim \Im(\delta_{A}^{0})=\dim \Im(\delta_{B}^{0})+1-\mathfrak{e}^4_1.$$
In addition,
$$\Ker(\delta_{A}^{1})=\langle a\pll  a, b\pll  b, \alpha_{i}\pll  \alpha_{j}\ |\ 1\le i, j \le n\rangle$$ is $(n^{2}+2)$-dimensional and
$$\Ker(\delta_{B}^{1})=\langle a'\pll  a', b'\pll  b', \alpha'_{i}\pll  \alpha'_{j}, b'\pll  b'\alpha'_{i}, a'\pll  \alpha'_{i}a'\ |\ 1\le i, j \le n\rangle$$ is $(n^{2}+2n+2)$-dimensional. Hence $$\dim \Ker(\delta_{B}^{1})=\dim \Ker(\delta_{A}^{1})+\mathfrak{ap}^4_1.$$
One can verify that $\HH^1(A)$ is isomorphic to $\pgl_n(k)$ and $\HH^1(B)$ contains a subalgebra isomorphic to $\gl_n(k)$. Note also that by the notations in the proof of Theorem \ref{Lie-strucutre-of-hh1}, in this example the subspace $Y$ of $\Ker(\delta_{B}^{1})$ is equal to $\Im(\delta_{B}^{0})\oplus \langle \sum\limits_{i=1}^n\alpha_{i}'\pll \alpha_{i}'\rangle$ and $Y$ is not a Lie ideal of $\Ker(\delta_{B}^{1})$.		
\end{Ex1}

The following example shows that in non-monomial case, the dimension of $V_{\E}$ is not equal to the number of extendable paths and the number of almost parallel pairs may be smaller than the number of extendable paths.
\begin{Ex1}\label{eg-not-equal}
The algebra $B$ is obtained from $A$ by gluing $e_{1}$ and $e_{4}$:
\begin{center}
\begin{tikzcd}
Q_A:\hspace{-2em} & e_2\bullet \arrow[r, "a"] & \bullet e_1 \arrow[r, "b_1", shift left=2] \arrow[r, "b_2"', shift right=2] & \bullet e_3 \arrow[r, "c"] & \bullet e_4 &  & Q_B:\hspace{-2em} & f_2\bullet \arrow[r, "a'"] & \bullet f_1 \arrow[r, "b_1'", shift left=4] \arrow[r, "b_2'" description] & \bullet f_3 \arrow[l, "c'", shift left=4]
\end{tikzcd}
\end{center}
Where $Z_A=\{cb_1a-cb_2a\}, Z_{new}=\{b_1'c',b_2'c' \}$ and $Z_B=Z_A\cup Z_{new}$. We fix the order on $(Q_A)_1$ by $c\prec b_2\prec b_1\prec a$. Then it is clear that $\mathcal{G}_A=Z_A$ and $\mathcal{G}_B=Z_B$. Moreover, with respect to gluing $e_1$ and $e_4$, the set of extendable paths and the set of almost parallel pairs are respectively given by \[\mathcal{E}_1^4=\{cb_1,cb_2\}\quad \text{and}\quad \mathcal{AP}_1^4=\{(a,cb_1a)=(a,cb_2a)\}. \] So $|\mathcal{E}_1^4|=2$ and $\mathcal{AP}_1^4=1$. Note that $c'b_1'a'=c'b_2'a'$ in $\B_B$, it follows that \begin{equation}\label{form-4} \delta_B^0(f_1\pll c'b_1')=-a'\pll c'b_1'a'=-a'\pll c'b_2'a'=\delta_B^0(f_1\pll c'b_2'). \end{equation} As a consequence, we have \[V_{\mathcal{E}_1^4}=\langle \delta_B^0(f_1\pll c'b_1'), \delta_B^0(f_1\pll c'b_2') \rangle=\langle a'\pll c'b_1'a' \rangle.\] Therefore, the number of almost parallel pairs $|\mathcal{AP}_1^4|=1$ is less than the number of extendable paths $|\mathcal{E}_1^4|=2$. In addition, $\mathfrak{e}_1^4=\dim V_{\mathcal{E}_1^4}=1<2=|\mathcal{E}_1^4|$. 

On the other hand, the above formula (\ref{form-4}) also shows that \[\Ker(\delta_B^0|_{k((Q_B)_0\pll\varphi(\mathcal{E}_1^4))})=\langle f_1\pll c'b_1'-f_1\pll c'b_2'\rangle\] has dimension $|\mathcal{E}_1^4|-\mathfrak{e}_1^4=2-1=1$. Note that there is no maximal path between $e_1$ and $e_4$, i.e., $\mathcal{M}_1^4=\emptyset$. By Proposition \ref{center-from-the-same-block}, there is an algebra monomorphism \[ Z(A)\simeq \Ker(\delta_A^0)\simeq\langle 1_A=\sum_{i=1}^4 e_i\pll e_i\rangle \hookrightarrow \langle 1_B=\sum_{i=1}^3f_i\pll f_i, f_1\pll c'b_1'-f_1\pll c'b_2' \rangle \simeq \Ker(\delta_B^0)\simeq Z(B),\] and $\dim Z(B)=\dim Z(A)+|\mathcal{E}_1^4|-\mathfrak{e}_1^4=\dim Z(A)+1$.
\end{Ex1}

By Corollary \ref{Kernel-delta-one-source-and-sink}, if $B$ is a radical embedding obtained by gluing a source vertex $e_1$ and a sink vertex $e_n$ of $A$ (in case $\chr(k)=2$, we assume that $B$ has no block isomorphic to $k[x]/(x^2)$), then $\Ker(\delta_B^1)\simeq \Ker(\delta_A^1)$. However, the converse of Corollary \ref{Kernel-delta-one-source-and-sink} is not true in general as the following example shows.

\begin{Ex1}\label{ex-converse-coro}
Let $B$ be obtained from $A$ by gluing $e_{1}$ and $e_{4}$:
\begin{center}
			\begin{tikzcd}
			&                             &                                                                                                             &            &                            &  &      & f_3 \bullet \arrow[rd, "b'", shift right] & {} \arrow[rd, "\ddots", phantom, shift right=-1]                                                                               &    \\
			Q_A: \hspace{-2em}& e_{2}\bullet \arrow[r, "a"] & e_1\bullet \arrow[r, "\alpha_n"', bend right] \arrow[r, "\alpha_1", bend left] \arrow[r, "\vdots", phantom] & e_4\bullet & e_3\bullet \arrow[l, "b"'] &  & Q_B:\hspace{-2em} & f_2\bullet \arrow[r, "a'"]                & \bullet f_1 \arrow["\alpha_n'", loop, distance=2em, in=325, out=35] \arrow["\alpha_1'", loop, distance=2em, in=55, out=125] & {}
			\end{tikzcd}	
\end{center}
Where $Z_A=\{\alpha_ia\ |\ 1\le i\le n\}$, $Z_{new}=\{\alpha_i'b', \alpha'_{i}\alpha'_{j}\ |\ 1\le i, j \le n\}$ and $Z_B=Z_A\cup Z_{new}$. Note that although $\mathcal{AP}_1^4=\{(\alpha_i,e_1),(\alpha_i,e_4)\ |\ 1\le i\le n\},$ we have \[V_{\mathcal{AP}_1^4}:=k(\mathcal{AP}_1^4)\cap \Ker(\delta_B^1)=\langle \alpha_i'\pll f_1 \mid 1\le i\le n \rangle \cap \Ker(\delta_B^1)=0.\]
 By Proposition \ref{Kernel-delta-one-difference}  we have $\dim \Ker(\delta_B^1)=\dim \Ker(\delta_A^1)$. In fact, a direct computation shows that both
$$\Ker(\delta_A^1)=\langle a\pll a, b\pll b, \alpha_i\pll \alpha_j\ |\ 1\le i,j\le n\rangle$$
and
$$\Ker(\delta_B^1)=\langle a'\pll a', b'\pll b', \alpha'_i\pll \alpha'_j\ |\ 1\le i,j\le n\rangle$$
are $(n^2+2)$-dimensional. Hence although we do not glue a source and a sink, we still have $\Ker(\delta_B^1)\simeq \Ker(\delta_A^1)$.	
\end{Ex1}

The following example shows various types of almost parallel pairs.

\begin{Ex1}\label{example-type-special-pair} In this example we always assume that $B$ is obtained from $A$ by gluing $e_1$ and $e_n$, and that $\alpha$ is an arrow in $Q_A$ and $p$ is a path in $\B_A$. It can be proved that the almost parallel pairs ($\alpha$,$p$) arise exclusively from the following seven cases and their dual cases:
	
	$(\romannumeral1):$ $\alpha$ is a loop at $e_1$ or $e_n$, assume that \begin{tikzcd}
	e_1\bullet \arrow["\alpha", loop, distance=2em, in=55, out=125].
	\end{tikzcd} (The case that \begin{tikzcd}
	e_n\bullet \arrow["\alpha", loop, distance=2em, in=55, out=125]
	\end{tikzcd} is dual.)
	
	Case 1: $p=a_n\cdots a_1$ is an oriented cycle at $e_n$ or $p=e_n$, such as:
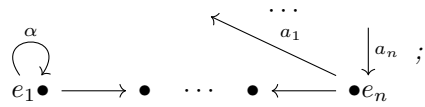
\begin{center}
	\begin{tikzcd}
	& \arrow[r, "\cdots", phantom]     & \arrow[d, "a_n"]      \\
	e_1\bullet \arrow["\alpha", loop, distance=2em, in=55, out=125] \arrow[r] & \bullet\quad\cdots\quad\bullet & \arrow[l]\bullet e_n \arrow[lu, "a_1"']
	\end{tikzcd} ;
\end{center}

	Case 2: $p=a_n\cdots a_1$ is a path between $e_1$ and $e_n$, such as:
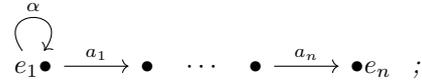
\begin{center}
	\begin{tikzcd}
	e_1\bullet \arrow["\alpha", loop, distance=2em, in=55, out=125] \arrow[r, "a_1"] & \bullet\quad\cdots\quad\bullet \arrow[r, "a_n"] & \bullet e_n
	\end{tikzcd} ;
\end{center}

	$(\romannumeral2):$ $\alpha$ is an arrow between $e_1$ and $e_n$, assume that \begin{tikzcd}
	e_1\bullet \arrow[r, "\alpha"] & \bullet e_n.
	\end{tikzcd} (The case that \begin{tikzcd}
	e_n\bullet \arrow[r, "\alpha"] & \bullet e_1
	\end{tikzcd} is dual.)
	
	Case 3: $p=a_n\cdots a_1$ is an oriented cycle at $e_1$ or $e_n$ or $p=e_1$ or $e_n$, such as:
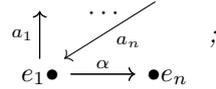
\begin{center}
	\begin{tikzcd}
	{} \arrow[r, "\cdots", phantom, shift right=3]                      & {} \arrow[ld, "a_n"] \\
	e_1\bullet \arrow[r, "\alpha"] \arrow[u] \arrow[u, "a_1"] \arrow[r] & \bullet e_n
	\end{tikzcd} ;
\end{center}

	Case 4: $p=a_n\cdots a_1$ is a path from $e_n$ to $e_1$, such as:
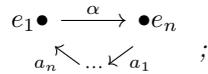
\begin{center}
	\begin{tikzcd}
	{e_1\bullet} & {\bullet e_n} \\
	{} & {}
	\arrow["\alpha", from=1-1, to=1-2]
	\arrow[""{name=0, anchor=center, inner sep=0}, "\cdots"{description}, shift right=4, draw=none, from=2-2, to=2-1]
	\arrow["{a_1}"{pos=0.3}, shift left=1, shorten >=4pt, from=1-2, to=0]
	\arrow["{a_n}"{pos=0.7}, shift left=1, shorten <=4pt, from=0, to=1-1]
	\end{tikzcd} ;
\end{center}

	$(\romannumeral3):$ Exactly one of the vertex of $\alpha$ is $e_1$ or $e_n$, assume that \begin{tikzcd}
	e_1\bullet \arrow[r, "\alpha"] \arrow[r] & \bullet .
	\end{tikzcd} (The other cases are dual.)
	
	Case 5: $p=a_n\cdots a_1$ is a path from $e_n$ to $t(\alpha)$, such as:
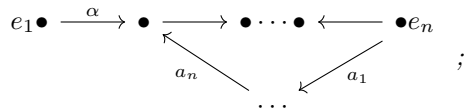
\begin{center}
	\begin{tikzcd}
	{} & {e_1\bullet} & \bullet & \bullet \cdots \bullet & {\bullet e_n} \\
	&&& \cdots
	\arrow["\alpha", from=1-2, to=1-3]
	\arrow[from=1-3, to=1-4]
	\arrow[from=1-5, to=1-4]
	\arrow["{a_1}", from=1-5, to=2-4]
	\arrow["{a_n}", from=2-4, to=1-3]
	\end{tikzcd} ;
\end{center}

	Case 6: $p=\alpha p_1$, where $p_1=a_n \cdots a_1$ is a path from $e_n$ to $e_1$, such as:
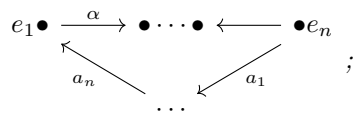
\begin{center}
	\begin{tikzcd}
	{e_1\bullet} & {\bullet \cdots \bullet} & {\bullet e_n} \\
	& \cdots
	\arrow["\alpha", from=1-1, to=1-2]
	\arrow[from=1-3, to=1-2]
	\arrow["{a_1}", from=1-3, to=2-2]
	\arrow["{a_n}", from=2-2, to=1-1]
	\end{tikzcd};
\end{center}

	Case 7: $p=p_2\alpha p_1$, where $p_1=a_n \cdots a_1$ is a path from $e_n$ to $e_1$ and $p_2=b_m\cdots b_1$ is a cycle at $t(\alpha)$, such as:
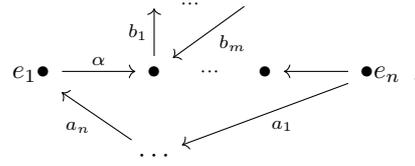
\begin{center}
	\begin{tikzcd}
	{} & {} & {} \\
	{e_1\bullet} & {\bullet } & \bullet & {\bullet e_n} \\
	& \cdots
	\arrow["\alpha", from=2-1, to=2-2]
	\arrow["{a_n}", from=3-2, to=2-1]
	\arrow["\cdots"{description, pos=0.7}, shift right=2, draw=none, from=1-1, to=1-3]
	\arrow["\cdots"{description}, draw=none, from=2-2, to=2-3]
	\arrow[from=2-4, to=2-3]
	\arrow["{b_1}", from=2-2, to=1-2]
	\arrow["{b_m}", shorten <=4pt, from=1-3, to=2-2]
	\arrow["{a_1}", from=2-4, to=3-2]
	\end{tikzcd}.
\end{center}

	After giving relations in specific examples, we can show that the almost parallel pair $(\alpha,p)$ in each of the above cases can appear. Indeed, the following example covers all the above 7 cases:

	
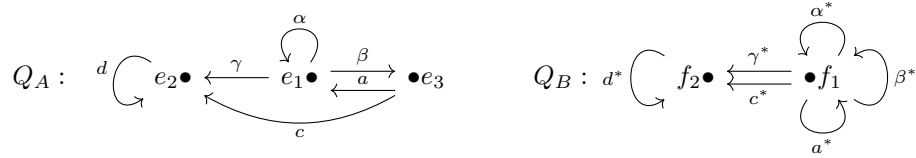
\begin{center}
		\begin{tikzcd}
		Q_A: & e_2\bullet \arrow["d"', loop, distance=2em, in=215, out=145] & e_1\bullet \arrow[r, "\beta", shift left] \arrow["\alpha"', loop, distance=2em, in=125, out=55] \arrow[l, "\gamma"'] & \bullet e_3 \arrow[l, "a"', shift left=2] \arrow[ll, "c", bend left] & Q_B: & f_2\bullet \arrow["d'"', loop, distance=2em, in=215, out=145] & \bullet f_1 \arrow["\alpha'"', loop, distance=2em, in=125, out=55] \arrow["\beta'"', loop, distance=2em, in=35, out=325] \arrow["a'"', loop, distance=2em, in=305, out=235] \arrow[l, "\gamma'"', shift right] \arrow[l, "c'", shift left]
		\end{tikzcd}
	\end{center}
Where $Z_A$ consists of all paths in $Q_A$ of length $3$ except $d\gamma a$,  $Z_B=Z_A\cup Z_{new}$ where $Z_{new}=\{a'\alpha', c'\alpha',\alpha'\beta',(\beta')^2,\gamma'\beta',(a')^2,c'a'\}$. We list all almost parallel pairs ($\alpha,p$) for each case as follows:

Case 1: ($\alpha,\beta a$), ($\alpha,e_3$);

Case 2: ($\alpha,a$), ($\alpha,\beta$), ($\alpha,\beta\alpha$), ($\alpha,\alpha a$);

Case 3: ($\beta,\alpha$), ($\beta,a\beta$), ($\beta,\beta a$), ($\beta,e_1$), ($\beta,e_3$) ($a,\alpha$), ($a,a\beta$), ($a,\beta a$), ($a,e_1$), ($a, e_3$);

Case 4: ($\beta,a$), ($a,\beta$);

Case 5: ($\gamma,c$), ($\gamma,dc$), ($c,\gamma \alpha$), ($c,d\gamma$);

Case 6: ($\gamma,\gamma a$), ($c,c\beta$);

Case 7: ($\gamma,d\gamma a$).

By check one by one, we have $\mathcal{AP}_1^3$ is the set consisting of these 25 almost parallel pairs and $k(\mathcal{AP}^3_1) = \langle \alpha' \pll  p' \ |\ (\alpha,p) \in \mathcal{AP}^3_1 \rangle $ and therefore
\begin{equation*}
	\begin{split}
		V_{\mathcal{AP}_1^3} & =k(\mathcal{AP}^3_1) \cap \Ker(\delta_B^1) \\ & =\langle a'\pll  \beta'a', \alpha'\pll  \beta'\alpha', \alpha'\pll \alpha'a', \beta'\pll  a'\beta', \beta'\pll  \beta' a', a'\pll a'\beta', \\ & a'\pll  \beta'a',\gamma' \pll  d'c', c'\pll  \gamma'\alpha', c'\pll  d'\gamma', \gamma'\pll \gamma'a', c'\pll c'\beta', \gamma'\pll d'\gamma'a' \rangle.
	\end{split}
\end{equation*}
Hence $\mathfrak{ap}^3_1=13$. Note also that the extendable paths in this example are $\beta$ and $a$, so $\mathfrak{e}^3_1=2$.
\end{Ex1}

It worth to mention that, although the $k$-space $k(\AP)$ is generated by the elements of the form $\alpha' \pll  p'$ (where $\alpha$ is an arrow and $p$ is a path), an element in $V_{\AP}$ is usually a $k$-linear combination of such elements.

\begin{Ex1}\label{linear-combination-zspp} Let $B$ be obtained from $A$ by gluing $e_1$ and $e_5$:
	
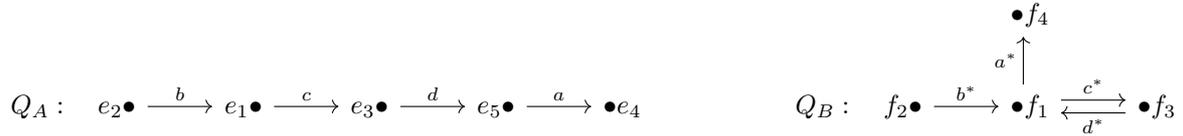
\begin{center}
		\begin{tikzcd}
		&&&&&&&&& {\bullet f_4} \\
		{Q_A:} \hspace{-3em} & {e_2\bullet} & {e_1\bullet} & {e_3 \bullet} & {e_5 \bullet} & {\bullet e_4} && {Q_B:} \hspace{-3em} & {f_2\bullet} & {\bullet f_1} & {\bullet f_3}
		\arrow["b", from=2-2, to=2-3]
		\arrow["c", from=2-3, to=2-4]
		\arrow["d", from=2-4, to=2-5]
		\arrow["a", from=2-5, to=2-6]
		\arrow["{b'}", from=2-9, to=2-10]
		\arrow["{a'}", shift left=1, from=2-10, to=1-10]
		\arrow["{c'}", shift left=1, from=2-10, to=2-11]
		\arrow["{d'}", shift left=1, from=2-11, to=2-10]
		\end{tikzcd}
	\end{center}

  Where $Z_A=\emptyset$ and $Z_B=Z_{new}=\{ a'b',c'd' \}$. It follows from a direct calculation that $$\Im(\delta_A^0)=\langle a\pll a,b\pll b,c\pll c,d\pll d \rangle = \Ker(\delta_A^1).$$ Hence $\HH^1(A)=0$. Similarly we have $$\Im(\delta_B^0)=\langle a'\pll a',b'\pll b',d'\pll d'-c'\pll c',a'\pll a'd'c'-b'\pll d'c'b' \rangle,$$ $$\Ker(\delta_B^1)=\langle a'\pll a',b'\pll b',c'\pll c',d'\pll d',a'\pll a'd'c'-b'\pll d'c'b' \rangle,$$ hence $\HH^1(B)\simeq \langle c'\pll c' \rangle$. Using the notation in Theorem \ref{Lie-strucutre2-of-hh1}, we get the ideal $\mathcal{I} \simeq \langle \varphi_1(\delta_A^0(e_1\pll e_1)) \rangle$ $=\langle c'\pll c'-b'\pll b' \rangle$ and $\HH^1(A)\simeq \HH^1(B)/\mathcal{I}$. It is clear that $\mathcal{AP}^5_1=\{ (a,adc),(b,dcb) \}$, therefore $k(\mathcal{AP}_1^5)=\langle a'\pll a'd'c',b'\pll d'c'b' \rangle$ and
  \[V_{\mathcal{AP}_1^5}=k(\mathcal{AP}_1^5)\cap \Ker(\delta_B^1)=\langle a'\pll a'd'c'-b'\pll d'c'b' \rangle.\]
\end{Ex1}

The following example shows that the difference between the dimensions of $\HH^1(A)$ and $\HH^1(B)$ can be arbitrarily large under gluing two idempotents of $A$.

\begin{Ex1}
\label{diffblockarblarg}
Let $A$ be given by two blocks $A_1$ and $A_2$ such that $A_1$ and $A_2$ are radical square zero local algebras having $m$-loops and $n$-loops respectively. If we exclude the case that $m=1$ and $n=1$ in $\chr(k)=2$ (for this case, see Example \ref{Diffblockinc}), then the dimension of $\HH^1(A)$ is the sum of the dimensions of $\HH^1(A_1)\simeq \mathfrak{gl}_{m}(k)$ and $\HH^1(A_2)\simeq \mathfrak{gl}_{n}(k)$, that is, $m^2+n^2$. Let $B$ be obtained by gluing the units of $A_1$ and $A_2$. Then $\HH^1(B)\simeq \mathfrak{gl}_{m+n}(k)$ and consequently has dimension $(m+n)^2$.
\end{Ex1}

We use the following example to show a particular case of Corollary \ref{product-decomp-source-sink}. 

\begin{Ex1}\label{eg-I-summand-of-HH}
Suppose $\chr(k)=0$. Let $B$ be obtained from $A$ by gluing $e_1$ and $e_4$:
    \begin{center}
    	\begin{tikzcd}
    	&            & e_2\bullet \arrow[rd, "\beta", shift left]                                                                    &             &  &      &            &                                                                                                                                                                    &                                                \\
    	Q_A:\hspace{-2em} & e_3\bullet & e_1\bullet \arrow[r, "\gamma"'] \arrow[u, "\alpha_1"] \arrow[u, "\alpha_2"', shift right=3] \arrow[l, "\eta"] & \bullet e_4 &  & Q_B:\hspace{-2em} & f_3\bullet &  f_1 \bullet \arrow["\gamma'"', loop, distance=2em, in=125, out=55] \arrow[r, "\alpha'_1", shift left=4] \arrow[r, "\alpha'_2" description] \arrow[l, "\eta'"] & \bullet f_2 \arrow[l, "\beta'", shift left=4]
    	\end{tikzcd}
    \end{center}
Where $Z_A=\{\beta\alpha_1\}$ and $Z_{new}=\{(\gamma')^2,\alpha'_i\gamma',\alpha'_i\beta',\gamma'\beta',\eta'\gamma',\eta'\beta'\ |\ i=1,2 \}$. From a straightforward computation we have $$\Im(\delta_A^0)=\langle \alpha_1\pll \alpha_1+\alpha_2\pll \alpha_2+\gamma\pll \gamma,\beta\pll \beta+\gamma\pll \gamma,\eta\pll \eta \rangle,$$ $$\Im(\delta_B^0)=\langle \alpha'_1\pll \alpha'_1+\alpha'_2\pll \alpha'_2-\beta'\pll \beta',\eta'\pll \eta' \rangle,$$ $$\Ker(\delta_A^1)=\langle \alpha_2\pll \alpha_1,\alpha_1\pll \alpha_1,\alpha_2\pll \alpha_2,\beta\pll \beta,\gamma\pll \gamma,\gamma\pll \beta\alpha_2,\eta\pll \eta \rangle.$$ Since we glue a source and a sink, Corollary \ref{Kernel-delta-one-source-and-sink} shows that $\Ker(\delta_B^1)\simeq\Ker(\delta_A^1)$. As a consequence, $$\HH^1(A)\simeq \langle \alpha_2\pll \alpha_1,\alpha_1\pll \alpha_1,\alpha_2\pll \alpha_2,\gamma\pll \beta\alpha_2\rangle,$$ $$\HH^1(B)\simeq \langle \alpha'_2\pll \alpha'_1,\alpha'_1\pll \alpha'_1,\alpha'_2\pll \alpha'_2,\gamma'\pll \gamma',\gamma'\pll \beta'\alpha'_2\rangle.$$ Using the notation in Theorem \ref{Lie-strucutre2-of-hh1}, we have $\HH^1(A)\simeq \HH^1(B)/\mathcal{I}$ with the 1-dimensional Lie ideal $\mathcal{I}$ of $\HH^1(B)$ given by \[\mathcal{I}=\langle \varphi_1(\delta_A^0(e_1\pll e_1)) \rangle=\langle \alpha'_1\pll \alpha'_1+\alpha'_2\pll \alpha'_2+\gamma'\pll \gamma'+\eta'\pll \eta' \rangle.\]

Note  that in this case $\Delta$ in Definition \ref{def-Delta} is  equal to $\{[\alpha],[\gamma] \}$, where $[\alpha]=\{\alpha_1,\alpha_2 \}$ and $[\gamma]=\{\gamma \}$. Since $\eta'\pll \eta'\in \Im(\delta_B^0)$, we can rewrite the generator of $\mathcal{I}$ as follows \[ \mathcal{I}=\langle \varphi_1(\mathcal{I}_{[\alpha]}+\mathcal{I}_{[\gamma]})\rangle=\langle \alpha'_1\pll \alpha'_1+\alpha'_2\pll \alpha'_2+\gamma'\pll \gamma' \rangle.\] By definition we have $L_0^{[\alpha']}=\langle \alpha'_2\pll \alpha'_1,\alpha'_1\pll \alpha'_1,\alpha'_2\pll \alpha'_2 \rangle$ and $L_0^{[\gamma']}=\langle \gamma'\pll \gamma' \rangle$, hence
 $$L_0=L_0^{[\alpha']}\oplus L_0^{[\gamma']}=\langle \alpha'_2\pll \alpha'_1,\alpha'_1\pll \alpha'_1,\alpha'_2\pll \alpha'_2 \rangle \oplus \langle \gamma'\pll \gamma'\rangle$$
$$=\langle \alpha'_2\pll \alpha'_1,\alpha'_1\pll \alpha'_1,\alpha'_2\pll \alpha'_2 \rangle \oplus \langle \alpha'_1\pll \alpha'_1+\alpha'_2\pll \alpha'_2+\gamma'\pll \gamma'\rangle= L_0^{[\alpha']} \oplus \mathcal{I}$$ as Lie algebras. Since $L_1=\langle \gamma'\pll \beta'\alpha'_2 \rangle$, we have the following Lie algebras isomorphism as shown in Corollary \ref{product-decomp-source-sink}, \[\HH^1(B)=L_0\oplus L_1=(L_0^{[\alpha']} \oplus \mathcal{I})\oplus L_1=(L_0^{[\alpha']} \oplus L_1)\oplus \mathcal{I}\simeq \HH^1(A)\oplus \mathcal{I}\simeq \HH^1(A)\oplus k.\]
\end{Ex1}

\medskip
{\bf Acknowledgements and Funding.} The authors would like to thank the referee for the very helpful comments and suggestions, which helped improve and simplify the exposition of this article. This first author was partially supported by NSFC (No. 12031014). The second author has been partially supported by the project PRIN 2017 - Categories, Algebras: Ring-Theoretical and Homological Approaches and  by the project REDCOM: Reducing complexity in algebra, logic,
combinatorics, financed by the programme Ricerca Scientifica di Eccellenza 2018 of the Fondazione Cariverona.  In addition, the second author has been funded by the INDAM- GNSAGA Project (CUP E53C24001950001) and also by the Department of Mathematics ‘Tullio Levi-Civita’ of the University of Padova through its SID 2025 - Investimento Strategico di Dipartimento, within the project `Homotopical and relative methods in Hochschild cohomology'. The second author participates in the INdAM group GNSAGA and he is  grateful to Beijing Normal University  for its hospitality.

\end{document}